\documentclass[12pt]{article}
\usepackage{xcolor}
\usepackage{graphicx} % Required for inserting images
\usepackage{amsfonts}
\usepackage{mathtools}
\usepackage{amsmath}
\usepackage{amssymb}
\usepackage{amsthm}
\usepackage{multicol}
\usepackage{subcaption}
\usepackage{wrapfig}
\usepackage[hidelinks]{hyperref}

\usepackage{setspace}
\usepackage{geometry}
\usepackage{amsmath}
\usepackage{amssymb}
\usepackage{amsthm}
\usepackage{biblatex}

\newcommand{\R}{\mathbb{R}}
\newcommand{\pt}{\partial_t}

\newcommand{\dx}{\nabla_\textbf{x}}
\newcommand{\dv}{\nabla_\textbf{v}}
\newcommand{\p}{\partial}
\newcommand{\E}{\textbf{E}}

\newcommand{\J}{\textbf{J}}
\newcommand{\xx}{\textbf{x}}
\newcommand{\vv}{\textbf{v}}
\newcommand{\uu}{\textbf{u}}

\newcommand{\intI}{\int_{I_{ij}}}
\newcommand{\fh}{f_h}
\newcommand{\psix}{\psi_x}

\newcommand{\intv}{\int^{v_{j+\frac{1}{2}}}_{v_{j-\frac{1}{2}}}}
\newcommand{\intx}{\int^{x_{i+\frac{1}{2}}}_{x_{i-\frac{1}{2}}}}
\newcommand{\fhvm}{\hat{f}_{x,j-\frac{1}{2}}}
\newcommand{\fhvp}{\hat{f}_{x,j+\frac{1}{2}}}
\newcommand{\fhxm}{\hat{f}_{i-\frac{1}{2},v}}
\newcommand{\fhxp}{\hat{f}_{i+\frac{1}{2},v}}
\newcommand{\phvp}{\hat{G}_{x,j+\frac{1}{2}}}
\newcommand{\phvm}{\hat{G}_{x,j-\frac{1}{2}}}

\makeatletter
\let\@fnsymbol\@arabic
\makeatother

\title{A Local Macroscopic Conservative Low-Rank Discontinuous Galerkin Method for the Vlasov-Poisson Equation with Dougherty-Fokker-Planck Collisions}
\author{Austin Nelson\thanks{Department of Mathematics and Statistics, Texas Tech University, Lubbock, TX, 70409. Email: austinne@ttu.edu.} , Wei Guo\thanks{Department of Mathematics and Statistics, Texas Tech University, Lubbock, TX, 70409. Email: weimath.guo@ttu.edu.} , Pierson Guthrey\thanks{Lawrence Livermore National Laboratory, Livermore, CA, 74550. Email: guthrey1@llnl.gov}}
\date{}

\begin{document}

\nocite{*}
\maketitle
\allowdisplaybreaks

\begin{abstract}
    In this paper, we construct a low-rank, structure preserving discontinuous Galerkin (DG) method to simulate the Vlasov-Poisson (VP) system coupled with the Dougherty Fokker-Planck (DFP) collision operator. When Coulomb collisions occur in dense or weakly-collisional plasmas, electrons get pushed to a low-rank steady state. In many cases, the plasma arrives to this steady state quickly, meaning that for most of the runtime, the plasma consists mainly of numerical low-rank structures. Our new low-rank scheme is constructed to exploit these numerical low-rank structures to greatly reduce the needed storage complexity of simulations for the VP-DFP system. It is constructed as an extension of the previously established Local Macroscopic Conservative (LoMaC) method by incorporating Coulomb collisions into the system. The LoMaC property ensures local conservation of macroscopic mass, momentum, and energy at the discrete level. Details of the new method are discussed in this paper. Numerical experiments are performed to show the efficacy of the method.
\end{abstract}

\begin{flushleft}\textbf{Keywords:} Vlasov-Poisson equation, Dougherty Fokker-Planck operator, local conservation, local Discontinuous Galerkin method.\end{flushleft}

\stepcounter{footnote}\stepcounter{footnote}\stepcounter{footnote}\stepcounter{footnote}\footnotetext{This work was performed under the auspices of the U.S. Department of Energy by
Lawrence Livermore National Laboratory under contract DE-AC52-07NA27344. LLNL-JRNL-2021257-DRAFT.}

\section{Introduction}
%%%%%%%%%%%%%%%%%%%%%%%%%%%%%%%%%%%%%%%%%%%%%%%%%%%%%%%%%%%%%%%%%%%%%%%%%%%%%%
%%% Introduction
%%%%%%%%%%%%%%%%%%%%%%%%%%%%%%%%%%%%%%%%%%%%%%%%%%%%%%%%%%%%%%%%%%%%%%%%%%%%%%
% Intro Paragraph 1
 The Vlasov–Poisson (VP) system is a fundamental kinetic model for describing the evolution of the probability distribution function of collisionless particles under self-consistent electrostatic fields in phase space, with its roots in statistical physics and Hamiltonian transport theory. However, in many modern applications, particularly those involving dense or weakly collisional plasmas, collisions play a non-negligible role and must be incorporated to obtain physically meaningful predictions. The Dougherty Fokker–Planck (DFP) model \cite{dougherty1964model} is widely utilized, offering a computationally tractable alternative to the full Landau operator, while retaining essential physical properties, including the conservation of mass, momentum, and energy. Despite significant research efforts devoted to developing effective numerical methods for the VP–DFP system, substantial challenges remain due to the inherent computational complexity of the system. These challenges include, but are not limited to, the high dimensionality of phase space, stiffness induced by the collision operator, multiscale dynamics, and the need to preserve critical physical structures.

%%%%%%%%%%%%%%%%%%%%%%%%%%%%%%%%%%%%%
% Intro Paragraph 2

Recently, low-rank approaches have emerged as effective tools for reducing storage and computational costs in kinetic simulations by exploiting the intrinsic low-rank structure of the solution through advanced matrix and tensor decomposition techniques. One such class is the dynamical low-rank (DLR) approach, which has been extensively used for kinetic simulations. In \cite{einkemmer2021mass}, this technique was employed to construct a low-rank scheme to model the Vlasov equation while conserving mass, momentum, and energy. Most recently, \cite{coughlin2024robust} used this approach to construct the first low-rank energy-conserving method for the VP-DFP system. Another notable class of low-rank methods is the step-and-truncate (SAT) approach. \cite{guo2024LoMaC,guo2024LoMaCDG} have employed this technique to construct a local macroscopic conservative (LoMaC) low-rank method to simulate the Vlasov-Poisson equation. However, there is currently no SAT scheme that uses this approach to model the VP-DFP system. In this paper, we fill this gap by constructing a new LoMaC scheme to effectively find low-rank solutions to the VP-DFP system while locally conserving mass, momentum, and energy at the discrete level.

%%%%%%%%%%%%%%%%%%%%%%%%%%%%%%%%%%%%%
% Intro Paragraph 3

The VP-DFP system is given by
\begin{gather*}
    \frac{\p f}{\p t} + \vv \cdot \dx f + \E(\xx,t) \cdot \dv f = \nu \dv \cdot (T \dv f + (\vv - \uu)f), \tag{1.1}\label{eq:1.1}\\
    \E(\xx,t) = -\dx \phi(\xx,t), \tag{1.2} \label{eq:1.2}\\ 
    -\Delta_\xx \phi(\xx,t) = \rho(\xx,t) - \rho_0,\tag{1.3}\label{eq:1.3}
\end{gather*}
which describes the dynamics of the probability distribution function $f(\xx,\vv,t)$ of electrons in a collisional plasma. $\E$ is the electric field and $\phi$ is the self-consistent electrostatic potential determined by Poisson's equation. $f$ couples to the long range fields via the charge density $\rho(\xx,t) = \int_{\Omega_v}f(\xx,\vv,t)d\vv,$ where we take the limit of uniformly distributed infinitely massive ions in the background. The effects of collision are determined by the collision frequency $\nu$, and the Dougherty Fokker-Planck collision operator is given by
\begin{gather*}
    C[f] = \dv \cdot (T \dv f + (\vv - \uu)f), \tag{1.4}\label{eq:1.4}\\
    T = \frac{1}{\rho d} \int_{\Omega_v} 
    |\vv - \uu|^2 fd\vv, \tag{1.5} \label{eq:1.5}\\
    \uu = \frac{\J}{\rho}, \tag{1.6}\label{eq:1.6}
\end{gather*}
which models the effects of Coulomb collisions on the particles \cite{dougherty1964model}. Here, $T$ is the temperature of the system, $\uu$ is the average velocity, d is the dimension of the velocity domain, and $\J$ is the current density given by $\J(\xx,t) = \int_{\Omega_v}f(\xx,\vv,t)\vv d\vv.$
%%%%%%%%%%%%%%%%%%%%%%%%%%%%%%%%%%%%%
% Intro Paragraph 4

When the collision frequency is zero ($\nu = 0$), the VP-DFP system becomes the collisionless VP system. In this setting, many problems, such as the strong Landau or two-stream problem, see the development of thin filamentation structures that cause rank increase. As shown in \cite{guo2024LoMaCDG}, although low-rank methods can reduce computational complexity, they cannot fully prevent the rank growth induced by these filamentation structures. However, a different story occurs for sufficient values of $\nu$ (e.g. $\nu \geq 0.1$). In this setting, the collisional effects of the Dougherty Fokker-Planck operator set in and push the electrons towards a steady state in the form of a low-rank Maxwellian. In most cases, this steady state is reached rapidly, so the computational runtime is largely dominated by the dynamics of these low-rank structures. \cite{ye2024energy} shows this well for the two-stream problem and bump-on-tail problem. This phenomenon exemplifies the need for low-rank methods such as \cite{coughlin2024robust} to avoid unnecessary computational complexity from redundant storage of numerical low-rank structures in the VP-DFP system. As demonstrated in this paper, our novel low-rank method specifically tackles this problem while maintaining the physical accuracy of the system.
%%%%%%%%%%%%%%%%%%%%%%%%%%%%%%%%%%%%%
% Intro Paragraph 5

Our method can be seen as an extension of the previously established LoMaC method \cite{guo2024LoMaCDG} with the incorporation of a low-rank discretization of the DFP operator. It consists of two main components: a discontinuous Galerkin (DG) discretization, and a conservative decomposition. The use of DG discretizations is due to their efficiency, high-order accuracy, and stability properties. \cite{guo2024LoMaCDG,cai2021high,cheng2014energy,de2012high} exemplify this well, especially \cite{ye2024energy} whose methods show the full capabilities of DG schemes on the VP-DFP system. For our method, we take inspiration from \cite{ye2024energy} to develop our own DG discretization that works in the low-rank format. The use of a conservative decomposition is necessary to maintain the physical accuracy of the system by ensuring local conservation of mass, momentum, and energy on the discrete level. Let 
\begin{align*}
    \text{charge density: }& \hspace{4mm}\rho(\xx,t) = \int_{\Omega_\vv} f(\xx,\vv,t)dv, \tag{1.6}\label{eq:1.6}\\
    \text{current density: }& \hspace{4mm} \textbf{J}(\xx,t) = \int_{\Omega_\vv} f(\xx,\vv,t)\vv d\vv,\tag{1.7}\label{eq:1.7}\\
    \text{kinetic energy density: }& \hspace{4mm} \kappa(\xx,t) = \int_{\Omega_\vv} \frac{1}{2}|\vv|^2f(\xx,\vv,t)d\vv,\tag{1.8}\label{eq:1.8}\\
    \text{energy density: }& \hspace{4mm} e(\xx,t) = \kappa(\xx,t) + \frac{1}{2}\textbf{E}(\xx)^2.\tag{1.9}\label{eq:1.9}
\end{align*}
We can derive the conservation laws by taking the first few moments of the VP-DFP system
\begin{align*}
        \pt \rho + \dx \cdot \J &= 0,\\
        \pt \J + \dx \cdot \sigma &= \rho \E, \tag{1.10} \label{eq:1.10}\\
        \pt e + \dx \cdot Q &= 0,
\end{align*}
where $\textbf{$ \sigma $}(t,\xx) = \int_{\Omega_\vv}(\vv \otimes \vv)f(\xx,\vv,t)d\vv$ and $\textbf{Q}(\xx,t) = \frac{1}{2}\int_{\Omega_\vv} \vv|\vv|^2 f(\xx,\vv,t)d\vv.$ Similar to the work in
\cite{guo2024LoMaC,guo2024LoMaCDG}, our method evolves the mass, momentum, and total energy from the conservation laws at each time step. Then, we perform a conservative decomposition where one part will hold the information of these quantities, and the other will be truncated to exploit low-rank structures. This process will allow our method to stay low-rank while ensuring physical accuracy on the discrete level.
%%%%%%%%%%%%%%%%%%%%%%%%%%%%%%%%%%%%%
% Intro Paragraph 6

The paper is organized as follows. In Section 2, we introduce the DG framework along with the nodal DG discretization. In Section 3, we discuss the low rank framework with a tensor product of nodal DG meshes, the weighted inner product spaces, and the corresponding macroscopic conservative projection and weighted SVD truncation. We outline the method in Section 4 and provide remarks about extensions to higher dimensions. In Section 5, we benchmark our method with several numerical examples. We conclude in Section 6.

%%%%%%%%%%%%%%%%%%%%%%%%%%%%%%%%%%%%%%%%%%%%%%%%%%%%%%%%%%%%%%%%%%%%%%%%%%%%%%
%%% DG Section
%%%%%%%%%%%%%%%%%%%%%%%%%%%%%%%%%%%%%%%%%%%%%%%%%%%%%%%%%%%%%%%%%%%%%%%%%%%%%%
\section{A nodal DG framework for the collisional Vlasov dynamics}

In this section we construct the discontinuous Galerkin (DG) discretization of the system. For simplicity of illustrating the main idea of the method, we only discuss the 1D1V case in the following section.
\subsection{DG discretization with nodal Lagrangian basis functions.}
%%%%%%%%%%%%%%%%%%%%%%%%%%%%%%%%%%%%%%%%%%%%%%%%%%%%%%%%%%%%%%%%%%%%%%%%%%%%%%
%%% DG Formulation
%%%%%%%%%%%%%%%%%%%%%%%%%%%%%%%%%%%%%%%%%%%%%%%%%%%%%%%%%%%%%%%%%%%%%%%%%%%%%%

We begin with a tensor product Cartesian partition on a truncated 1D1V domain $\Omega = [x_{\text{min}},x_{\text{max}}] \times [-v_{\text{max}},v_{\text{max}}]$, denoted by $\Omega_{h}$ with
\begin{align*}
    x_{\text{min}} &= x_{\frac{1}{2}} < x_{\frac{3}{2}} < \cdots < x_{N_x + \frac{1}{2}} = x_{\text{max}},\\
    -v_{\text{max}} &= v_{\frac{1}{2}} < v_{\frac{3}{2}} < \cdots < v_{N_{v} + \frac{1}{2}} = v_{\text{max}}.
\end{align*}
Define an element as $I_{ij} = [x_{i - \frac{1}{2}}, x_{i + \frac{1}{2}}] \times [v_{j - \frac{1}{2}}, v_{j + \frac{1}{2}}] \in \Omega_h.$ Each element has size $h_{x,i}h_{v,j}$ and a center of $(x_i,v_i) = (\frac{1}{2}[x_{i-\frac{1}{2}} + x_{i+\frac{1}{2}}], \frac{1}{2}[v_{j-\frac{1}{2}} + v_{j+\frac{1}{2}}]).$ We define the finite dimensional discrete space of piecewise polynomials as 
\begin{equation*}
    Q^k_{h} = \{ p(x,v) \in L^2(\Omega) : p|_{I_{ij}}\in Q^k(I_{ij}), \forall I_{ij} \in \Omega_h \}, \tag{2.1}\label{eq:2.1}
\end{equation*}
where the local space $Q^k(I)$ consists of polynomials with terms of the form $x^mv^n$ with max($m,n) \leq k$ on $I \in \Omega_h.$ To distinguish the left and right limits of a function $p \in Q^k_h$ at $(x_{i + \frac{1}{2}},v),$ we let \textcolor{blue}{$p^{\pm}_{i+\frac{1}{2},v} = \text{lim}_{\delta \to \pm 0} p(x_{i+\frac{1}{2}} + \delta,v).$}
%%%%%%%%%%%%%%%%%%%%%%%%%%%%%%%%%%%%%

Similar to the work in \cite{ye2024energy}, we start by rewriting the original system containing second-order derivatives into an equivalent first-order system.
\begin{gather*}
    \pt f + v \cdot \p_x f + E \cdot \p_v f = \nu(\p_v \cdot (TG) + \p_v \cdot((v - u)f)) \tag{2.2} \label{eq:2.2}\\
    G = \p_v f \tag{2.3} \label{eq:2.3}
\end{gather*}
where $G$ is an auxiliary variable. Using this new system, we may construct our semi-discrete local DG method as the following: find $f_h(\cdot,\cdot,t),G_h(\cdot,\cdot,t) \in Q^k_h$, such that $\forall \psi \in Q^k_h$ and $\forall I_{ij} \in \Omega_h,$

{\small
\begin{align}
        \intI \pt \fh \psi dxdv &= \intI v\fh \psix dxdv - \intv v\left( \fhxp \psi^-_{i+\frac{1}{2},v} - \fhxm \psi^+_{i-\frac{1}{2},v} \right)dv \tag{2.4a} \label{eq:2.4a}\\
        &+ \intI Ef_h \psi_v dxdv - \intx E(x)\left( \fhvp \psi^-_{x,j+\frac{1}{2}} - \fhvm \psi^+_{x,j-\frac{1}{2}} \right)dx \tag{2.4b} \label{eq:2.4b}\\
        &- \nu\intI TG_h \psi_v dxdv + \nu\intx T\left( \phvp \psi^-_{x,j+\frac{1}{2}} - \phvm \psi^+_{x,j-\frac{1}{2}} \right) dx \tag{2.4c} \label{eq:2.4c}\\
        &- \nu\intI vf_h\psi_v dxdv + \nu\intx \left( v_{j+\frac{1}{2}} \tilde{f}_{x,j+\frac{1}{2}} \psi^-_{x,j+\frac{1}{2}} - v_{j - \frac{1}{2}} \tilde{f}_{x,j-\frac{1}{2}} \psi^+_{x,j-\frac{1}{2}} \right) dx \tag{2.4d} \label{eq:2.4d}\\
        &+ \nu\intI uf_h\psi_v dxdv - \nu\intx \left( u_{i+\frac{1}{2}} \tilde{\tilde{f}}_{x,j+\frac{1}{2}} \psi^-_{x,j+\frac{1}{2}} - u_{i-\frac{1}{2}}\tilde{\tilde{f}}_{x,j-\frac{1}{2}} \psi^+_{x,j-\frac{1}{2}}\right) dx, \tag{2.4e} \label{eq:2.4e}\\
     G_h &= \int_{I_{ij}} f_h \psi_v dxdv - \int^{v_{j+1}}_{v_{j-1}}\left(\hat{\hat{f}}_{x,j+{\frac{1}{2}}}\psi^{-}_{x,j+\frac{1}{2}} - \hat{\hat{f}}_{x,j-\frac{1}{2}}\psi^{+}_{x,j-\frac{1}{2}}\right)dx. \tag{2.4f} \label{eq:2.4f} 
    \end{align}
    }
    
\noindent Here, $\hat{f},\hat{\hat{f}},\tilde{f},\tilde{\tilde{f}}$, and $\hat{G}$ are numerical fluxes that will be defined momentarily. For now, it can be verified that this system is $L^2$ stable for sufficient $f_h$ and $G_h$, which we denote in the following proposition.
\begin{flushleft}
    \textbf{Proposition 2.1:} The system defined by \eqref{eq:2.4a}-\eqref{eq:2.4f} is $L^2$ stable for sufficient $f_h,G_h \in Q^k_h.$
\end{flushleft}
%%%%%%%%%%%%%%%%%%%%%%%%%%%%%%%%%%%%%

To define our discrete scheme, we use a nodal Lagrangian basis to represent our functions in $Q^k_n$ and approximate the integrals by numerical quadrature. Define a reference cell $I = [-\frac{1}{2},\frac{1}{2}] \times [-\frac{1}{2},\frac{1}{2}]$ and a tensor product of Gaussian quadrature points in each direction $\{\xi_{ig},\eta_{jg}\}^k_{ig,jg=0}$. Let $\{\omega_l\}^k_{l=0}$ be the corresponding quadrature weights on the reference element. The local nodal Lagrangian basis on the reference cell is $\{L_{ig,jg}(\xi,\eta)\}^k_{ig,jg=0} \in Q^k_h$ where 
\begin{equation*}
    L_{ig,jg}(\xi_{ig'},\eta_{jg'}) = L_{ig}(\xi_{ig'})L_{jg}(\eta_{jg'}) =  \delta_{ig,ig'}\delta_{jg,jg'}, \hspace{1cm} ig,ig',jg,jg' = 0,...,k. \tag{2.5} \label{eq:2.5}
\end{equation*}
Here, $\delta_{\cdot,\cdot'}$ is the Kronecker delta function. We can use this basis on each computational cell $I_{ij}$ by taking the transformations $\xi = \frac{x - x_i}{h_{x,i}}$ and $\eta = \frac{v - v_i}{h_{v,i}}.$
%%%%%%%%%%%%%%%%%%%%%%%%%%%%%%%%%%%%%

 Using the nodal basis functions, we can equivalently replace the test functions of \eqref{eq:2.4a}-\eqref{eq:2.4f} with $L_{ig',jg'}(\xi,\eta),ig',jg'=0,...,k.$ And the integrals can be approximated using the quadrature points $\{\xi_{ig},\eta_{jg}\}^k_{ig,jg=0}$ and quadrature weights $\{\omega_l\}^k_{l=0}.$ We then look for the DG solution expressed in the form of $f_{h,i,j}(x,v,t) = \sum^k_{ig,jg=0} f^{ig,jg}_{h,i,j}L_{ig,jg}(\xi(x),\eta(v))$, with its nodal values satisfying the following equations:
 
{\small
\begin{align*}
    &h_{x,i}h_{v,j}\omega_{ig}\omega_{jg}\left(\frac{d}{dt}f^{ig,jg}_{h,i,j}(t)\right) \tag{2.6} \label{eq:2.6}\\
    =&h_{x,i}h_{v,j}\omega_{j,g}v_{j,jg}\sum_{ig''}\omega_{ig''}\left( \frac{d}{dx}L_{ig}(\xi_{ig''})f^{ig'',jg}_{h,i,j}(t) \right) - h_{v,j}\omega_{jg}v_{j,jg}\left( \hat{f}_{i+\frac{1}{2},jg} L_{ig}\left(\frac{1}{2}\right) - \hat{f}_{i - \frac{1}{2},jg} L_{ig}\left(-\frac{1}{2}\right) \right)\\
    +&h_{x,i}h_{v,j}\omega_{ig}E_{i,ig}\sum_{jg''}\omega_{jg''}\left( \frac{d}{dv}L_{jg}(\eta_{jg''})f^{ig,jg''}_{h,i,j}(t)\right) - h_{x,i}\omega_{ig}E_{i,ig}\left( \hat{f}_{ig,j+\frac{1}{2}}L_{jg}\left(\frac{1}{2}\right) - \hat{f}_{ig,j-\frac{1}{2}}L_{jg}\left(-\frac{1}{2}\right)\right)\\
    -&\nu h_{x,i}h_{v,j}\omega_{ig}T_{i,ig}\sum_{jg''}\omega_{jg''}\left( \frac{d}{dv}L_{jg}(\eta_{jg''})G^{ig,jg''}_{h,i,j}(t) \right) - \nu h_{x,i}\omega_{ig}T_{i,ig}\left( \hat{G}_{ig,j+\frac{1}{2}}L_{jg}\left(\frac{1}{2}\right) - \hat{G}_{ig,j-\frac{1}{2}}L_{jg}\left(-\frac{1}{2}\right)\right)\\
    -&\nu h_{x,i}h_{v,j}\omega_{ig}\sum_{jg''}\omega_{jg''}\left( \frac{d}{dv}L_{jg}(\eta_{jg''})v_{j,jg''}f^{ig,jg''}_{h,i,j}(t)\right) - \nu h_{x,i}\omega_{ig}\left( v_{j+\frac{1}{2}}\tilde{f}_{ig,j+\frac{1}{2}}L_{jg}\left(\frac{1}{2}\right) - v_{j-\frac{1}{2}}\tilde{f}_{ig,j-\frac{1}{2}}L_{jg}\left(-\frac{1}{2}\right) \right)\\
    +&\nu h_{x,i}h_{v,j}\omega_{ig}u_{i,ig}\sum_{jg''}\omega_{jg''}\left( \frac{d}{dv}L_{jg}(\eta_{jg''})f^{ig,jg''}_{h,i,j}(t)\right) - \nu h_{x,i}\omega_{ig}u_{i,ig}\left(\tilde{\tilde{f}}_{ig,j+\frac{1}{2}}L_{jg}\left(\frac{1}{2}\right) - \tilde{\tilde{f}}_{ig,j-\frac{1}{2}}L_{jg}\left(-\frac{1}{2}\right)\right).
\end{align*}
}
A similar (and shorter) construction is also made for the nodal values of $G_{h,i,j}$ in \eqref{eq:2.4f}.
%%%%%%%%%%%%%%%%%%%%%%%%%%%%%%%%%%%%%

All that is left now is to compute the numerical fluxes. We compute the numerical flux $\hat{G}$ using an alternating flux: first, G is computed using an upwind flux, then $\hat{G}$ is obtained using a downwind flux. As for $\hat{f},\tilde{f},\tilde{\tilde{f}},$ we take them as monotone upwind fluxes. Let $v^+ = $ max($v,0$), $v^- = $ min($v,0$), $E^+ = $ max($E,0$), $E^- = $ min($E,0$), $u^+ = $ max($u,0$), and $u^- = $ min($u,0$). We can then take $\hat{f},\tilde{f},$ and$\tilde{\tilde{f}}$ to be upwind fluxes based on $v$, $E$, and $u$, respectively. After computing these numerical fluxes and performing minor simplifications of the coefficients, the scheme \eqref{eq:2.6} can be written as follows.

{\small
\begin{align*}
    &\pt f^{ig,jg}_{h,i,j}(t)  \tag{2.7} \label{eq:2.7}\\
    =&\frac{v^+_{j,jg}}{\omega_{ig}h_{x,i}}\left( \sum_{ig''} \omega_{ig''}\frac{dL_{ig}}{d\xi}(\xi_{ig'})f^{ig''.jg}_{h,i,j} - f^{ig'',jg}_{h,i,j}L_{ig''}\left(\frac{1}{2}\right)L_{ig}\left(\frac{1}{2}\right) + f^{ig'',jg}_{h,i-1,j}L_{ig''}\left(\frac{1}{2}\right)L_{ig}\left(-\frac{1}{2}\right)\right)\\
    +&\frac{v^-_{j,jg}}{\omega_{ig}h_{x,i}}\left( \sum_{ig''} \omega_{ig''}\frac{dL_{ig}}{d\xi}(\xi_{ig'})f^{ig''.jg}_{h,i,j} - f^{ig'',jg}_{h,i+1,j}L_{ig''}\left(-\frac{1}{2}\right)L_{ig}\left(\frac{1}{2}\right) + f^{ig'',jg}_{h,i,j}L_{ig''}\left(-\frac{1}{2}\right)L_{ig}\left(-\frac{1}{2}\right)\right)\\
    +&\frac{E^+_{i,ig}}{\omega_{jg}h_{v,j}}\left( \sum_{jg''}\omega_{jg''}\frac{dL_{jg}}{d\eta}(\eta_{jg''})f^{ig,jg''}_{h,i,j} - f^{ig,jg''}_{h,i,j}L_{jg''}\left(\frac{1}{2}\right)L_{jg}\left(\frac{1}{2}\right) + f^{ig,jg''}_{h,i,j-1}L_{jg''}\left(\frac{1}{2}\right)L_{jg}\left(-\frac{1}{2}\right)\right)\\
    +&\frac{E^-_{i,ig}}{\omega_{jg}h_{v,j}}\left( \sum_{jg''}\omega_{jg''}\frac{dL_{jg}}{d\eta}(\eta_{jg''})f^{ig,jg''}_{h,i,j} - f^{ig,jg''}_{h,i,j+1}L_{jg''}\left(-\frac{1}{2}\right)L_{jg}\left(\frac{1}{2}\right) + f^{ig,jg''}_{h,i,j}L_{jg''}\left(-\frac{1}{2}\right)L_{jg}\left(-\frac{1}{2}\right)\right)\\
    -&\frac{\nu T_{i,ig}}{\omega_{jg}h_{v,j}}\left( \sum_{jg''}\omega_{jg''}\frac{dL_{jg}}{d\eta}(\eta_{jg''})G^{ig,jg''}_{h,i,j} - G^{ig'',jg}_{h,i,j+1}L_{jg''}\left(-\frac{1}{2}\right)L_{jg}\left(\frac{1}{2}\right) + G^{ig,jg''}_{h,i,j}L_{jg''}\left(-\frac{1}{2}\right)L_{jg}\left(-\frac{1}{2}\right)\right)\\
    -&\frac{\nu}{\omega_{jg}h_{v,j}}\left( \sum_{jg''}\omega_{jg''}\frac{dL_{jg}}{d\eta}(\eta_{jg''})v^+_{j,jg''}f^{ig,jg''}_{h,i,j} - v^+_{j,jg''}f^{ig,jg''}_{h,i,j}L_{jg''}(-\frac{1}{2})L_{jg}(\frac{1}{2}) + v^+_{j,jg''}f^{ig,jg''}_{h,i,j-1}L_{jg''}(-\frac{1}{2})L_{jg}(-\frac{1}{2})\right)\\
    -&\frac{\nu}{\omega_{jg}h_{v,j}}\left( \sum_{jg''}\omega_{jg''}\frac{dL_{jg}}{d\eta}(\eta_{jg''})v^-_{j,jg''}f^{ig,jg''}_{h,i,j} - v^-_{j,jg''}f^{ig,jg''}_{h,i,j+1}L_{jg''}(-\frac{1}{2})L_{jg}(\frac{1}{2}) + v^-_{j,jg''}f^{ig,jg''}_{h,i,j}L_{jg''}(-\frac{1}{2})L_{jg}(-\frac{1}{2})\right)\\
    +&\frac{\nu u^+_{i,ig}}{\omega_{jg}h_{v,j}}\left( \sum_{jg''}\omega_{jg''}\frac{dL_{jg}}{d\eta}(\eta_{jg''})f^{ig,jg''}_{h,i,j} - f^{ig,jg''}_{h,i,j}L_{jg''}\left(\frac{1}{2}\right)L_{jg}\left(\frac{1}{2}\right) + f^{ig,jg''}_{h,i,j-1}L_{jg''}\left(\frac{1}{2}\right)L_{jg}\left(-\frac{1}{2}\right)\right)\\
    +&\frac{\nu u^-_{i,ig}}{\omega_{jg}h_{v,j}}\left( \sum_{jg''}\omega_{jg''}\frac{dL_{jg}}{d\eta}(\eta_{jg''})f^{ig,jg''}_{h,i,j} - f^{ig,jg''}_{h,i,j+1}L_{jg''}\left(-\frac{1}{2}\right)L_{jg}\left(\frac{1}{2}\right) + f^{ig,jg''}_{h,i,j}L_{jg''}\left(-\frac{1}{2}\right)L_{jg}\left(-\frac{1}{2}\right)\right)\\
\end{align*}
}
%%%%%%%%%%%%%%%%%%%%%%%%%%%%%%%%%%%%%
Finally, we can re-express our expression \eqref{eq:2.7} into the following simpler system.
\begin{align*}
    \p_t f^{ig,jg}_{h,i,j}(t) &= v^{+}_{j,jg} \cdot D^{+}_{x,i,ig}f^{+,:,jg}_{h,i,j} + v^{-}_{j,jg} \cdot D^{-}_{x,i,ig}f^{-,:,jg}_{h,i,j} \tag{2.8a} \label{eq:2.8a}\\
    &+ E^{+}_{i,ig} \cdot D^{+}_{v,j,jg}f^{+,ig,:}_{h,i,j} + E^{-}_{i,ig} \cdot D^{-}_{v,j,jg}f^{-,ig,:}_{h,i,j} \tag{2.8b} \label{eq:2.8b}\\
    &- \nu T_{i,ig} \cdot D^{-}_{v,j,jg}G^{-,ig,:}_{h,i,j} \tag{2.8c} \label{eq:2.8c}\\
    &- \nu\left(D^{+}_{v,j,jg} (v^{+}_{j,jg} \cdot f^{+,ig,:}_{h,i,j}) - D^{-}_{v,j,jg} (v^{-}_{j,jg} \cdot f^{-,ig,:}_{h,i,j})\right) \tag{2.8d} \label{eq:2.8d}\\
    &+ \nu\left(u^{+}_{i,ig} \cdot D^{+}_{v,j,jg} f^{+,ig,:}_{h,i,j} - u^{-}_{i,ig} \cdot D^{-}_{v,j,jg} f^{+,ig,:}_{h,i,j}\right). \tag{2.8e} \label{eq:2.8e}\\
    G^{ig,:}_{h,i,j} &= D^+_{v,j,jg} f^{+,ig,:}_{h,i,j} \tag{2.8f} \label{eq:2.8f}
\end{align*}
Here, $\textbf{f}$ is given as
\begin{align*}
    \textbf{f}^{+,:,jg}_{h,i,j} &= (f^{0,jg}_{h,i-1,j},...,f^{k,jg}_{h,i-1,j},f^{0,jg}_{h,i,j},...,f^{k,jg}_{h,i,j}),\\
    \textbf{f}^{-,:,jg}_{h,i,j} &= (f^{0,jg}_{h,i,j},...,f^{k,jg}_{h,i,j},f^{0,jg}_{h,i+1,j},...,f^{k,jg}_{h,i+1,j}),\\
    \textbf{f}^{+,ig,:}_{h,i,j} &= (f^{ig,0}_{h,i,j-1},...,f^{ig,k}_{h,i,j-1},f^{ig,0}_{h,i,j},...,f^{ig,k}_{h,i,j}),\\
    \textbf{f}^{-,ig,:}_{h,i,j} &= (f^{ig,0}_{h,i,j},...,f^{ig,k}_{h,i,j},f^{ig,0}_{h,i,j+1},...,f^{ig,k}_{h,i,j+1}).
\end{align*}
$\textbf{G}$ is constructed similarly.

A major advantage to this scheme is that differentiation is no longer an element-by-element process, but rather a dimension-by-dimension one. This allows our method to be more tensor-friendly to tensor decompositions such as the hierarchical Tucker or tensor train decomposition. This is especially beneficial when extending the system to higher dimensions (e.g. 2d2v), since these decompositions can add considerable compression to the required storage of the system. Further remarks of this can be found in section 4.

 We would now like to highlight the primary differences between the construction of our DG scheme and that of \cite{ye2024energy}. Although we follow the initial idea of using a local DG framework on the diffusive term $\p_v \cdot (T \p_v f)$ in \eqref{eq:2.2}, our approach diverges at the semi-discrete level. In particular, the average velocity $u = u(x,t)$ depends on position and time, not on velocity. Constructing a numerical flux for the term $\p_v \cdot ((v-u)f)$ based on $v$ and $u$ together would couple variables from different dimensions of the phase space which may destroy the numerical low-rank structure we desire. To keep our method low-rank, we instead separate this term entirely into the components shown in \eqref{eq:2.4d} and \eqref{eq:2.4e}.

%%%%%%%%%%%%%%%%%%%%%%%%%%%%%%%%%%%%%%%%%%%%%%%%%%%%%%%%%%%%%%%%%%%%%%%%%%%%%%
%%% LoMaC Projection
%%%%%%%%%%%%%%%%%%%%%%%%%%%%%%%%%%%%%%%%%%%%%%%%%%%%%%%%%%%%%%%%%%%%%%%%%%%%%%
\section{A LoMaC low rank tensor approach with DG discretization}

In this section, we discuss the low rank approach and conservative decomposition of the method. There is not much change in the approach compared to our previous work in \cite{guo2024LoMaC}. We will only provide a brief summary of the low rank approach and LoMaC projection step as necessary to our method. We refer to \cite{guo2024LoMaC} for further discussions.

\subsection{A low-rank representation with DG discretization}

The low rank tensor approach \cite{guo2022low} is based on the assumption that the solution at time $t$ has a Schmidt decomposition of the form
\begin{equation*}
    f(x,v,t) = \sum^r_{l=1} \left( C_l(t)U^{(1)}_{l}(x,t)U^{(2)}_{l}(v,t)\right), \tag{3.1} \label{eq:3.1}
\end{equation*}
where $\left\{ U^{(1)}_{l}(x,t) \right\}^r_{l=1}$ and $\left\{ U^{(2)}_l (v,t)\right\}^r_{l=1}$ are sets of time-dependent low rank orthonormal basis in x and v dimensions, respectively, $C_l$ is the coefficient for the basis $U^{(1)}_{l}(x,t)U^{(2)}_{l}(v,t)$, and r is the representation rank. 
%%%%%%%%%%%%%%%%%%%%%%%%%%%%%%%%%%%%%

 Using the DG framework constructed above, we now define the nodal grid points for the DG discretization. These are given as the following tensor product of $(k+1)N_x \times (k+1)N_v$ points from $N_x \times N_v$ computational cells,
\begin{align*}
    x_{grid} : \hspace{3mm} x_{\text{min}} < ... < (x_{i,0} < ... < x_{i,k}) ... < x_{\text{max}}, \tag{3.2} \label{eq:3.2}\\
    v_{grid} : \hspace{3mm} -v_{\text{max}} < ... < (v_{j,0}< ... < v_{j,k}) ... < v_{\text{max}}. \tag{3.3} \label{eq:3.3}
\end{align*}
Here, $\{x_{i,ig}\}^k_{ig=0}$ and $\{v_{j,jg}\}^k_{jg=0}$ are the shifted Gaussian points on the cell $[x_{i-\frac{1}{2}},x_{i+\frac{1}{2}}]$ and $[v_{j - \frac{1}{2}}, v_{j + \frac{1}{2}}]$ respectively. Thus, our DG nodal solutions are organized as $\textbf{f} \in \R^{(k+1)N_x \times (k+1)N_v}$ with each of its component $f^{ig,jg}_{h,i,j}(t)$ being an approximation to point values of the solution on the tensor product of grids \eqref{eq:3.2}-\eqref{eq:3.3}. Connecting this to our low-rank approach, each nodal value $f^{ig,jg}_{h,i,j}$ has the following low rank approximation
\begin{equation*}
    f^{ig,jg}_{h,i,j} = \sum^{r}_{l = 1} C_l U^{(1)}_{l,i,ig}U^{(2)}_{l,j,jg}. \tag{3.4} \label{eq:3.4}
\end{equation*}
Using this low-rank nodal DG construction of the solution, we can construct the discrete macroscopic quantities of $f.$ The discrete macroscopic charge, current, and kinetic energy density $\rho,\J,\kappa \in \R^{(k+1)N_x}$ are computed as
\begin{equation*}
    \begin{pmatrix}
        \rho\\ \J\\ \kappa
    \end{pmatrix}
    = \sum^r_{l=1} C_l \left\langle U^{(2)}_{l}, \begin{pmatrix} \textbf{1}_v\\ \vv \\ \frac{1}{2}\vv^2\end{pmatrix}\right\rangle_v U^{(1)}_{l} \tag{3.5} \label{eq:3.5}
\end{equation*}
where $\langle \textbf{f},\textbf{g}\rangle_v = \sum_{j,jg} f_{j,jg}g_{j,jg}\omega_{v,j,jg},$ $\textbf{f},\textbf{g} \in \R^{(k+1)N_v}.$

\subsection{ A macroscopic conservative decomposition with DG discretization}
The central idea behind the conservative projection is to project the kinetic solution f onto the subspace
\begin{equation*}
    \mathcal{N} = \text{span}\{ \textbf{1}_v, \textbf{v}, \textbf{v}^2\}. \tag{3.6} \label{eq:3.6}
\end{equation*}
Here $\textbf{1}_v \in \R^{(k+1)N_v}$ is the vector of all ones, $\vv$ is the v-grid in \eqref{eq:3.3}, and $\vv^2 \in \R^{(k+1)N_v}$ is the element-wise square of $\vv$. To attain this projection, we first introduce the weight function $\textbf{w}_M = $ exp$(-v^2/2)$ with exponential decay to ensure proper decay of the projected function as $v \to \infty.$ Along with the weight function, we introduce the weighted inner product and associated norm as 
\begin{equation*}
    \langle \textbf{f},\textbf{g}\rangle_{\textbf{w}_M} = \sum_{j,jg}f_{j,jg}g_{j,jg}w_{M,j,jg}\omega_{v,j,jg}, \hspace{6mm} \|\textbf{f}\|_{\textbf{w}_M} = \sqrt{\langle \textbf{f},\textbf{f}\rangle_{\textbf{w}_M}}, \tag{3.7} \label{eq:3.7}
\end{equation*}
where $\textbf{w}_M \in \R^{(k+1)N_v}$ with $w_{M,j,jg} = \textbf{w}_{M}(v_{j,jg})$.

To perform the conservative decomposition, we first scale the low-rank DG solution f with the weight function to ensure proper decay of v. This gives
\begin{equation*}
    \tilde{\textbf{f}} = \frac{1}{\textbf{w}_M}\odot_v \textbf{f} = \sum^r_{l=1}\left(C_l U^{(1)}_{l} \otimes \left(\frac{1}{\textbf{w}_M} \odot_v U^{(2)}_{l}\right)\right). \tag{3.8} \label{eq:3.8}
\end{equation*}
 Then we find an orthogonal projection of $\tilde{\textbf{f}}$ with respect to the inner product defined in \eqref{eq:3.7} onto the subspace $\mathcal{N}.$ That is, we search for an orthogonal projection such that
\begin{equation*}
    \langle P_{\mathcal{N}}(\tilde{\textbf{f}}),\textbf{g} \rangle_{\textbf{w}_M} = \langle \tilde{\textbf{f}},\textbf{g} \rangle_{\textbf{w}_M}, \text{ } \forall\textbf{g} \in \mathcal{N}. \tag{3.9} \label{eq:3.9}
\end{equation*}
With this orthogonal projection, we then follow the ideas of \cite{guo2024conservative} and compute a conservative decomposition as
\begin{equation*}
    \textbf{f} = \textbf{w}_M \star (P_\mathcal{N}(\tilde{f}) + (I - P_\mathcal{N})(\tilde{f})) = \textbf{w}_M \star (\tilde{\textbf{f}}_1 + \tilde{\textbf{f}}_2) = \textbf{f}_1 + \textbf{f}_2. \tag{3.10} \label{eq:3.10}
\end{equation*}
\cite{guo2024conservative} proved that this decomposition exists. Moreover, $\textbf{f}_1$ takes the form
\begin{equation*}
    \textbf{f}_1(\rho,\J,\kappa) = \frac{\rho}{\|\textbf{1}\|^2_{\textbf{w}_M}} \otimes (\textbf{w}_M \star \textbf{1}_v) + \frac{\J}{\| \vv \|^2_{\textbf{w}_M}} \otimes (\textbf{w}_M \star \vv) + \frac{2\kappa - c\rho}{\| \vv^2 - c \textbf{1}_v \|^2_{\textbf{w}_M}} \otimes (\textbf{w}_M \star (\vv^2 - c\textbf{1}_v)), \tag{3.11} \label{eq:3.11}
\end{equation*}
where $c = \frac{\langle \textbf{1}_v, \vv^2\rangle_{\textbf{w}_M}}{\| \textbf{1}_v \|^2_{\textbf{w}_M}}$. This new set $\{ \textbf{1}_v, \vv, \vv^2 - c\textbf{1}_v\}$ forms an orthogonal set of basis defined by the discrete mass, momentum, and kinetic energy density, $\rho,\J,\kappa$. And by construction, $\textbf{f}_1$ preserves the discrete mass, momentum, and kinetic energy of $\textbf{f},$ while the remainder $\textbf{f}_2 = \textbf{f} - \textbf{f}_1$ does not hold any information on them.

Since $\textbf{f}_2$ does not contain information about the discrete macroscopic quantities, we may truncate this term to remove redundancy in basis representation. To accomplish this, we simply perform a weighted SVD truncation. The weights consist of the quadrature weights associated with the quadrature nodes together with the weight function $\textbf{w}_M$ evaluated at those nodes. We define the weighted inner product for the SVD procedure as
\begin{equation*}
    \langle \textbf{f},\textbf{g} \rangle = \sum_{i,ig;j,jg} f^{ig,jg}_{i,j} g^{ig,jg}_{i,j} \omega_{x,i,ig} \omega_{v,j,jg} w_{M,j,jg}, \text{ } \textbf{f},\textbf{g} \in \R^{(k+1)N_x \times (k+1)N_v}. \tag{3.12} \label{eq:3.12}
\end{equation*}
The weighted SVD procedure is then given by
\begin{equation*}
    \sqrt{\textbf{$\omega$}\star \textbf{w}_M}\mathcal{T}_{\epsilon,\textbf{$\omega$}\star \textbf{w}_M}\left(\frac{\textbf{f}_2}{\sqrt{\textbf{$\omega$}\star \textbf{w}_M}}\right), \tag{3.13} \label{eq:3.13}
\end{equation*}
where $\omega = \omega_x \otimes \omega_v.$ Here, the SVD procedure is determined by a truncation threshold rather than fixing the number of singular values. While the truncation threshold depends on the problem considered, our results in section 5 show that for many standard test problems the range for truncation is set to be between $10^{-6}-10^{-4}$.

The end results of this conservative decomposition gives us two terms: (1) $\textbf{f}_1$, computed by \eqref{eq:3.11}, which holds the information on the discrete mass, momentum, and kinetic energy; and (2) the remainder $\textbf{f}_2 = \textbf{f} - \textbf{f}_1$ which is truncated using \eqref{eq:3.13} to remove any redundancies in the stored information of the solution. 
\newline

 \noindent \textbf{Remark 3.1:} We emphasize that although the above conservative decomposition resembles a micro-macro decomposition, it differs fundamentally from standard constructions such as those in \cite{coughlin2024robust}. Classical micro-macro approaches typically take the form $\mathcal{M} + \varepsilon g$, where $\mathcal{M}$ is a Maxwellian and $g$ is the microscopic component of the solution. Here $\varepsilon$ is a small parameter that vanishes at the asymptotic limit. This decomposition is particularly effective for asymptotic preserving methods such as \cite{coughlin2024robust}, especially when taking the magnetohydrodynamic (MHD) limit. In contrast, our decomposition, which takes the form $\bf{f}_1 + f_2$ as in \eqref{eq:3.10}, is a moment based decomposition. Although this does not set up our method to capture the correct asymptotic behavior in the vanishing Knudsen regime, it is specifically designed to guarantee local conservation of mass, momentum, and total energy at the discrete level.

%%%%%%%%%%%%%%%%%%%%%%%%%%%%%%%%%%%%%%%%%%%%%%%%%%%%%%%%%%%%%%%%%%%%%%%%%%%%%%
%%% Method & L2 Stability
%%%%%%%%%%%%%%%%%%%%%%%%%%%%%%%%%%%%%%%%%%%%%%%%%%%%%%%%%%%%%%%%%%%%%%%%%%%%%%
\section{The LoMaC method for collisional Vlasov dynamics}

In this section, we outline the proposed LoMaC low-rank approach equipped with DG discretization. The overall structure closely follows the framework introduced in \cite{guo2024LoMaCDG}, with key modifications introduced to incorporate the Dougherty Fokker-Planck collision operator.
\subsection{Method Outline at $t^n$}

Consider the solution $f^n$ of the form \eqref{eq:3.1} at timestep $t^n.$ The method to compute the next iteration $f^{n+1}$ is as follows.\\

\begin{flushleft}\textbf{\underline{Step 1. Add basis and obtain an intermediate solution $\textbf{f}^{n+1,*}$.}}\end{flushleft}

We first perform a second order multi-step discretization of the time derivative in (2.1),
\begin{equation*}
    f^{n+1,*} = \frac{1}{4}f^{n-2} + \frac{3}{4}f^n - \frac{3}{2}\Delta t(v\p_x(f^n) + E^n \p_v(f^n)) + \frac{3}{2} \nu \Delta t(\p_v (T (\p_v f^n)) + \p_v ((v-u)f^n)). \tag{4.1} \label{eq:4.1}
\end{equation*}
Here, the electric field $E^n$ is computed via a suitable Poisson solver. Assuming $f$ also has a low-rank representation at $t^{n-2},$ $\textbf{f}^{n+1,*}$ can then be represented in the following low-rank format: 
\begin{align*}
    f^{n+1,*} = &\frac{1}{4}\sum^{r^{n-2}}_{t = 1} C^{n-2}_{l}(U^{(1),n-2}_{l} \otimes U^{(2),n-2}_{l}) + \frac{3}{4}\sum^{r^n}_{l=1} C^n_l (U^{(1),n}_l \otimes U^{(2),n}_l) \tag{4.2a}\label{eq:4.2a}\\
    &-\frac{3}{2}\Delta t (D_x U^{(1),n}_l \otimes \vv \odot_v U^{(2),n}_l + \E^n \odot_x U^{(1),n}_l \otimes D_v U^{(2),n}_l) \tag{4.2b} \label{eq:4.2b}\\
    &+\frac{3}{2}\nu\Delta t(U^{(1),n}_l \otimes D^-_v[T(D^+_v U^{(2),n}_l)]) \tag{4.2c} \label{eq:4.2c}\\
    &+ \frac{3}{2}\nu\Delta t(U^{(1),n}_l \otimes D_v[\vv \odot_v U^{(2),n}_l] - \textbf{u} \odot_x U^{(1),n}_l \otimes D_v[U^{(2),n}_l]) \tag{4.2d}\label{eq:4.2d}
\end{align*}
Here, $\vv \in \R^{N_v}$ denotes the coordinates of $v_{grid}$ introduced in \eqref{eq:3.3}. $D_x$ and $D_v$ represent high order spatial differentiations, and $\star$ denotes an element-wise multiplication operation. For example the discretization of the terms $D_x U^{(1),n}_l \otimes \vv \odot_v U^{(2),n}_l$ and $\E^n \odot_x U^{(1),n}_l \otimes D_v U^{(2),n}_l$ from the Vlasov equation follows
\begin{align*}
    D^+_x \textbf{U}^{(1),n}_l \otimes \vv^+ \odot_v \textbf{U}^{(2),n}_l &+ D^-_x\textbf{U}^{(1),n}_l \otimes \vv^- \odot_v \textbf{U}^{(2),n}_l,\tag{4.3}\label{eq:4.5}\\
    \E^{n,+} \odot_x \textbf{U}^{(1),n}_l \otimes D^+_v \textbf{U}^{(2),n}_l &+ \E^{n,-}\odot_x \textbf{U}^{(1),n}_l \otimes D^-_v \textbf{U}^{(2),n}_l.\tag{4.4}\label{eq:4.6}
\end{align*}
where $D^\pm_x$ and $D^\pm_v$ are defined as in \eqref{eq:2.4a}-\eqref{eq:2.4f}.
As for computing $T$ and $\bf{u}$, since we are simultaneously updating the macroscopic quantities $\rho,\J,k$, we can instead express $T$ and $\uu$ in terms of these quantities. Here, $T = 2\textbf{$k$}/\textbf{$\rho$} - (\J/\textbf{$\rho$})^2$ and $\uu = \J/\textbf{$\rho$}$.

\begin{flushleft}
    \textbf{\underline{Step 2. Perform a macroscopic conservative decomposition.}}
\end{flushleft}

Once the intermediate solution has been computed, we then perform the macroscopic conservative decomposition
\begin{equation*}
    \textbf{f}^{n+1,*} = \textbf{f}_1 + \textbf{f}_2.\tag{4.5}\label{eq:4.7}
\end{equation*}
Here, $\textbf{f}_1$ is computed from \eqref{eq:3.10}. $\textbf{f}_2 = \textbf{f} - \textbf{f}_1$ is the remainder term, where we will apply a weighted SVD truncation.

\begin{flushleft}
    \textbf{\underline{Step 3. Conservative update of macroscopic variables.}}
\end{flushleft}

We update our macroscopic observables using the reinterpreted macroscopic system
\begin{equation*}
    U_t + F_x = S, \tag{4.6}\label{eq:4.8}
\end{equation*}
where $U = (\rho,J,e)^\top$, $F = (J,\sigma,\textbf{Q})^\top$, and $S = (0,\rho E,0)^\top$. Denote the numerical solutions for $U$ as $\textbf{$\rho$}^M,\J^M,\textbf{$\kappa$}^M$ ($M$ stands for "Macroscopic variables"). We can advance the system \eqref{eq:1.10} similar to Step 1. Using a second order SSP multi-step time integrator, we look for an updated $U^{n+1}$ whose nodal values satisfy the following system:
\begin{equation*}
    U^{n+1}_{i,ig} = \frac{1}{4}U^{n-2}_{i,ig} + \frac{3}{4}U^{n}_{i,ig} + \frac{3}{2}\Delta t(D^+_{x,i,ig}\textbf{F}^{n,+}_{i,:} + D^-_{x,i,ig}\textbf{F}^{n,-}_{i,:}+ S^{n}_{i,ig}),\tag{4.7}\label{eq:4.9}
\end{equation*}
where $U^n_{i,ig} = (\rho^n_{i,ig}, J^n_{i,ig},e^n_{i,ig})^\top$ and $S^n_{i,ig}=(0,\rho^n_{i,ig}E^n_{i,ig},0)^\top, i=1,...,N_x, ig = 0,...,k.$ $\textbf{F}^{n,\pm}\in \R^{(k+1)N_x}$ are given by the kinetic flux vector splitting scheme \cite{guo2024LoMaC} with 
\begin{align*}
    \textbf{F}^{n,+} &= \sum^{r^n}_{l=1}C^n_l \left \langle \textbf{U}^{(2),n}_l , \begin{pmatrix} \vv^+ \\ (\vv^+)^2 \\ \frac{1}{2}(\vv^+)^3 \end{pmatrix}\right\rangle_v \textbf{U}^{(1),n}_l,\tag{4.8}\label{eq:4.10}\\
    \textbf{F}^{n,-} &= \sum^{r^n}_{l=1}C^n_l \left \langle \textbf{U}^{(2),n}_l , \begin{pmatrix} \vv^- \\ (\vv^-)^2 \\ \frac{1}{2}(\vv^-)^3 \end{pmatrix}\right\rangle_v \textbf{U}^{(1),n}_l.\tag{4.9}\label{eq:4.11}
\end{align*}
Here, $\vv^+ = \text{max}(v,0)$, $\vv^- = \text{min}(v,0)$, and the inner product is defined as in \eqref{eq:3.5}. $D^{\pm}_{x,i,ig}$ are defined as in step 1, and

\begin{align*}
    \textbf{F}^{n,+}_{i,:} &= (F^{n,+}_{i-1,0},...,F^{n,+}_{i-1,k},F^{n,+}_{i,0},...,F^{n,+}_{i,k}),\\
    \textbf{F}^{n,-}_{i,:} &= (F^{n,-}_{i-1,0},...,F^{n,-}_{i-1,k},F^{n,-}_{i,0},...,F^{n,-}_{i,k}).
\end{align*}
Once $U^{n+1}_{i,ig}$ has been obtained, we can then compute the updated kinetic energy
\begin{equation*}
    \kappa^{n+1,M}_{i,ig} = e^{n+1,M}_{i,ig} - \frac{1}{2}|E^{n+1,M}_{i,ig}|^2.\tag{4.10}\label{eq:4.12}
\end{equation*}
Here, $\E^{n+1,M}$ has been directly computed from $\textbf{$\rho$}^{n+1,M}$ via Poisson's equation using the local DG method \cite{coughlin2022efficient}. We finally construct $\textbf{f}^{M}_1$ using equation \eqref{eq:3.10} with the updated macroscopic observables $\textbf{$\rho$}^{n+1,M}, \J^{n+1,M}$, and $\kappa^{n+1,M}$.

\begin{flushleft}
    \textbf{\underline{Step 4. Update the low-rank solution.}}
\end{flushleft}

Lastly, we update our low-rank solution

\begin{equation*}
    \textbf{f}^{n+1} = \textbf{f}^M_1 + \mathcal{T}_{\varepsilon, \omega
\star \textbf{w}_M} (\textbf{f}_2).\tag{4.11}\label{eq:4.13}\end{equation*}
Here, $\textbf{f}^M_1$ is computed from Step 3 and $\mathcal{T}_{\varepsilon, \omega
\star \textbf{w}_M} (\textbf{f}_2)$ is the truncated weighted SVD of $\textbf{f}_2$ as in \eqref{eq:3.13}. $\textbf{f}^M_1$ is a correction to $\textbf{f}_1$ from step 2 and is constructed to contain all information of the mass, momentum, and energy densities at a fixed small rank. Whereas $\mathcal{T}_{\varepsilon, \omega\star \textbf{w}_M} (\textbf{f}_2)$ holds no information on the macroscopic quantities and is truncated to remove any redundancies in the basis representation of the solution.\\

\begin{flushleft}\textbf{Remark 4.1:} Similar to our previous work in \cite{guo2024LoMaCDG,guo2024LoMaC}, the above DG algorithm can be generalized to higher dimensions (e.g. 2D2V) using the hierarchical Tucker (HT) format \cite{tucker1966some,grasedyck2010hierarchical}. This was initially discussed in \cite{guo2024LoMaC} with extensions to the DG framework discussed in \cite{guo2024LoMaCDG}. We refer to them for more in depth discussions of extensions to higher dimensions.
\end{flushleft}

\section{Numerical Results}
%%%%%%%%%%%%%%%%%%%%%%%%%%%%%%%%%%%%%%%%%%%%%%%%%%%%%%%%%%%%%%%%%%%%%%%%%%%%%%
%%% Weak Landau Damping
%%%%%%%%%%%%%%%%%%%%%%%%%%%%%%%%%%%%%%%%%%%%%%%%%%%%%%%%%%%%%%%%%%%%%%%%%%%%%%
In all numerical examples, we primarily investigate the influence of the collision frequency $\nu.$ All simulations are performed using $P^2$ discontinuous Galerkin (DG) polynomials on a $(x \times v)$ mesh of $50 \times 100$ elements.
\begin{flushleft}
\textbf{Example 5.1.} (Weak Landau Damping.) 
\end{flushleft}

\begin{figure}
    \centering
    \begin{subfigure}[t]{.45\textwidth}
        \centering
        \includegraphics[scale=0.4]{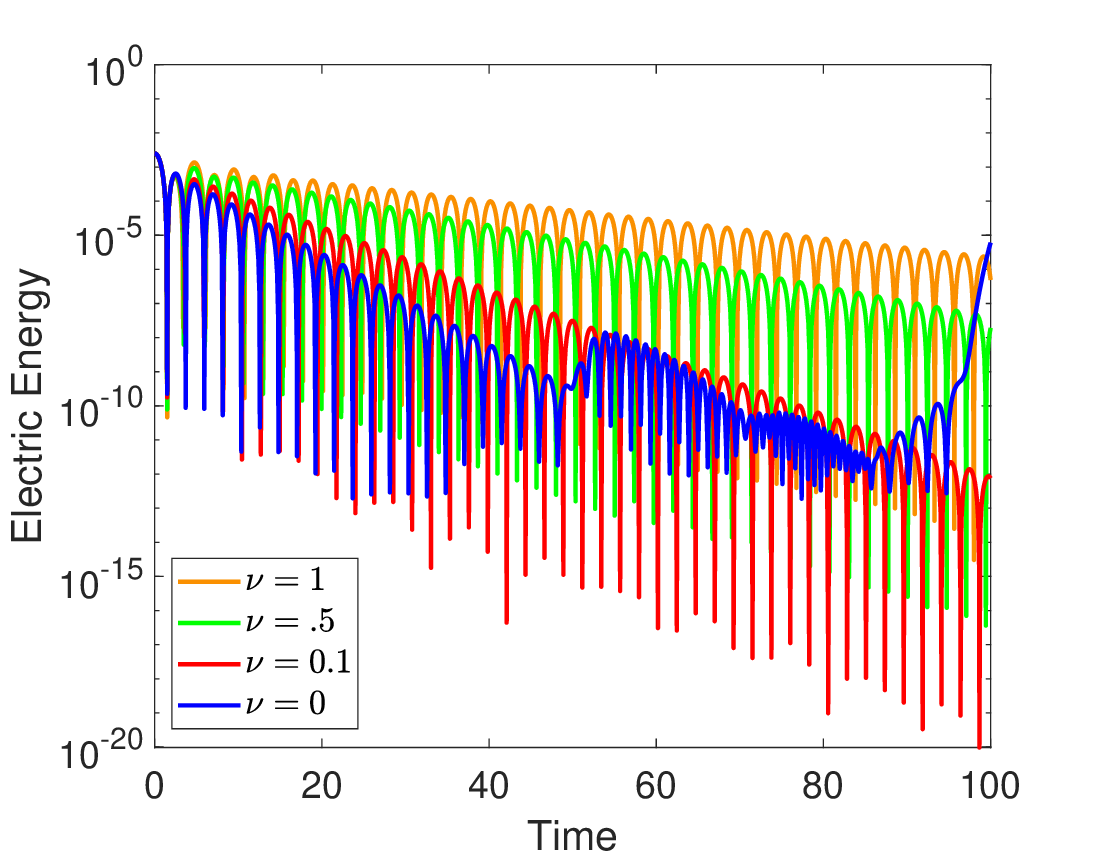}
        \caption{}
        \label{fig1:subim1}
    \end{subfigure}
    \begin{subfigure}[t]{0.45\textwidth}
        \centering
        \includegraphics[scale=0.4]{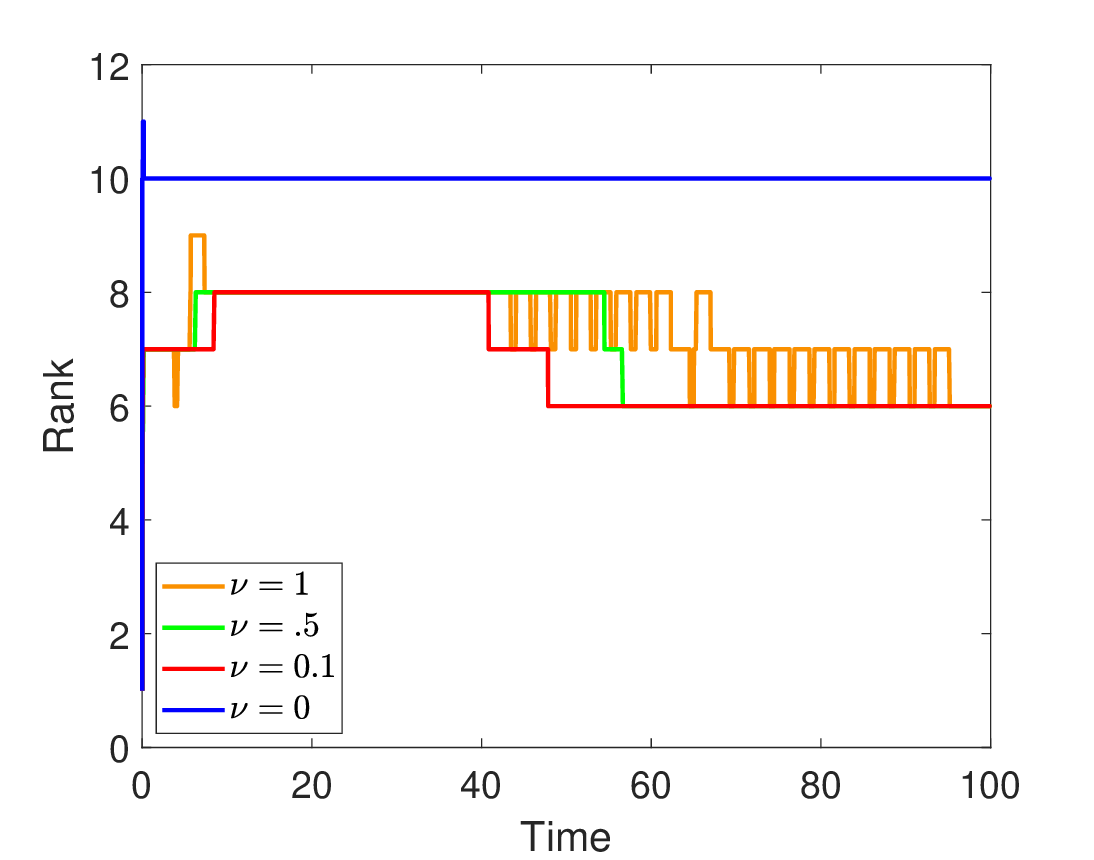}
        \caption{}
        \label{fig1:subim2}
    \end{subfigure}
    \hfill
    \begin{subfigure}[t]{0.45\textwidth}
        \centering
        \includegraphics[scale=0.4]{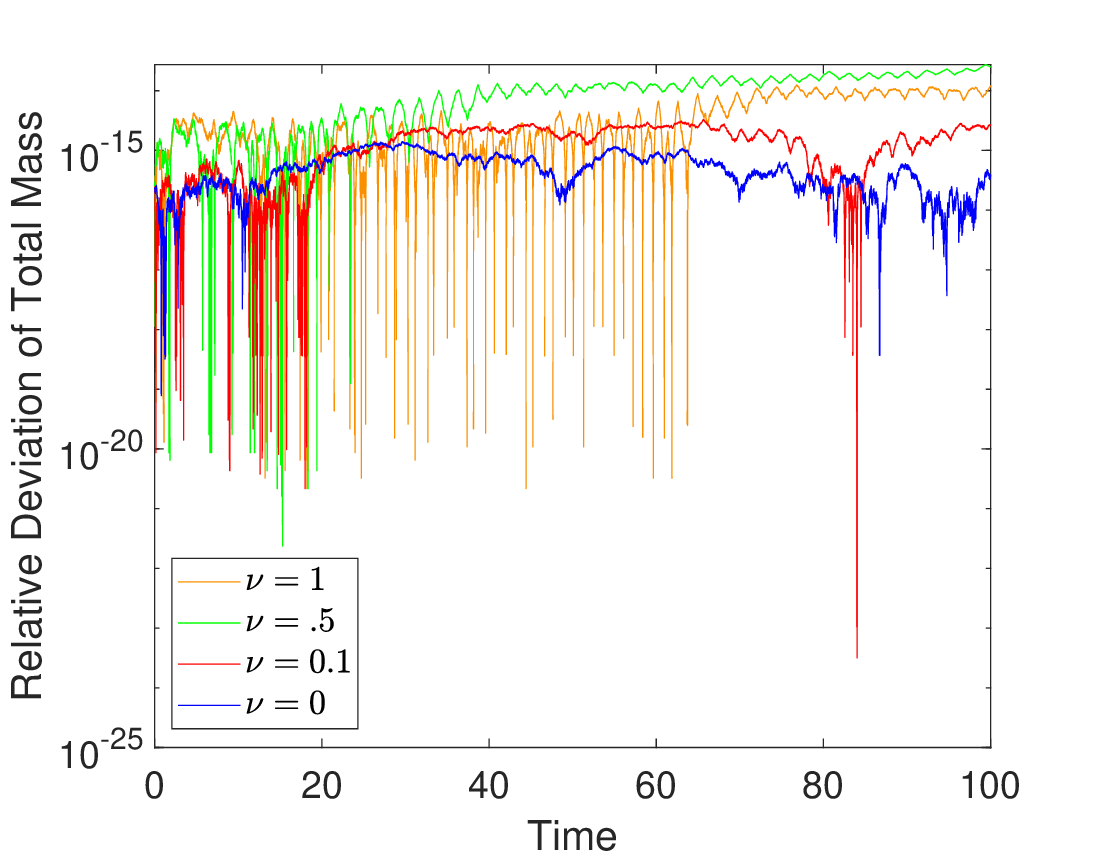}
        \caption{}
        \label{fig1:subim3}
    \end{subfigure}
    \begin{subfigure}[t]{0.45\textwidth}
        \centering
        \includegraphics[scale=0.4]{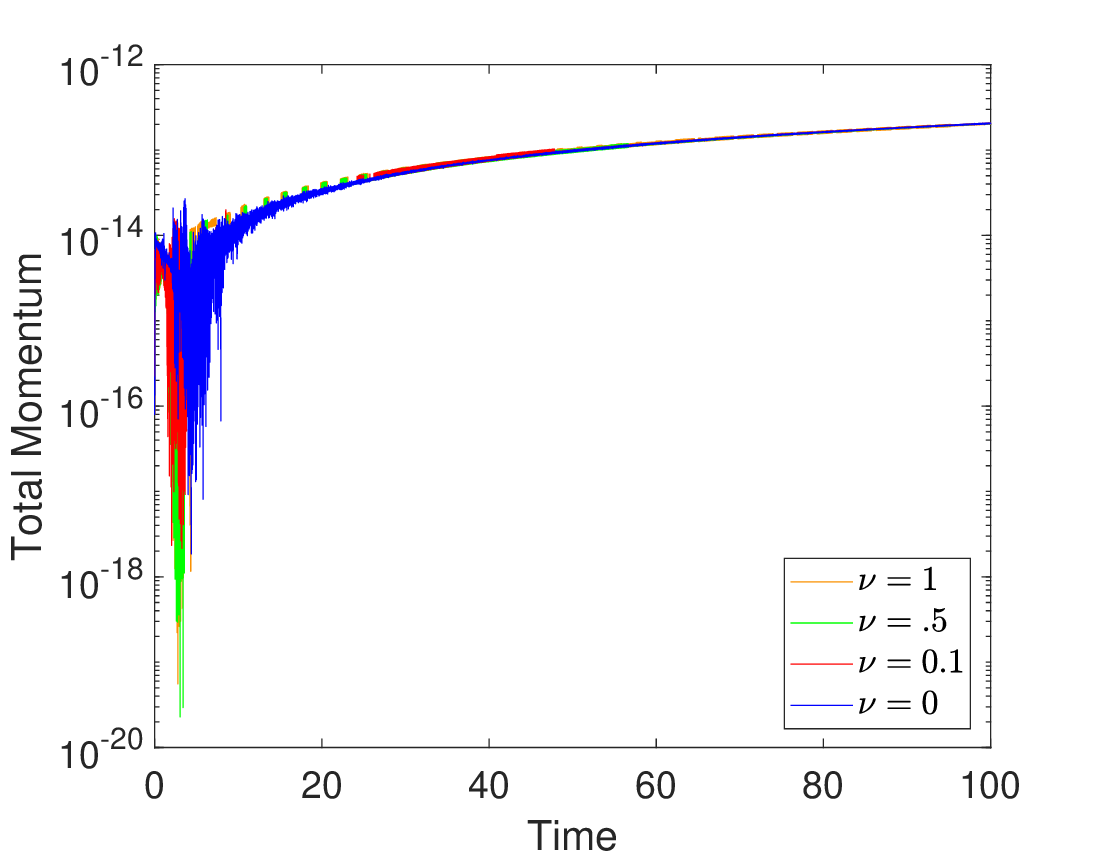}
        \caption{}
        \label{fig1:subim4}
    \end{subfigure}
    \hfill
    \begin{subfigure}[t]{0.45\textwidth}
        \centering
        \includegraphics[scale=0.4]{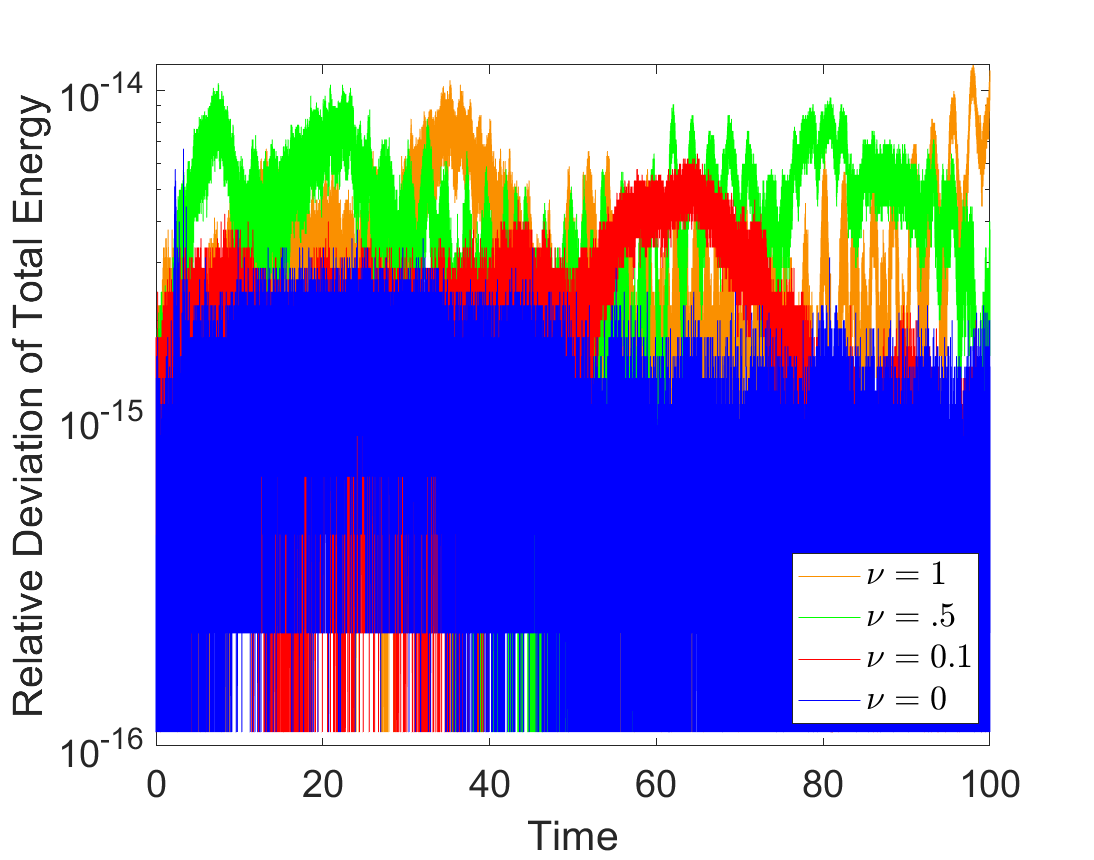}
        \caption{}
        \label{fig1:subim5}
    \end{subfigure}
    \hfill
    \caption{Weak Landau damping (Example 5.1). Time histories of the electric energy and numerical rank for $\nu = 0, 0.1, 0.5,$ and $1$ are shown in (a)-(b). As the collision frequency increases, the damping rate of the electric energy decreases. Unlike later benchmark problems, filamentary structures do not significantly develop in this setting, and the plasma remains inherently low rank throughout the simulation. As a result, the collisional effects have a comparatively small impact on the numerical ranks. The relative deviation of mass and total energy and the absolute total momentum are plotted in (c)-(e), confirming the conservation properties of the method.}
\end{figure}

We first simulate the weak Landau damping test with initial condition
\begin{equation*}
    f(x,v,t = 0) = \frac{1}{\sqrt{2\pi}}(1 + \alpha \cos(\kappa x))\text{exp}\left ( -\frac{v^2}{2} \right )
\end{equation*}
where $\alpha = 0.01$ and $\kappa = 0.5$. The computational domain is set to be $[0,L_x]$ $\times$ $[-L_v,L_v]$ with $L_x = 2\pi / k$ and $L_v = 6.$ We set the truncation threshold to be $\varepsilon = 10^{-5}$. Figures \hyperref[fig1:subim1]{1a-1b} present the time histories of the electric energy and numerical ranks of the low-rank DG solutions for collision frequencies $\nu = 0,0.1,0.5,$ and $1$. As the collision frequency increases, the damping rate decreases, which is consistent with the results reported in \cite{ye2024energy,}. In Figures \hyperref[fig1:subim3]{1c-1e}, we further report the time history of the relative deviation of the total mass and total energy, together with the absolute total momentum.

%%%%%%%%%%%%%%%%%%%%%%%%%%%%%%%%%%%%%%%%%%%%%%%%%%%%%%%%%%%%%%%%%%%%%%%%%%%%%%
%%% Strong Landau Damping
%%%%%%%%%%%%%%%%%%%%%%%%%%%%%%%%%%%%%%%%%%%%%%%%%%%%%%%%%%%%%%%%%%%%%%%%%%%%%%
\begin{flushleft}
\textbf{Example 5.2.} (Strong Landau Damping.) 
\end{flushleft}

\begin{figure}
    \centering
    \begin{subfigure}[t]{.45\textwidth}
        \centering
        \includegraphics[scale=0.4]{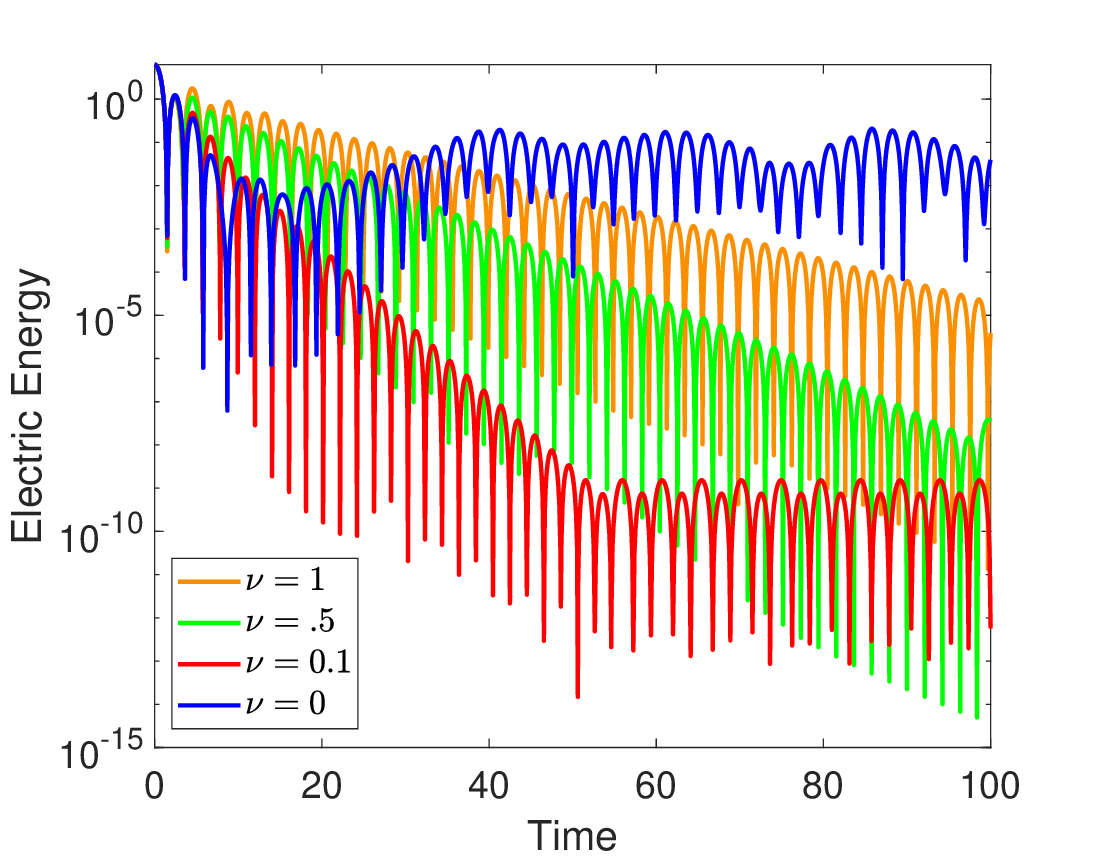}
        \caption{}
        \label{fig2:subim1}
    \end{subfigure}
    \begin{subfigure}[t]{0.45\textwidth}
        \centering
        \includegraphics[scale=0.4]{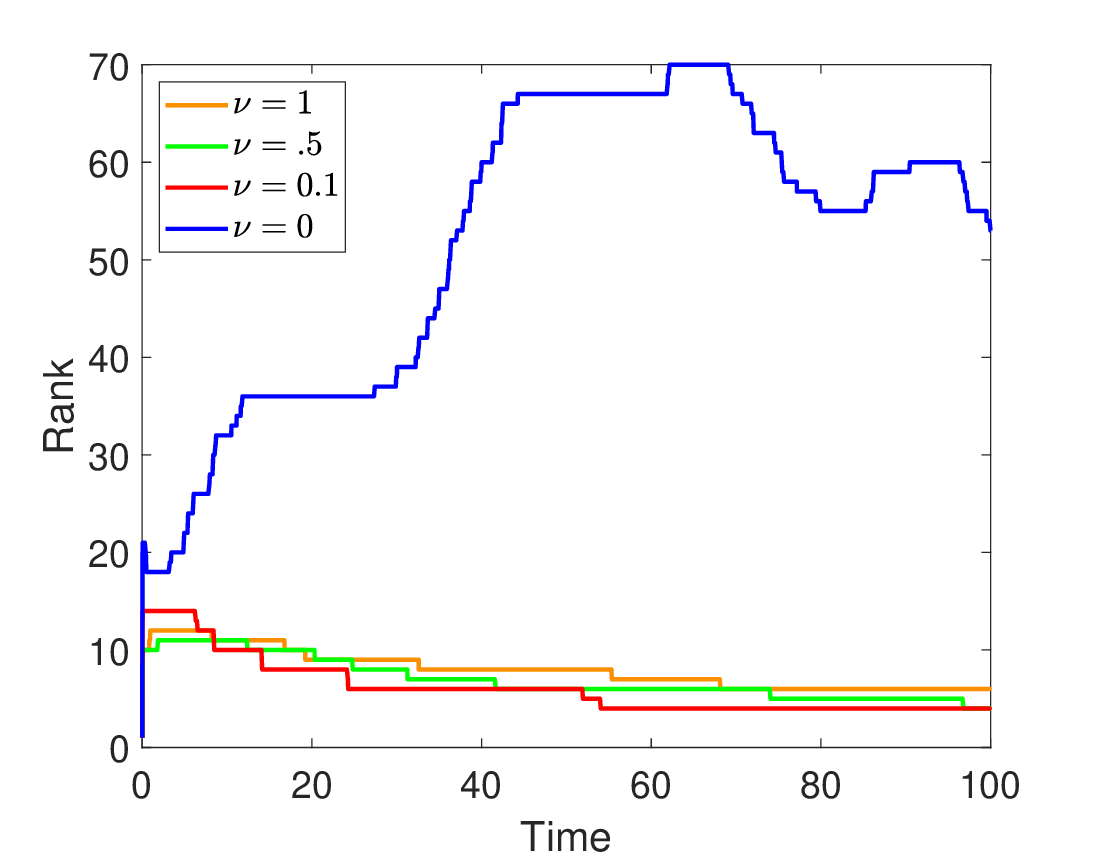}
        \caption{}
        \label{fig2:subim2}
    \end{subfigure}
    \hfill
    \begin{subfigure}[t]{0.45\textwidth}
        \centering
        \includegraphics[scale=0.4]{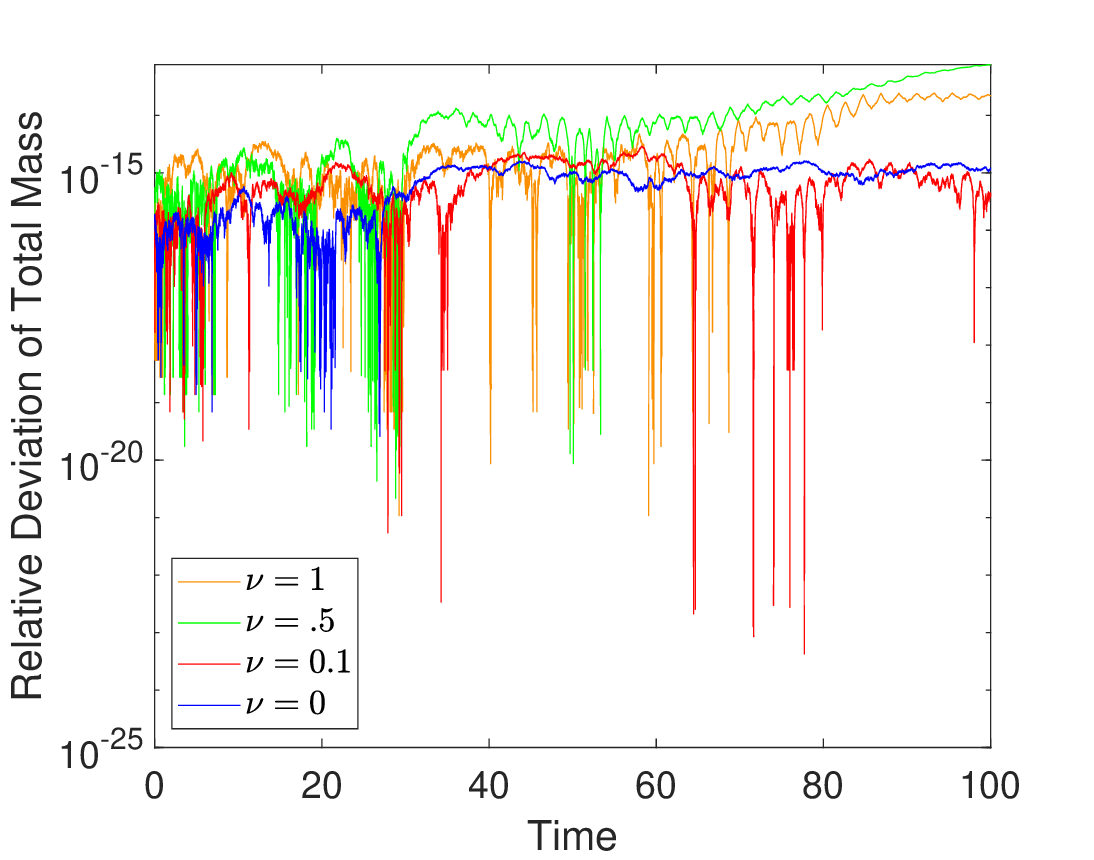}
        \caption{}
        \label{fig2:subim3}
    \end{subfigure}
    \begin{subfigure}[t]{0.45\textwidth}
        \centering
        \includegraphics[scale=0.4]{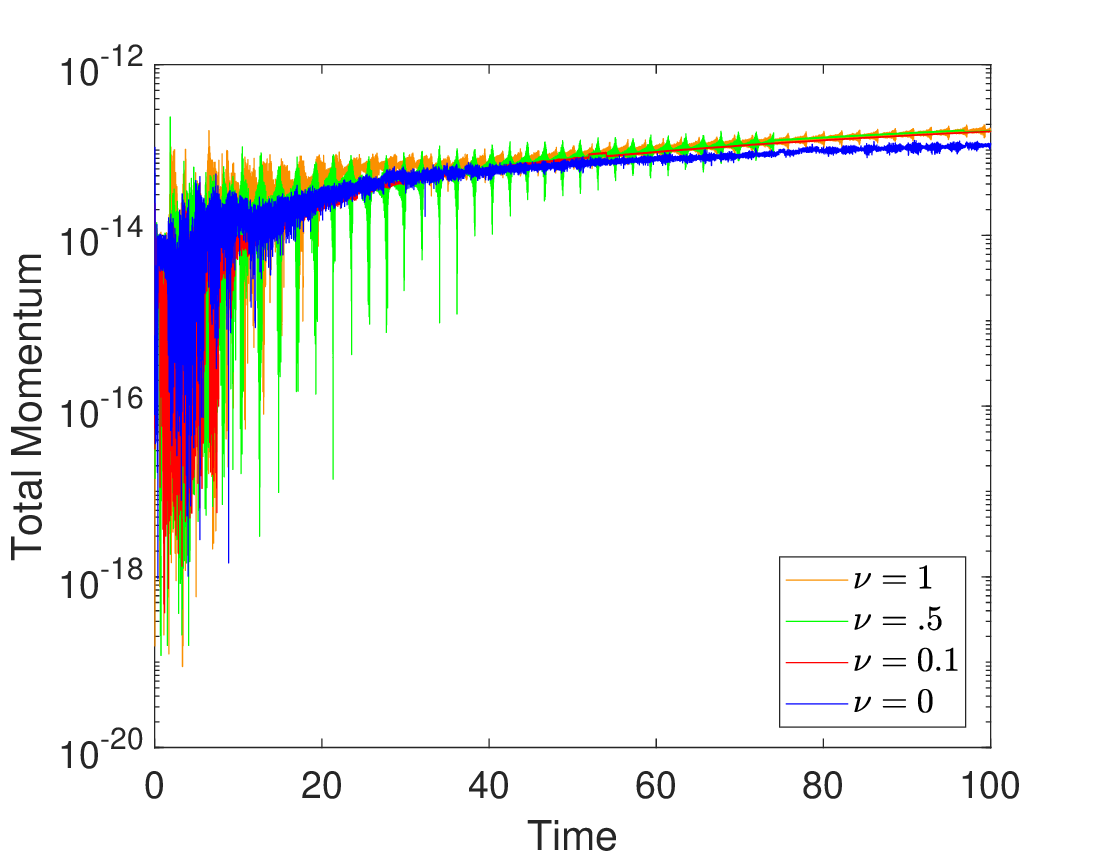}
        \caption{}
        \label{fig2:subim4}
    \end{subfigure}
    \hfill
    \begin{subfigure}[t]{0.45\textwidth}
        \centering
        \includegraphics[scale=0.4]{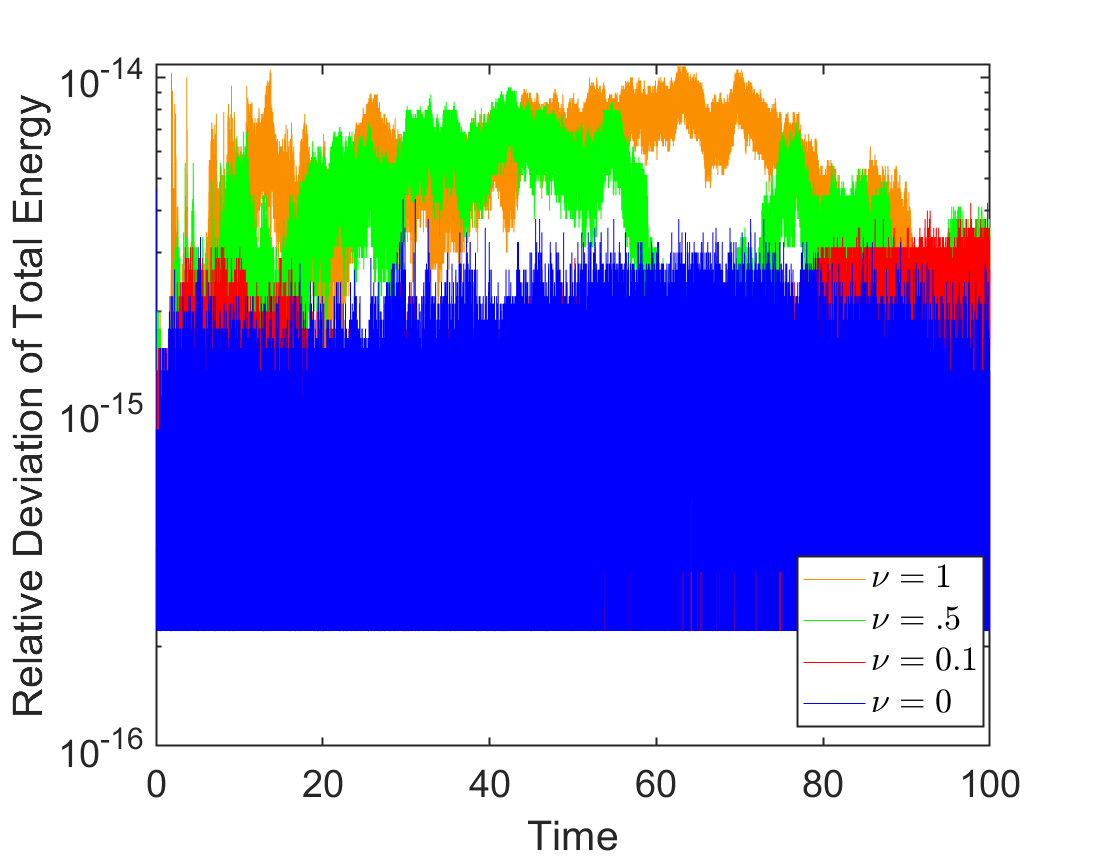}
        \caption{}
        \label{fig2:subim5}
    \end{subfigure}
    \hfill
    \caption{Strong Landau damping (Example 5.2). Time histories of the electric energy and numerical rank for $\nu = 0, 0.1, 0.5,$ and $1$ are shown in (a)-(b). As the collision frequency increases, the damping rate of the electric energy decreases. Most importantly, the collisional effects greatly suppress the formation of filamentary structures, thereby significantly limiting rank growth. The relative deviation of mass and total energy and the absolute total momentum are plotted in (c)-(e).}
    \label{fig2}
\end{figure}

\begin{figure}
    \centering
    \begin{subfigure}[t]{.45\textwidth}
        \centering
        \includegraphics[scale=0.4]{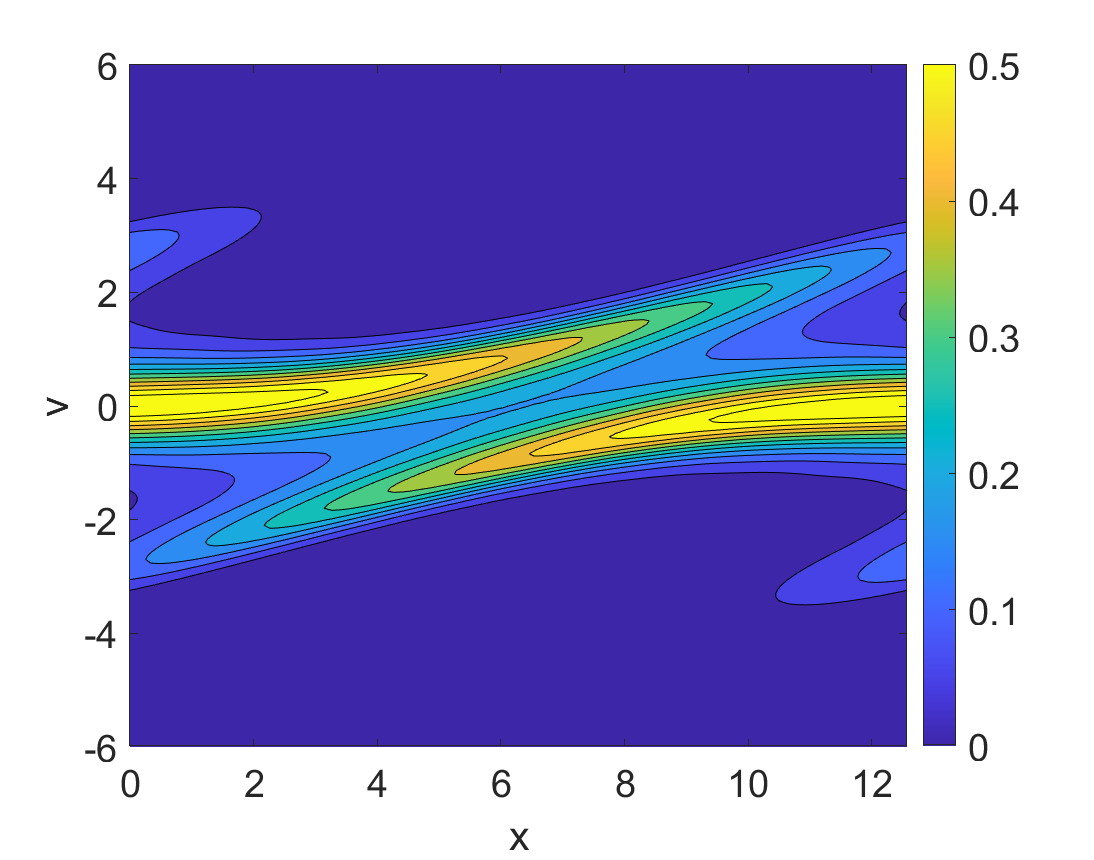}
        \caption{T = 1.}
        \label{fig3:subim1}
    \end{subfigure}
    \begin{subfigure}[t]{0.45\textwidth}
        \centering
        \includegraphics[scale=0.4]{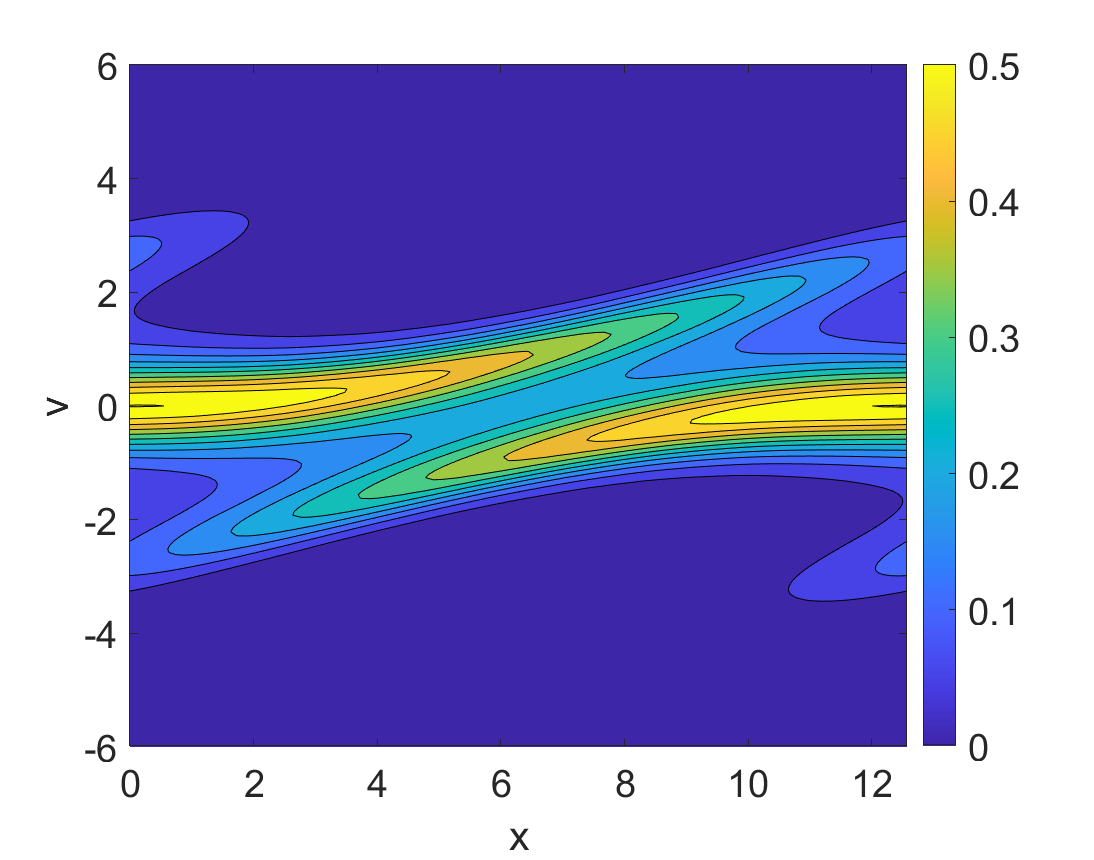}
        \caption{T = 1.}
        \label{fig3:subim2}
    \end{subfigure}
    \hfill
    \begin{subfigure}[t]{0.45\textwidth}
        \centering
        \includegraphics[scale=0.4]{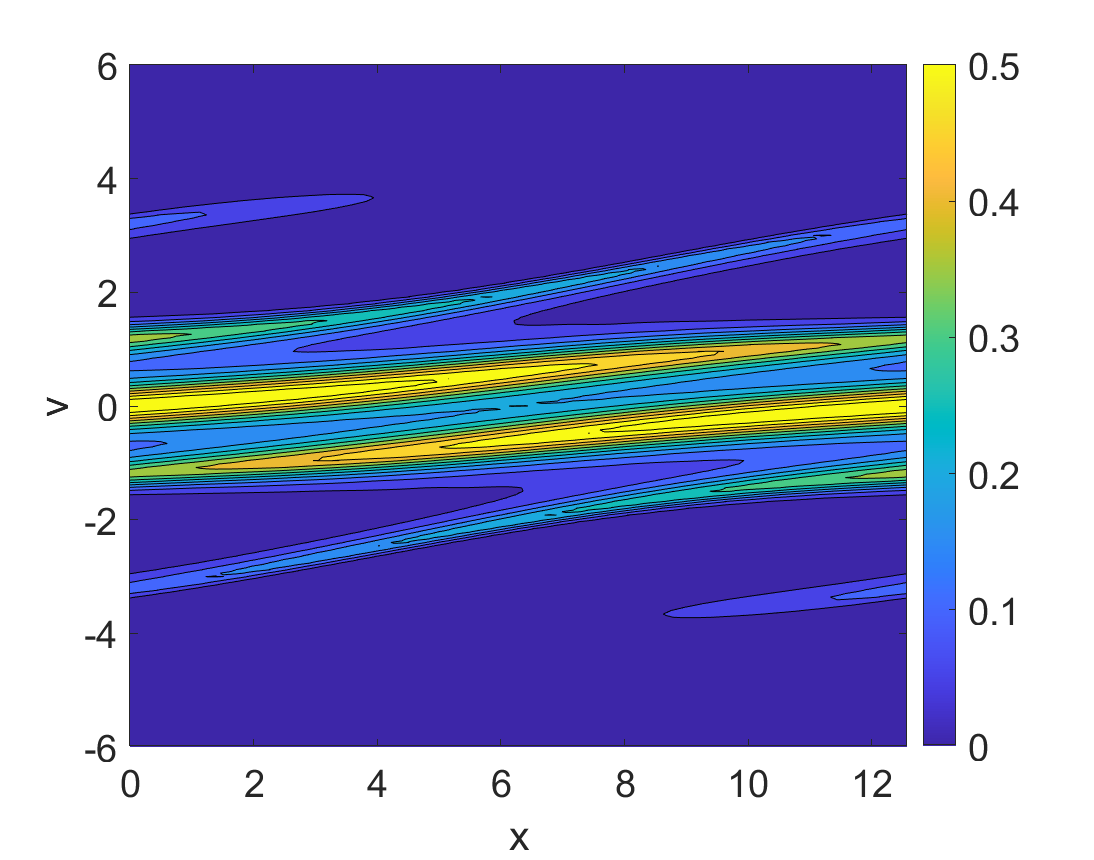}
        \caption{T = 8.}
        \label{fig3:subim3}
    \end{subfigure}
    \begin{subfigure}[t]{0.45\textwidth}
        \centering
        \includegraphics[scale=0.4]{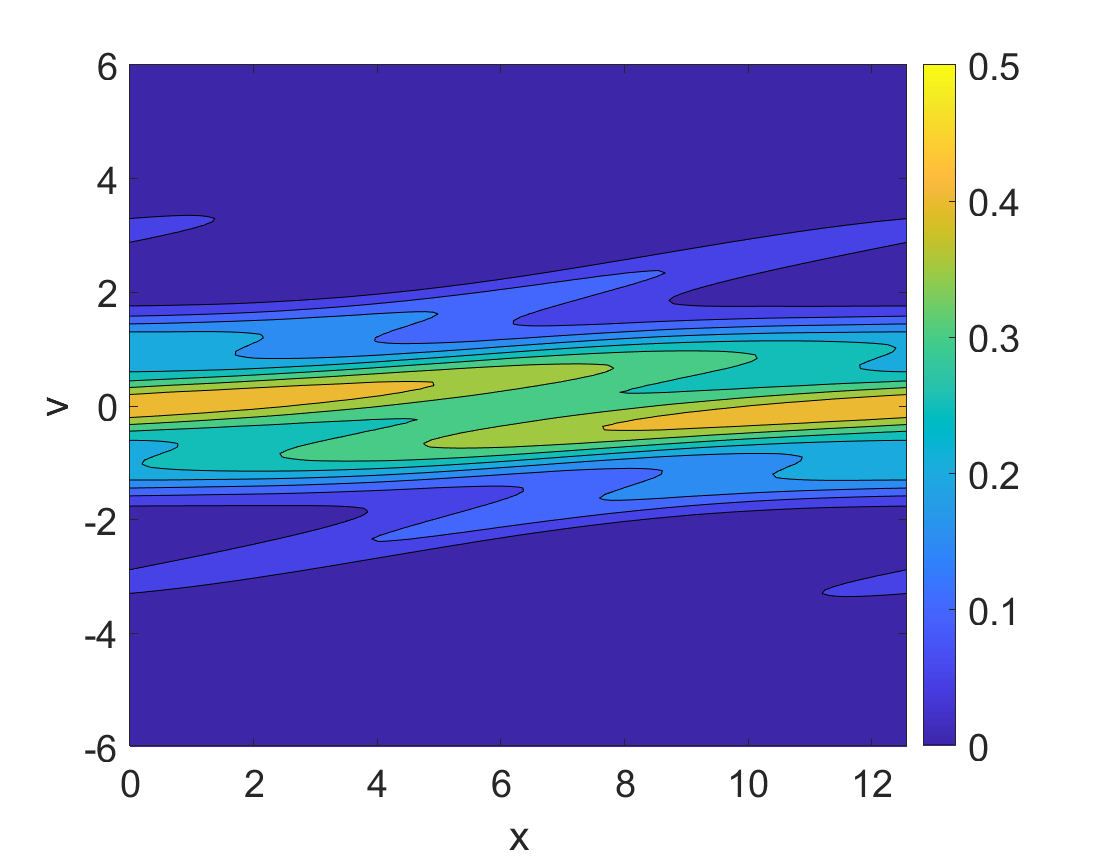}
        \caption{T = 8.}
        \label{fig3:subim1}
    \end{subfigure}
    \hfill
    \begin{subfigure}[t]{0.45\textwidth}
        \centering
        \includegraphics[scale=0.4]{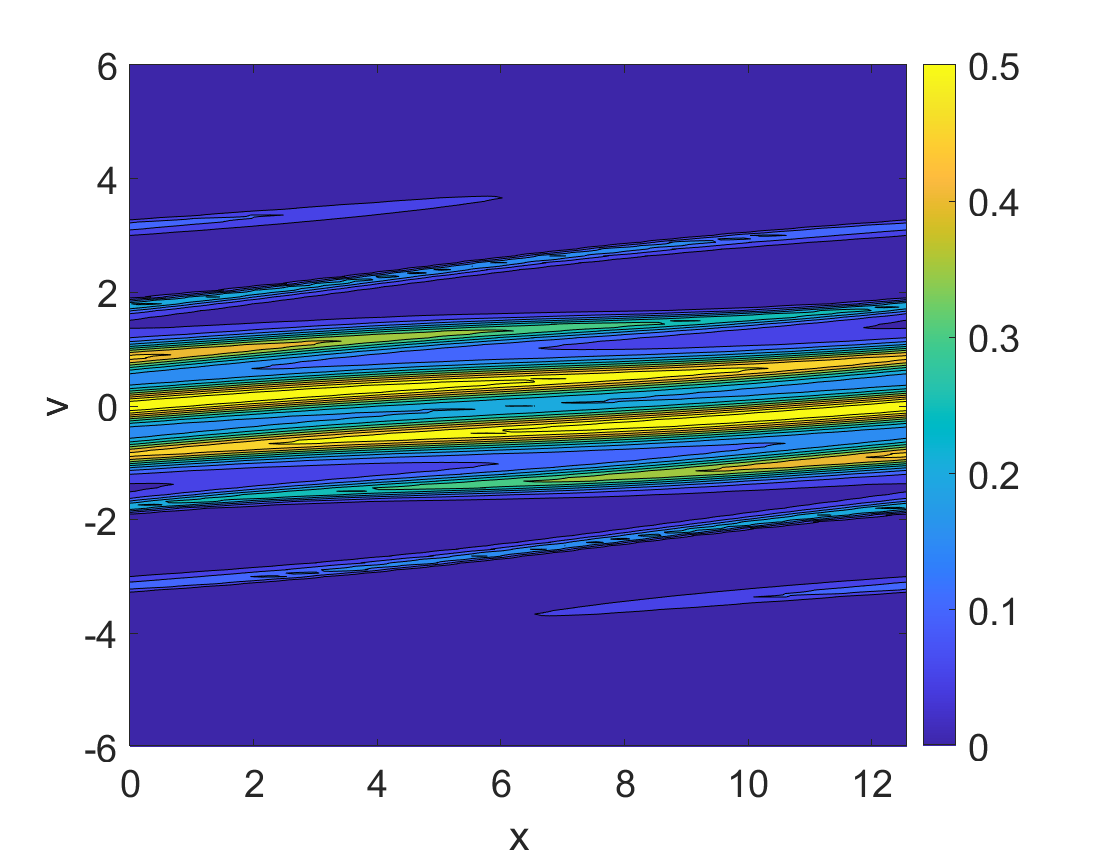}
        \caption{T = 12.}
        \label{fig3:subim2}
    \end{subfigure}
    \begin{subfigure}[t]{0.45\textwidth}
        \centering
        \includegraphics[scale=0.4]{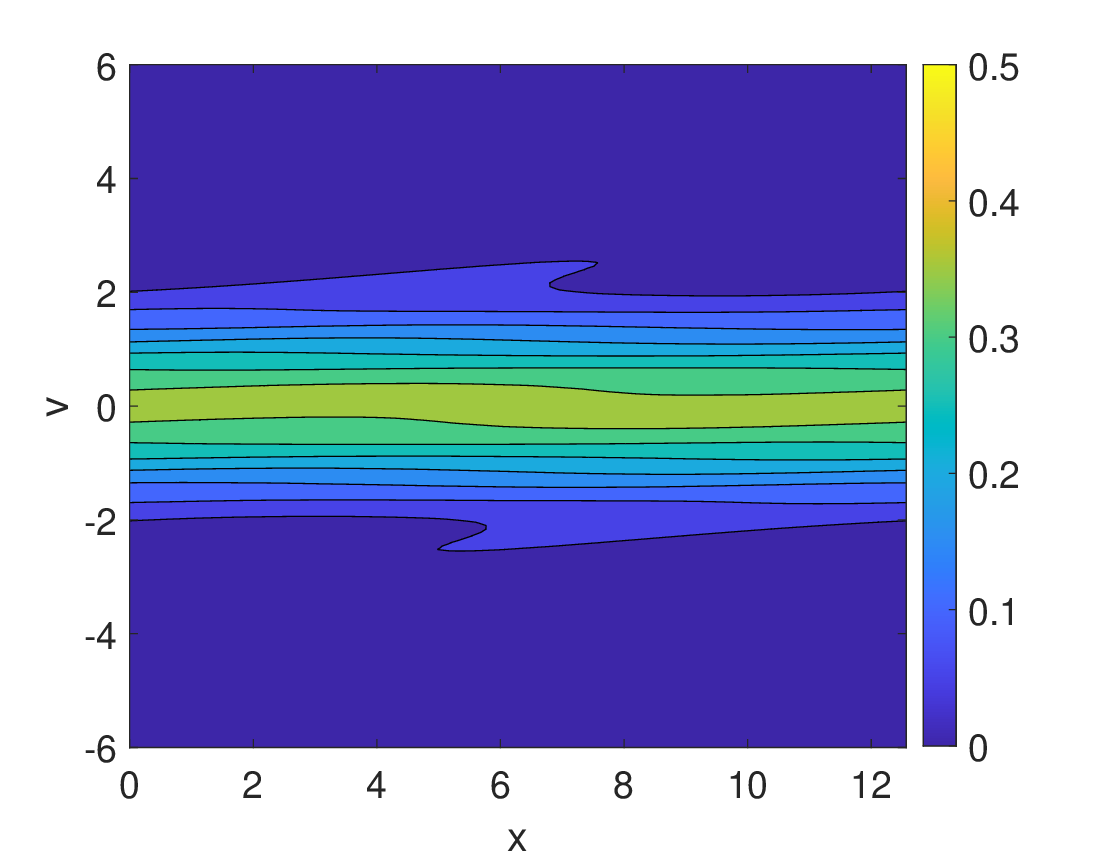}
        \caption{T = 12.}
        \label{fig3:subim3}
    \end{subfigure}
    \caption{Strong Landau damping example with $\nu = 0$ on the left (a,c,e) and $\nu = 0.01$ on the right (b,d,f). Phase-space contour plots are shown at different times. For the collisionless case, $\nu = 0$, filamentary structures continually develop, leading to significant rank growth. In contrast, even the weak collisional effects associated with $\nu = 0.01$ are sufficient to gradually suppress these structures and drive the plasma toward equilibrium.}
    \label{fig3}
\end{figure}

Next, we consider the strong Landau damping test. The initial conditions are identical to those in Example 5.1, with $\alpha = 0.5$, and the truncation threshold set to $10^{-4}.$ We again use collision frequencies $\nu = 0, 0.1, 0.5, $ and $1$. Figures \hyperref[fig2:subim1]{2a-2b} show the time histories of the electric energy and numerical ranks of the DG solutions. As the collisional effects of the Fokker-Planck operator increase, the plasma is driven more rapidly toward equilibrium, and the resulting low-rank structures dominate the solution for most of the simulation. Consequently, the numerical ranks remain significantly lower than in the collisionless case. Even for relatively small collision frequencies, the collisional effects have a substantial influence on the evolution of the distribution function. We illustrate this behavior in \hyperref[fig3]{Figure 3}. In the collisionless case, (3a,3c,3e), thin filamentation structures continually develop and persist throughout the simulation. In contrast, for $\nu = 0.01,$ these filamentary structures are gradually suppressed, and the plasma begins to relax toward equilibrium. Furthermore, Figures \hyperref[fig2:subim1]{2c-2e} show that the method conserves the physical invariants as expected.

%%%%%%%%%%%%%%%%%%%%%%%%%%%%%%%%%%%%%%%%%
%%% Strong Landau EE and Ranks
%%%%%%%%%%%%%%%%%%%%%%%%%%%%%%%%%%%%%%%%%

%%%%%%%%%%%%%%%%%%%%%%%%%%%%%%%%%%%%%%%%%
%%% Strong Landau nu = 0.01.
%%%%%%%%%%%%%%%%%%%%%%%%%%%%%%%%%%%%%%%%%

%%%%%%%%%%%%%%%%%%%%%%%%%%%%%%%%%%%%%%%%%%%%%%%%%%%%%%%%%%%%%%%%%%%%%%%%%%%%%%
%%% 2-Stream Instability
%%%%%%%%%%%%%%%%%%%%%%%%%%%%%%%%%%%%%%%%%%%%%%%%%%%%%%%%%%%%%%%%%%%%%%%%%%%%%%
\begin{flushleft}
    \textbf{Example 5.3.} (Two-Stream Instability.) 
\end{flushleft} 
    
\begin{figure}
    \centering
    \begin{subfigure}[t]{.45\textwidth}
        \centering
        \includegraphics[scale=.4]{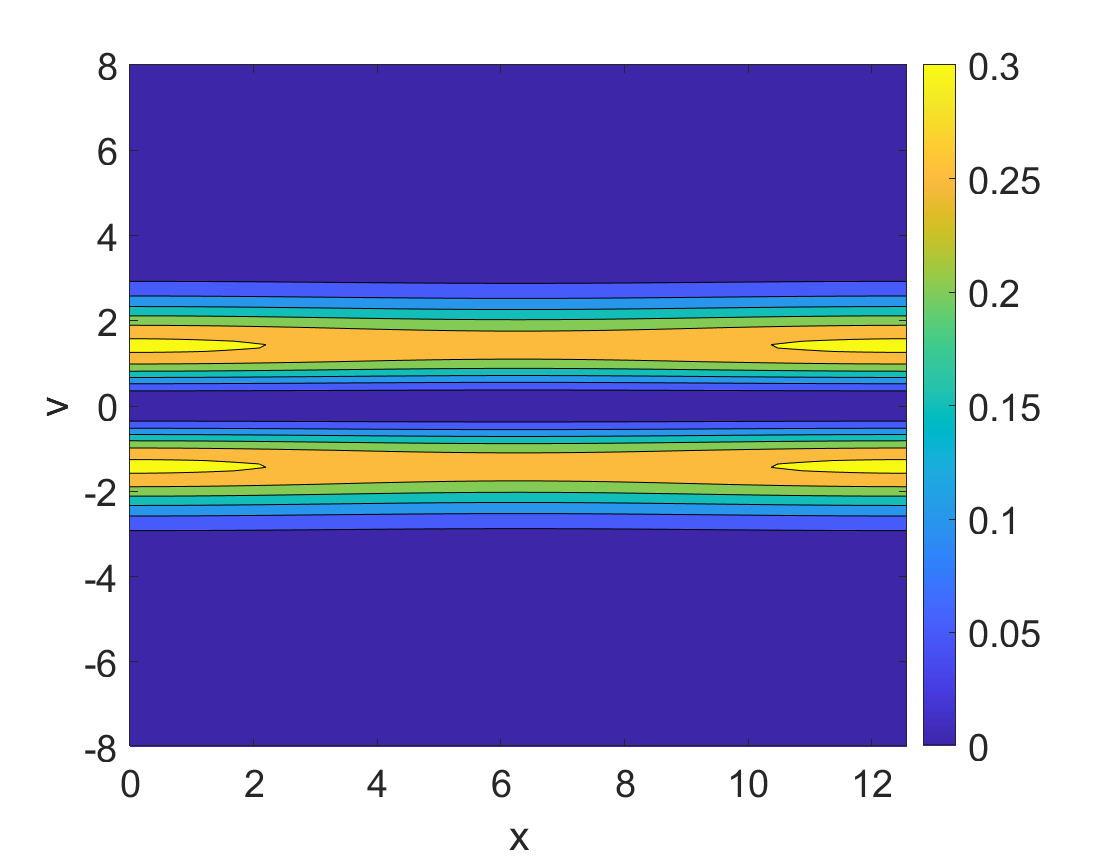}
        \caption{T = 0.}
        \label{fig4:subim1}
    \end{subfigure}
    \begin{subfigure}[t]{.45\textwidth}
        \centering
        \includegraphics[scale=.4]{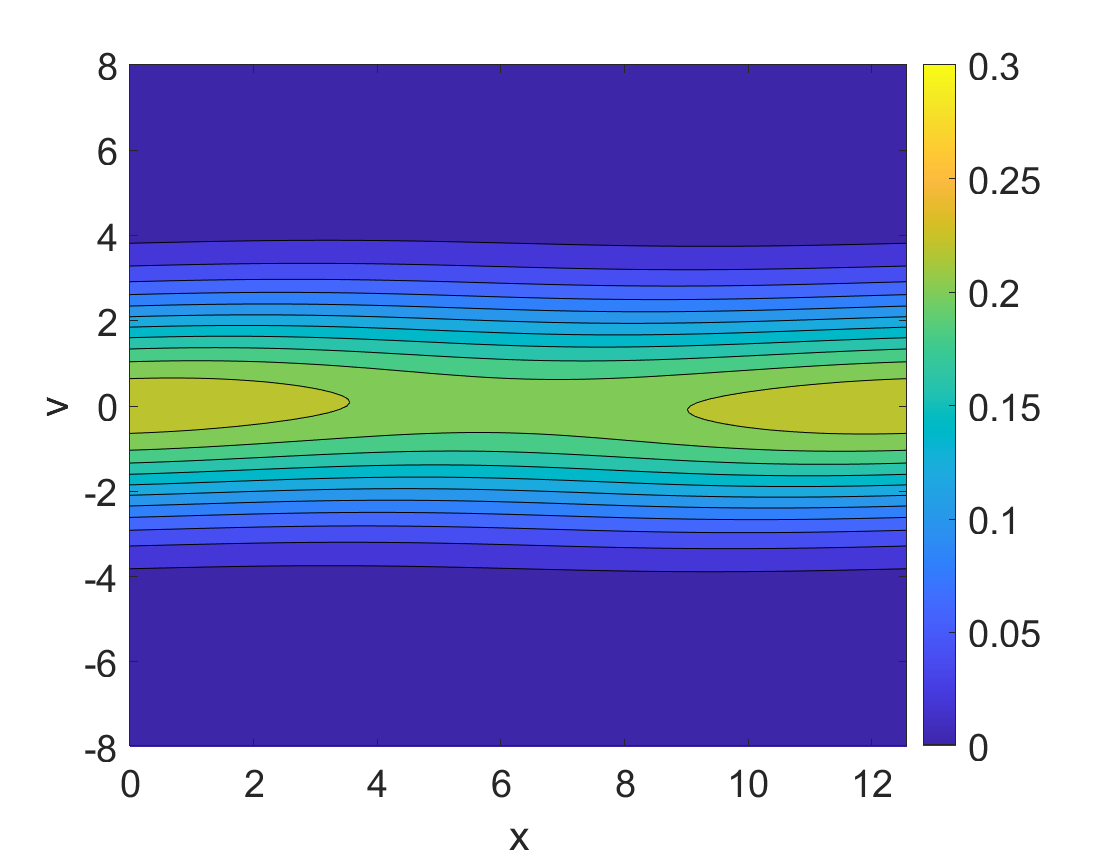}
        \caption{T = 0.5.}
        \label{fig4:subim2}
    \end{subfigure}
    \begin{subfigure}[t]{.45\textwidth}
        \centering
        \includegraphics[scale=.4]{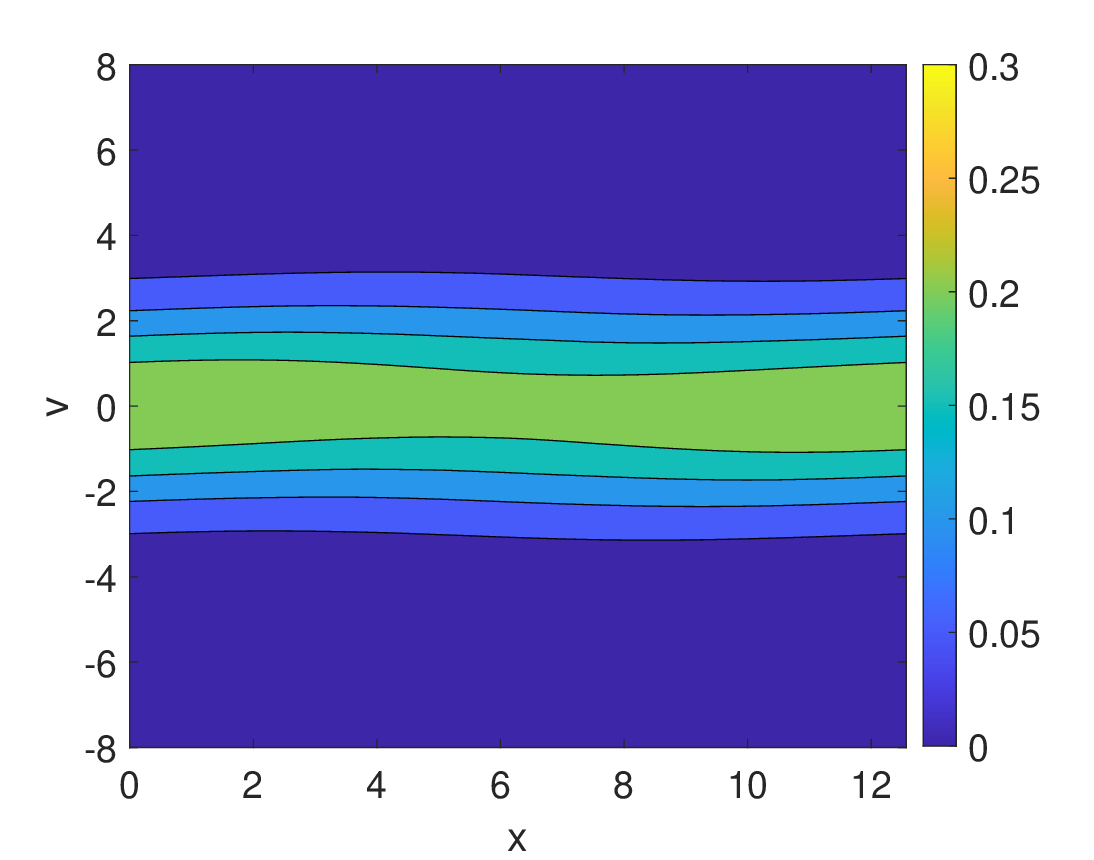}
        \caption{T = 1.}
        \label{fig4:subim3}
    \end{subfigure}
    \begin{subfigure}[t]{.45\textwidth}
        \centering
        \includegraphics[scale=.4]{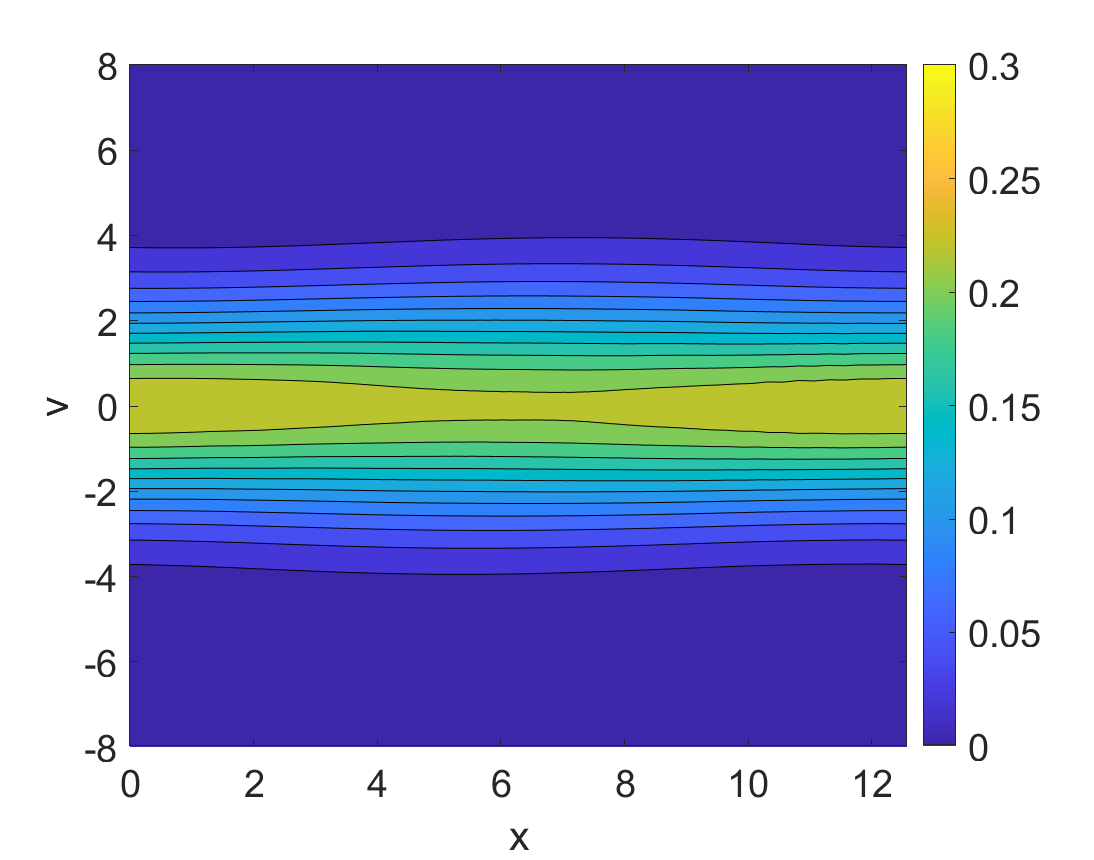}
        \caption{T = 2.}
        \label{fig4:subim4}
    \end{subfigure}
    \caption{Two-stream instability (Example 5.3) with collision frequency $\nu = 1$. Phase-space contour plots are shown at different times. The strong collisional effects rapidly suppress the instability generated by the interaction of the two plasma streams and prevent the formation of a singularity. As a result, filamentary structures are largely eliminated, and the plasma quickly relaxes toward equilibrium.}
    \label{fig4}
\end{figure}

\begin{figure}
    \centering
    \begin{subfigure}[t]{.45\textwidth}
        \centering
        \includegraphics[scale=.4]{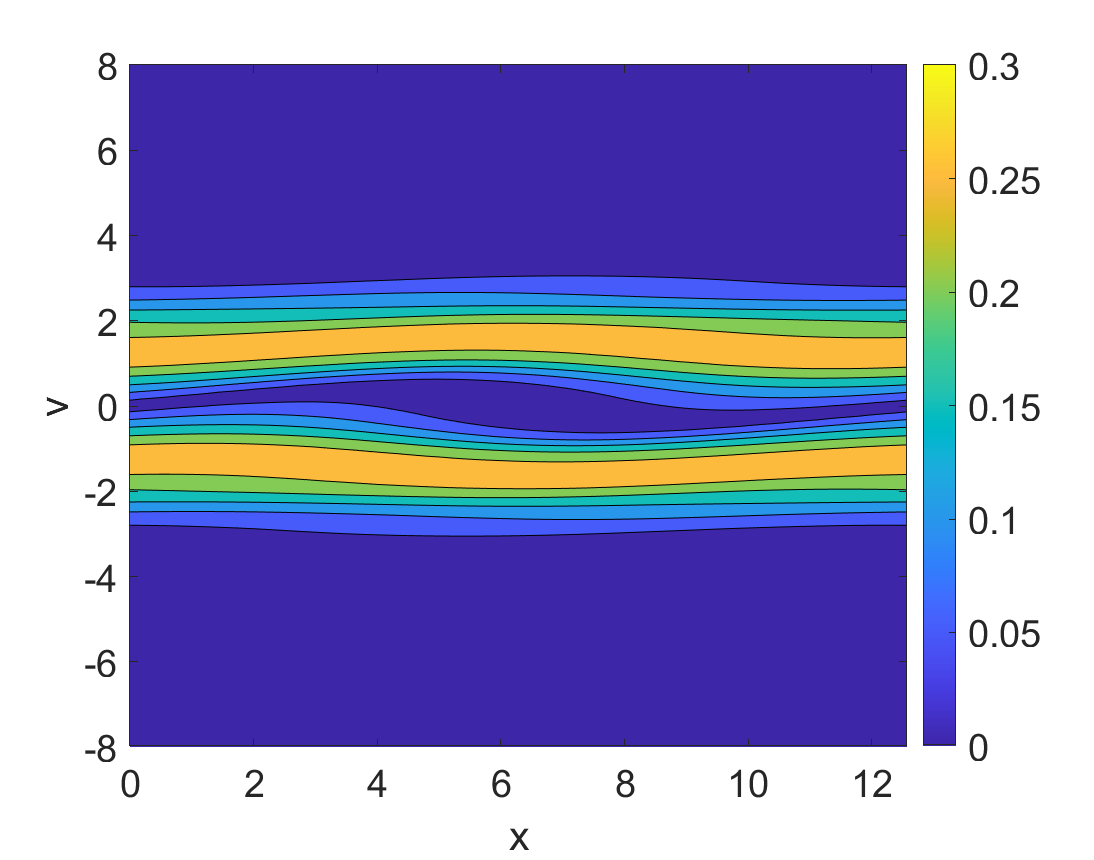}
        \caption{T = 10.}
        \label{fig5:subim1}
    \end{subfigure}
    \begin{subfigure}[t]{.45\textwidth}
        \centering
        \includegraphics[scale=.4]{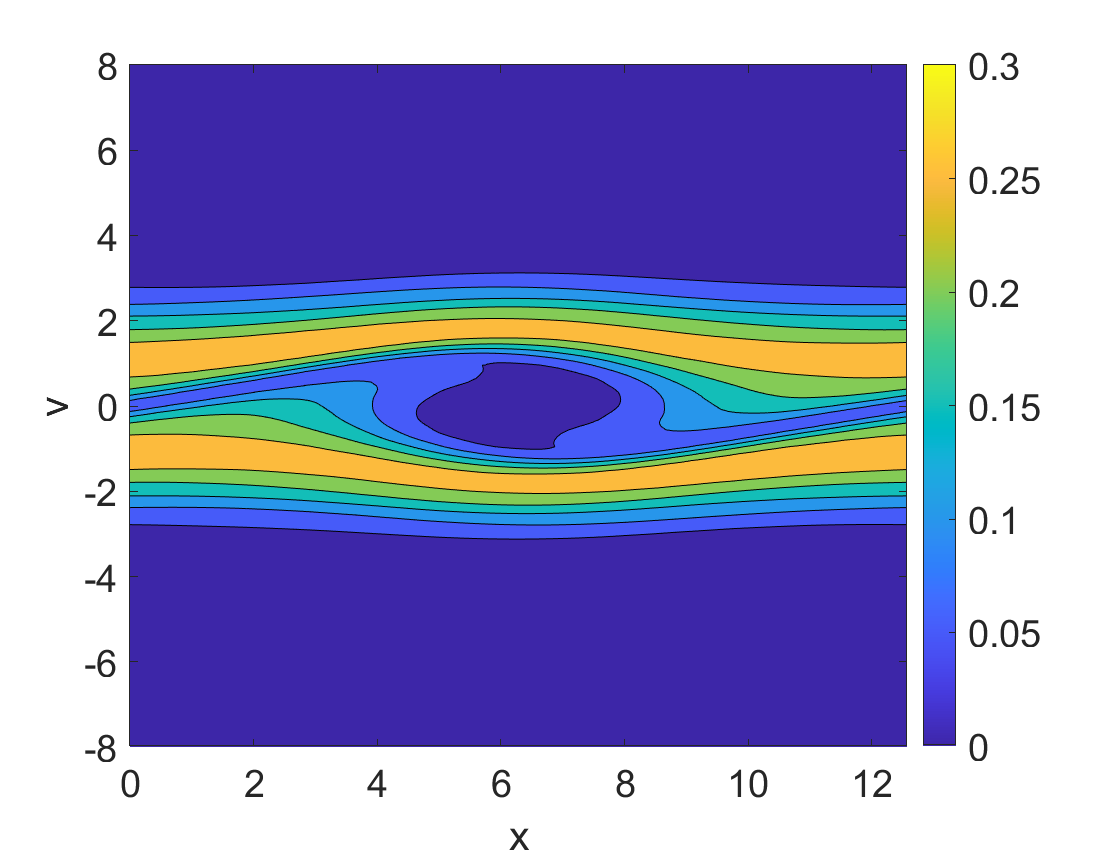}
        \caption{T = 15.}
        \label{fig5:subim2}
    \end{subfigure}
    \hfill
    \begin{subfigure}[t]{.45\textwidth}
        \centering
        \includegraphics[scale=.4]{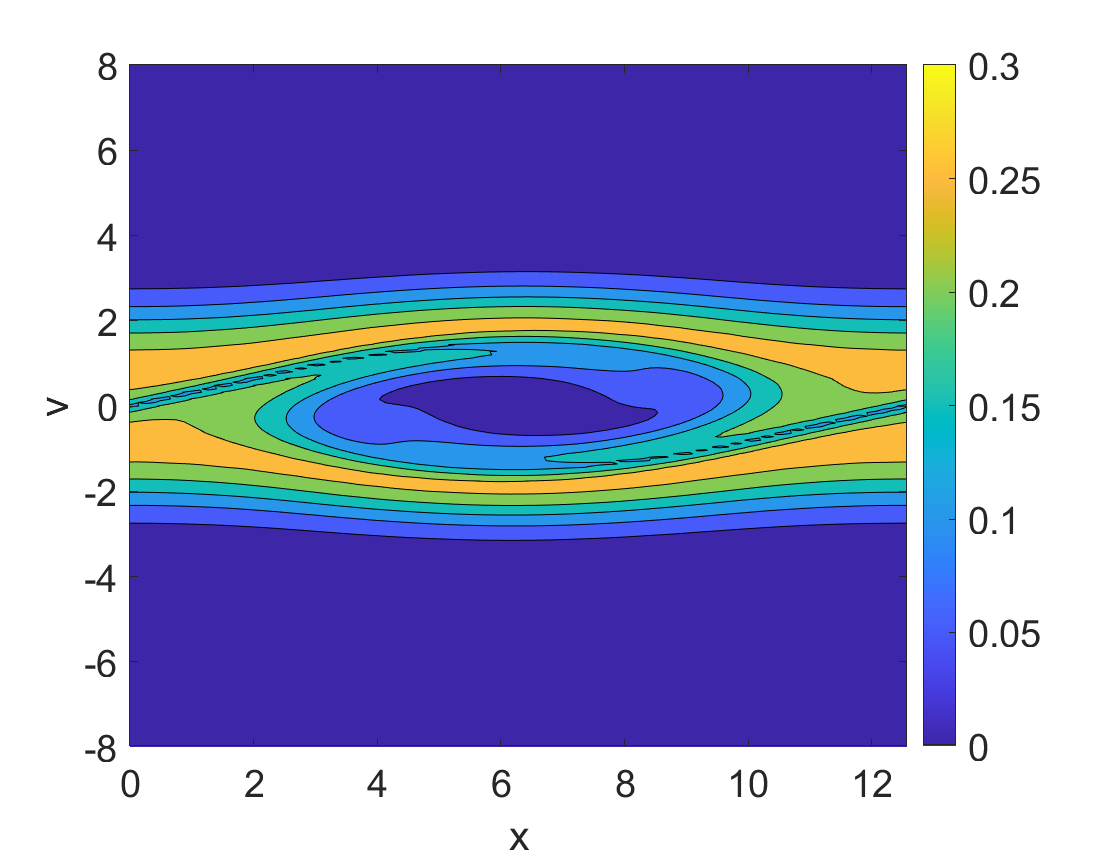}
        \caption{T = 20.}
        \label{fig5:subim3}
    \end{subfigure}
    \begin{subfigure}[t]{.45\textwidth}
        \centering
        \includegraphics[scale=.4]{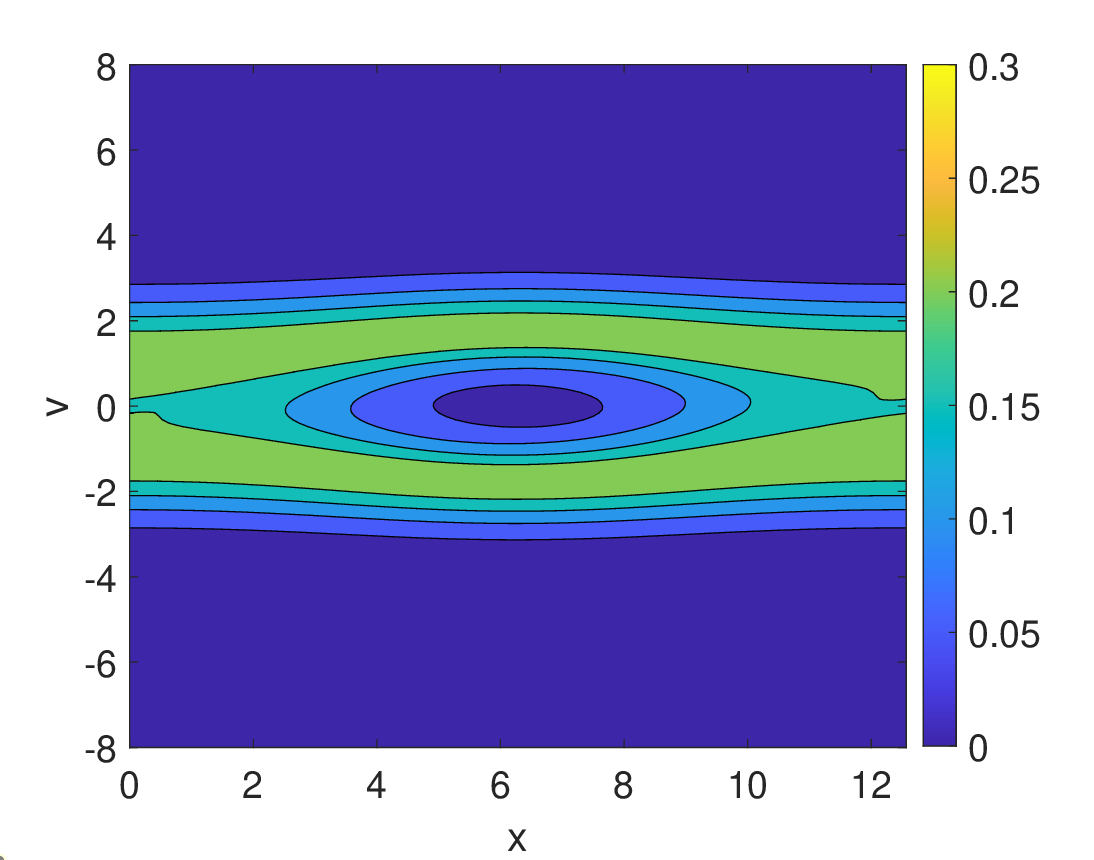}
        \caption{T = 40.}
        \label{fig5:subim4}
    \end{subfigure}
    \hfill
    \begin{subfigure}[t]{.45\textwidth}
        \centering
        \includegraphics[scale=.4]{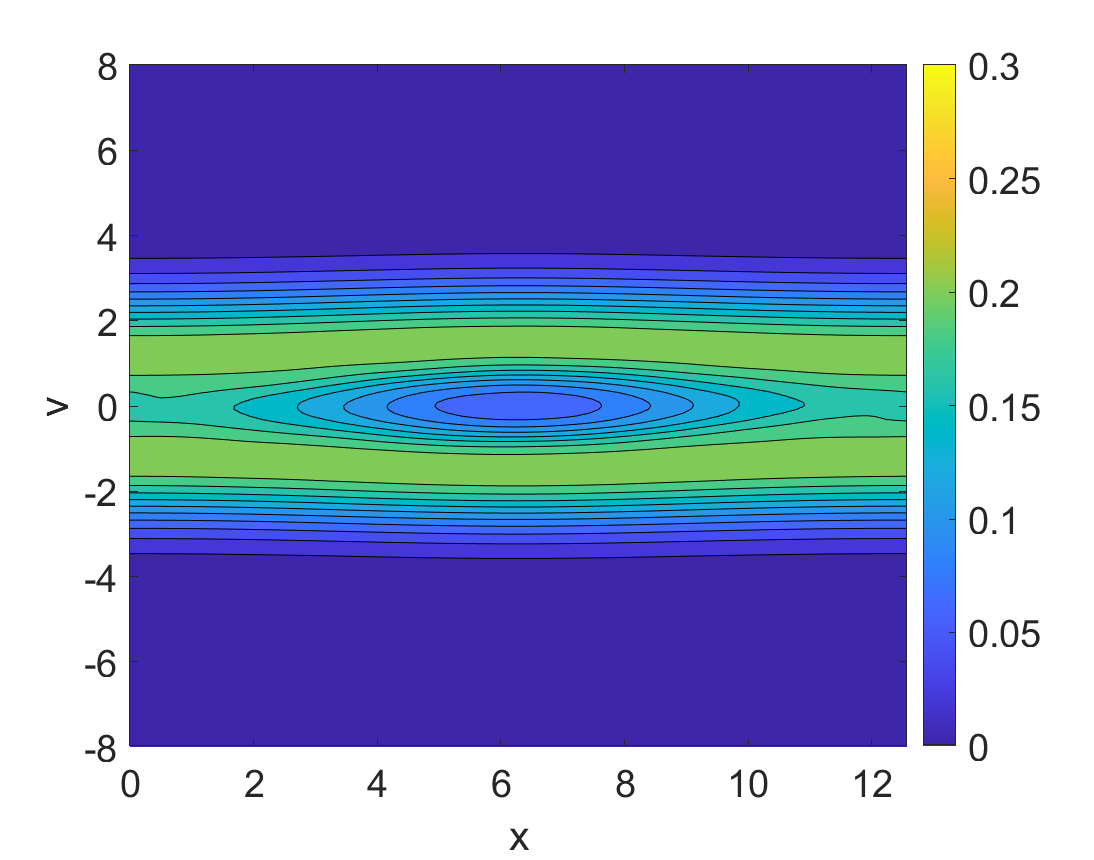}
        \caption{T = 80.}
        \label{fig5:subim5}
    \end{subfigure}
    \begin{subfigure}[t]{.45\textwidth}
    \centering
        \includegraphics[scale=.4]{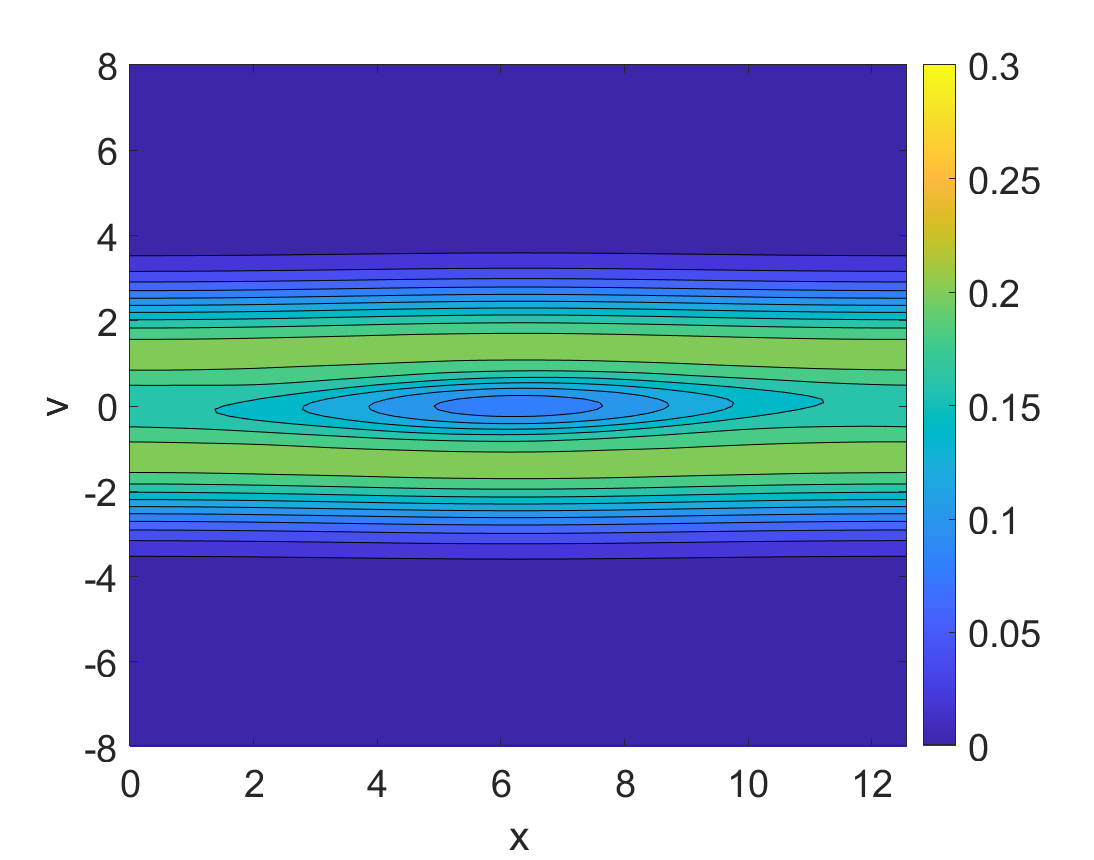}
        \caption{T = 100.}
        \label{fig5:subim6}
    \end{subfigure}
\caption{Two-stream instability (Example 5.3) with collision frequency $\nu = 0.001$. Phase-space contour plots are shown at different times. Since the collisional effects are weak, the instability generated by the interaction of the two plasma streams is only mildly suppressed. By $T = 20$, a singularity begins to form and develops filamentary structures that cause rank growth. However, as time progresses, the collisional effects smooth these structures and drive the plasma toward equilibrium.}
\label{fig5}
\end{figure}

\begin{figure}
    \centering
    \begin{subfigure}[t]{.45\textwidth}
        \centering
        \includegraphics[scale=0.4]{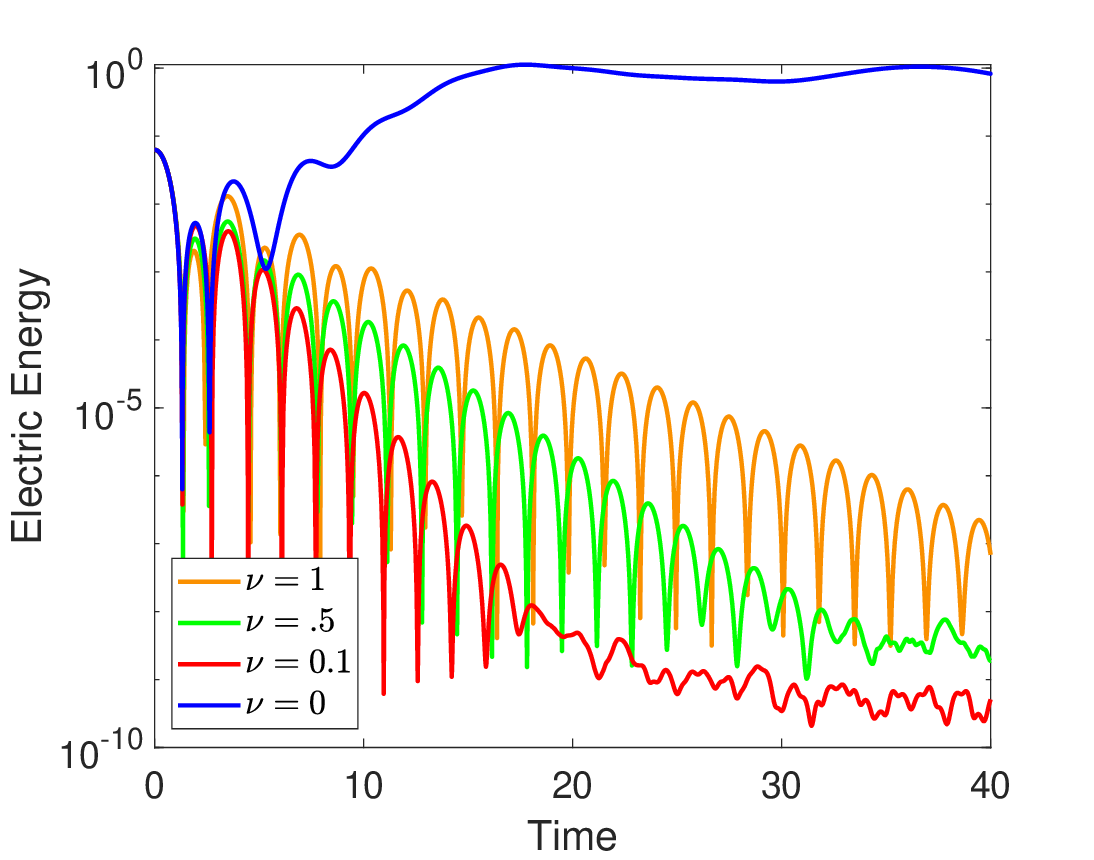}
        \caption{}
        \label{fig6:subim1}
    \end{subfigure}
    \begin{subfigure}[t]{0.45\textwidth}
        \centering
        \includegraphics[scale=0.4]{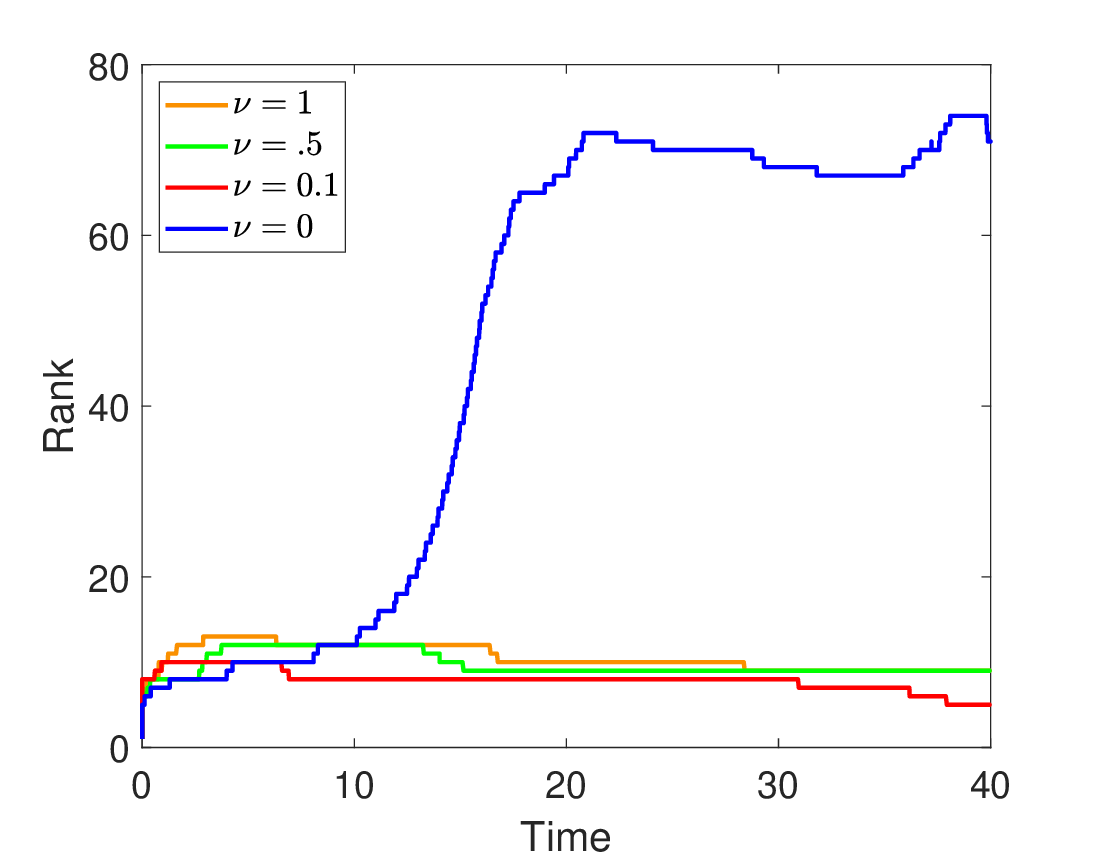}
        \caption{}
        \label{fig6:subim2}
    \end{subfigure}
    \hfill
    \begin{subfigure}[t]{0.45\textwidth}
        \centering
        \includegraphics[scale=0.4]{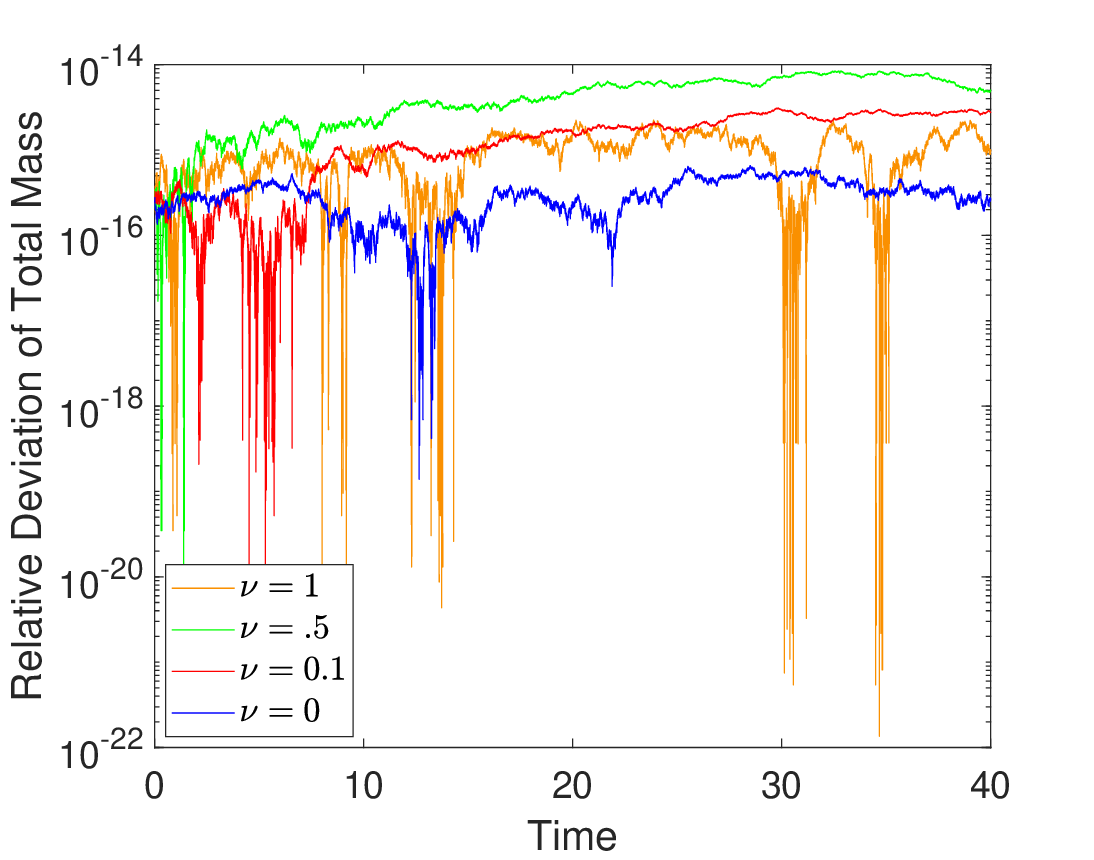}
        \caption{}
        \label{fig6:subim3}
    \end{subfigure}
    \begin{subfigure}[t]{0.45\textwidth}
        \centering
        \includegraphics[scale=0.4]{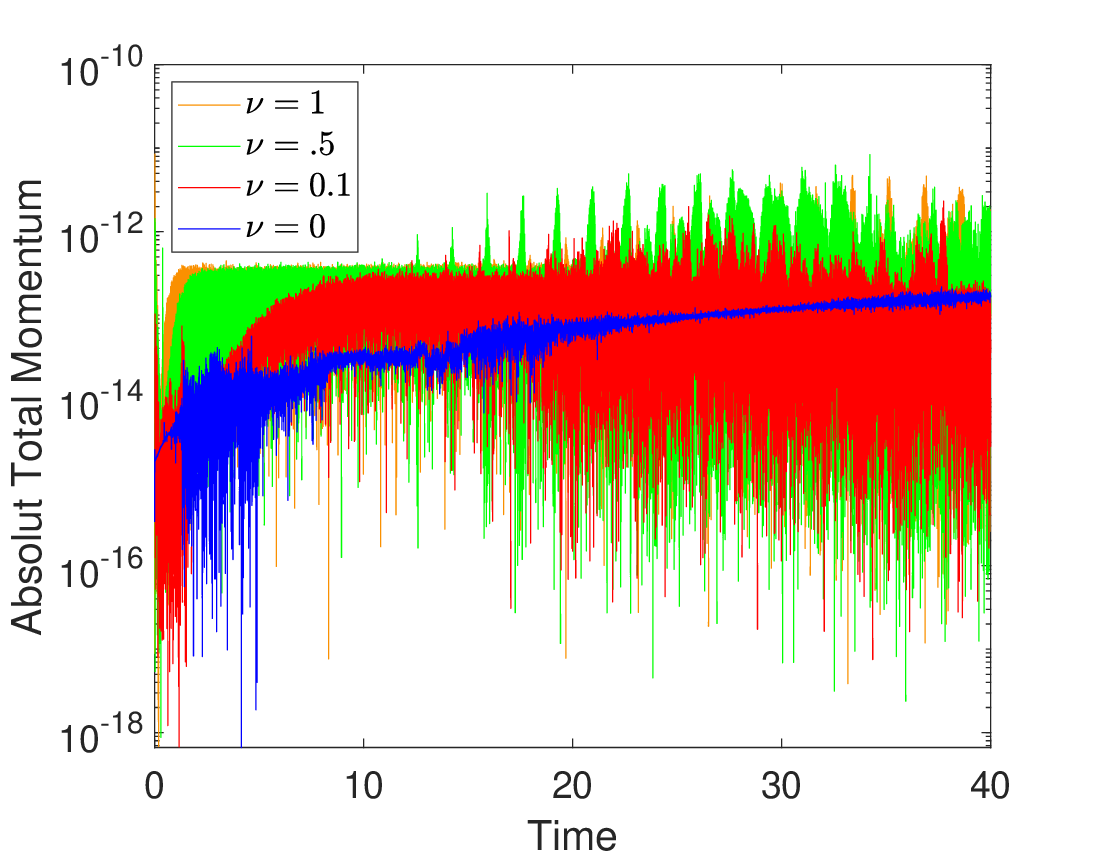}
        \caption{}
        \label{fig6:subim4}
    \end{subfigure}
    \hfill
    \begin{subfigure}[t]{0.45\textwidth}
        \centering
        \includegraphics[scale=0.4]{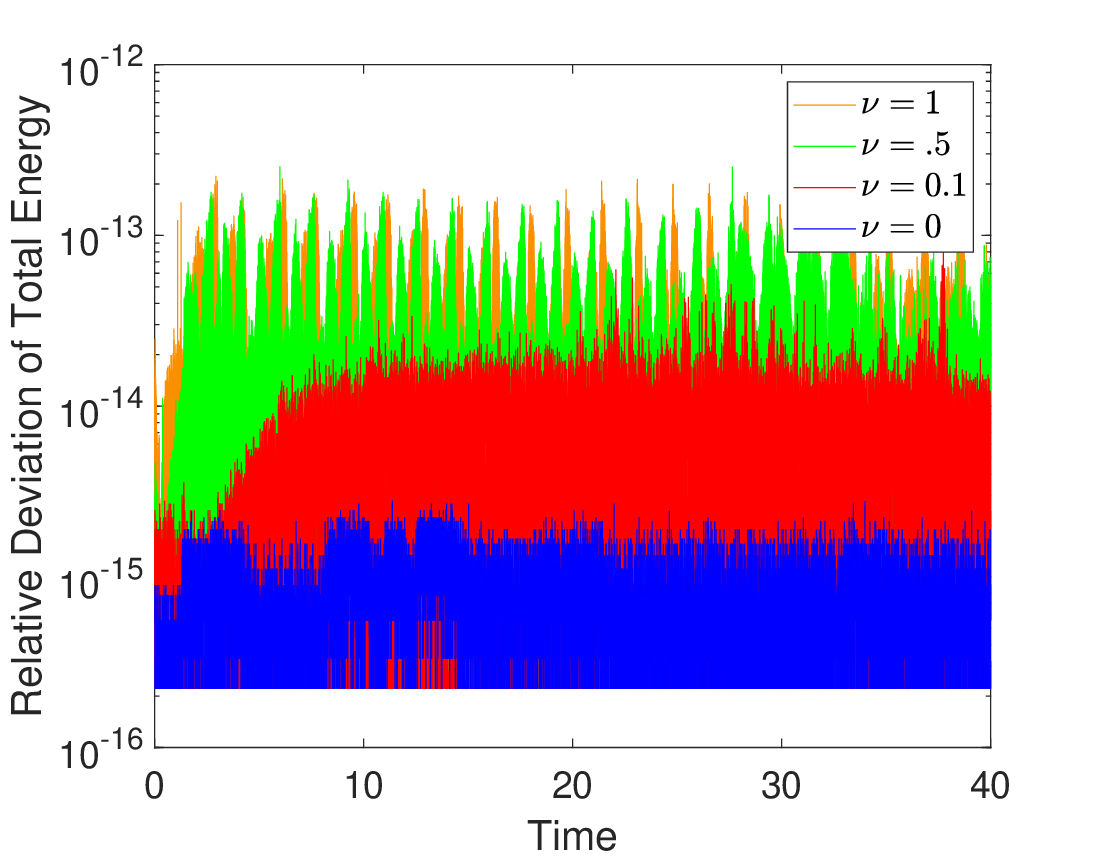}
        \caption{}
        \label{fig6:subim5}
    \end{subfigure}
    \hfill
    \caption{Two-stream instability (Example 5.3). Time histories of the electric energy and numerical rank for $\nu = 0, 0.1, 0.5,$ and $1$ are shown in (a)-(b). Increasing the collision frequency enhances the damping of the electric energy, almost mirroring the behavior observed in the Landau damping problems. Most notably, the collisional effects suppress the instability generated by the interaction of the two plasma streams, eliminating the singularity from forming and significantly limiting rank growth. The relative deviations of total mass and total energy, together with the absolute total momentum, are plotted in (c)-(e), confirming the conservation properties of the method.}
\end{figure}

    Here, we simulate the two-stream instability problem.
    \begin{gather*}
        f_0(x,v) = \frac{1}{\sqrt{2\pi}}v^2e^{-v^2/2}(1 + \alpha \text{cos}(\kappa x)), \hspace{3mm} x \in [0,L_x], \hspace{3mm} v \in [-L_v,L_v],
    \end{gather*}
    where $\alpha = .05,$ $\kappa = 0.5$, $L_x = 4\pi$, and $L_v = 8.$ For the truncation threshold, we use $\varepsilon = 10^{-4}$ for $\nu = 0$ and $\varepsilon = 10^{-6}$ for $\nu = 0.1, 0.5$ and $1$. In the collisionless case, a singularity would form between the two plasma streams, causing an instability as seen in \cite{guo2024LoMaC}. However, when collisions are introduced, this singularity is significantly suppressed. For sufficiently large collision frequency, e.g. $\nu \geq 0.1$, the singularity is entirely eliminated. These two cases are illustrated in Figures \hyperref[fig4]{4} and \hyperref[fig5]{5} for $\nu = 1$ and $\nu = 0.001.$ These results also agree with \cite{ye2024energy}. In Figures \hyperref[fig6:subim1]{6a-6b} we show the time history of the electric energy and numerical ranks. In Figures \hyperref[fig6:subim3]{6c-6e} we show that the method also conserves the mass, momentum, and energy of the system, as expected.

%%%%%%%%%%%%%%%%%%%%%%%%%%%%%%%%%%%%%%%%%
%%% 2-Stream nu = 1.
%%%%%%%%%%%%%%%%%%%%%%%%%%%%%%%%%%%%%%%%%

%%%%%%%%%%%%%%%%%%%%%%%%%%%%%%%%%%%%%%%%%
%%% 2-Stream nu = 0.001.
%%%%%%%%%%%%%%%%%%%%%%%%%%%%%%%%%%%%%%%%%

%%%%%%%%%%%%%%%%%%%%%%%%%%%%%%%%%%%%%%%%%%%%%%%%%%%%%%%%%%%%%%%%%%%%%%%%%%%%%%
%%% Bump-On-Tail Instability
%%%%%%%%%%%%%%%%%%%%%%%%%%%%%%%%%%%%%%%%%%%%%%%%%%%%%%%%%%%%%%%%%%%%%%%%%%%%%%
\begin{flushleft}
    \textbf{Example 5.4.} (Bump-On-Tail Instability.) 
\end{flushleft}
    
\begin{figure}
    \centering
    \begin{subfigure}[t]{.45\textwidth}
        \centering
        \includegraphics[scale=.4]{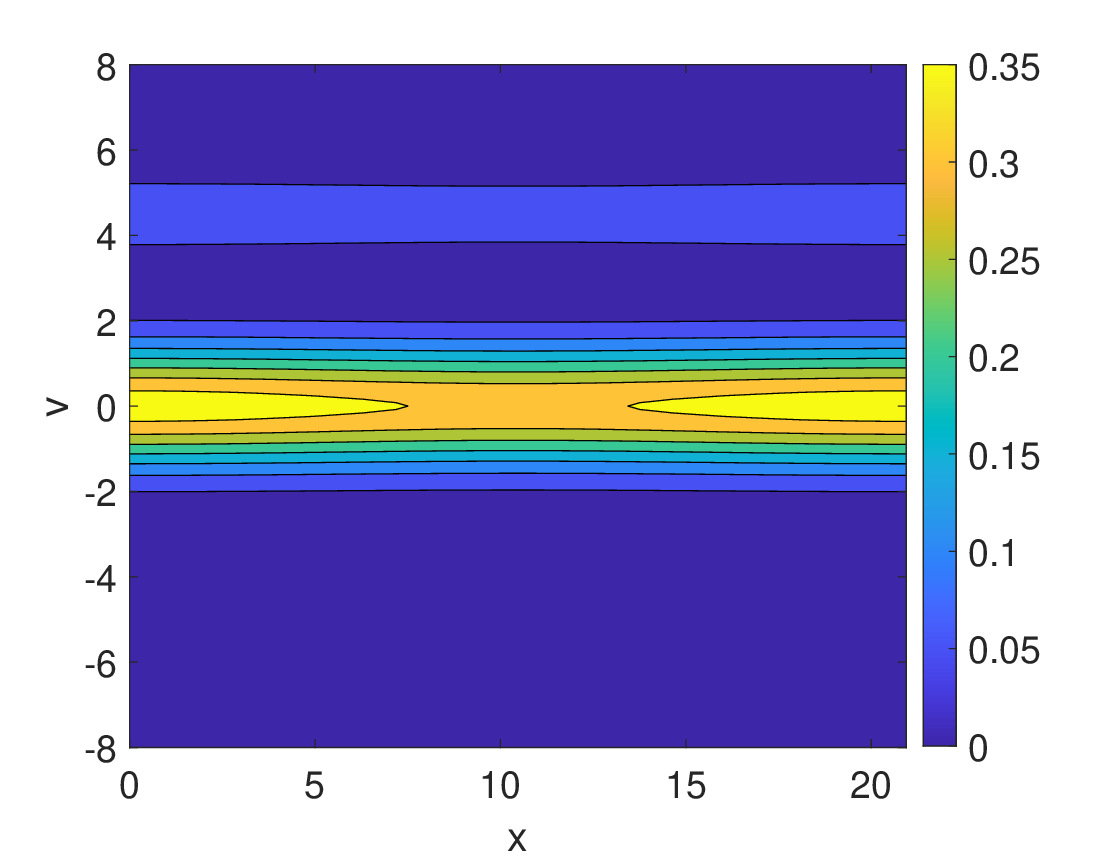}
        \caption{T = 0.}
        \label{fig7:subim1}
    \end{subfigure}
    \begin{subfigure}[t]{.45\textwidth}
        \centering
        \includegraphics[scale=.4]{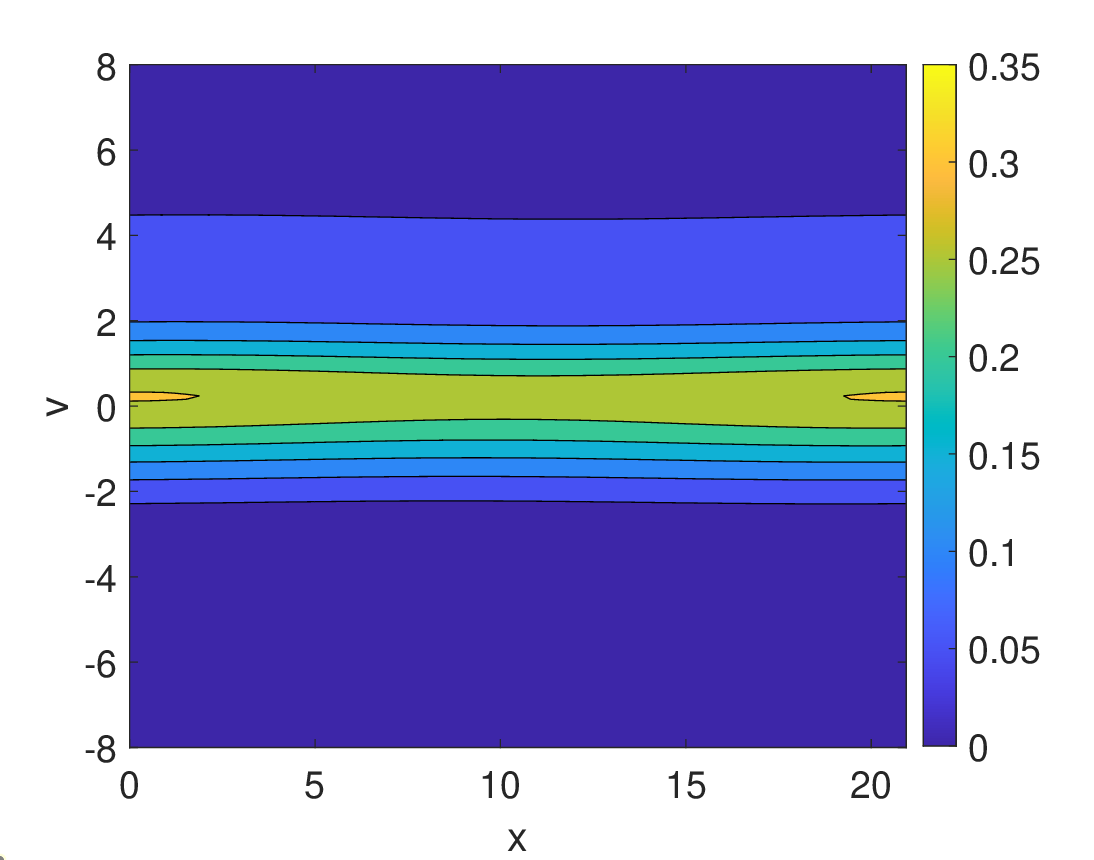}
        \caption{T = 0.1.}
        \label{fig7:subim2}
    \end{subfigure}
    
    \begin{subfigure}[t]{.45\textwidth}
        \centering
        \includegraphics[scale=.4]{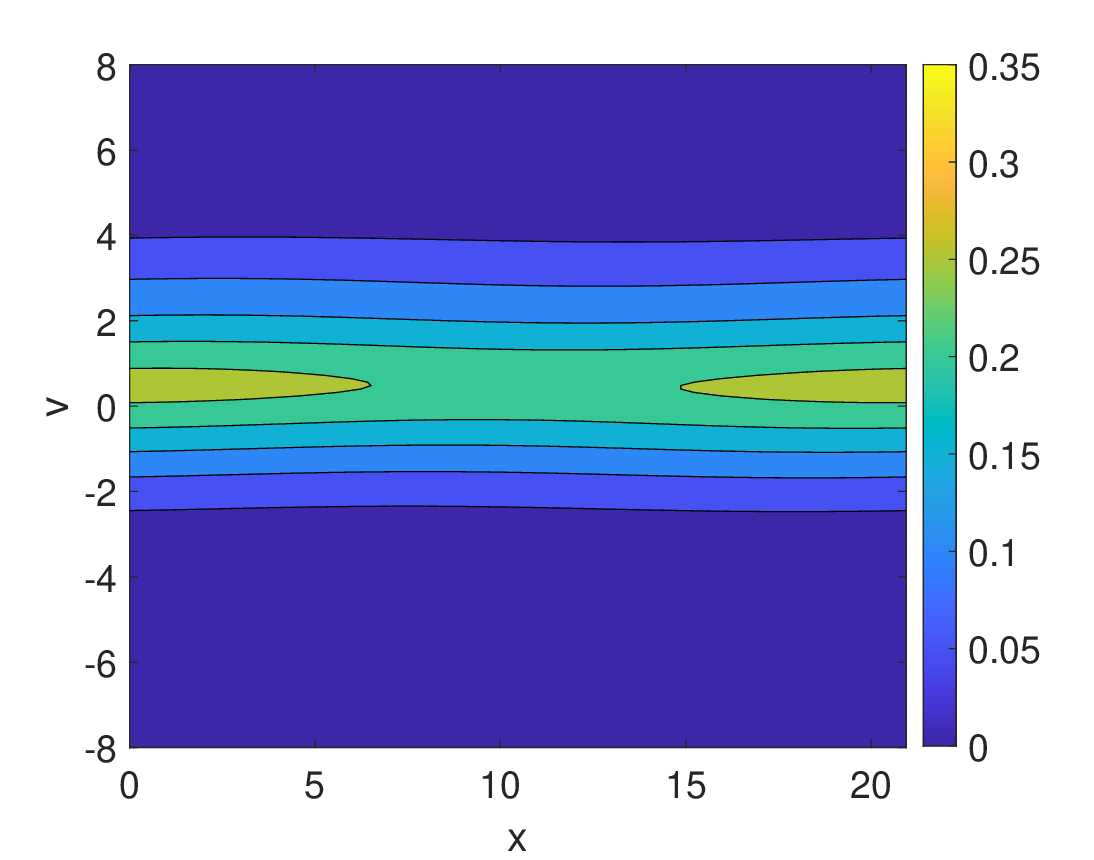}
        \caption{T = 0.25.}
        \label{fig7:subim3}
    \end{subfigure}
    \begin{subfigure}[t]{.45\textwidth}
        \centering
        \includegraphics[scale=.4]{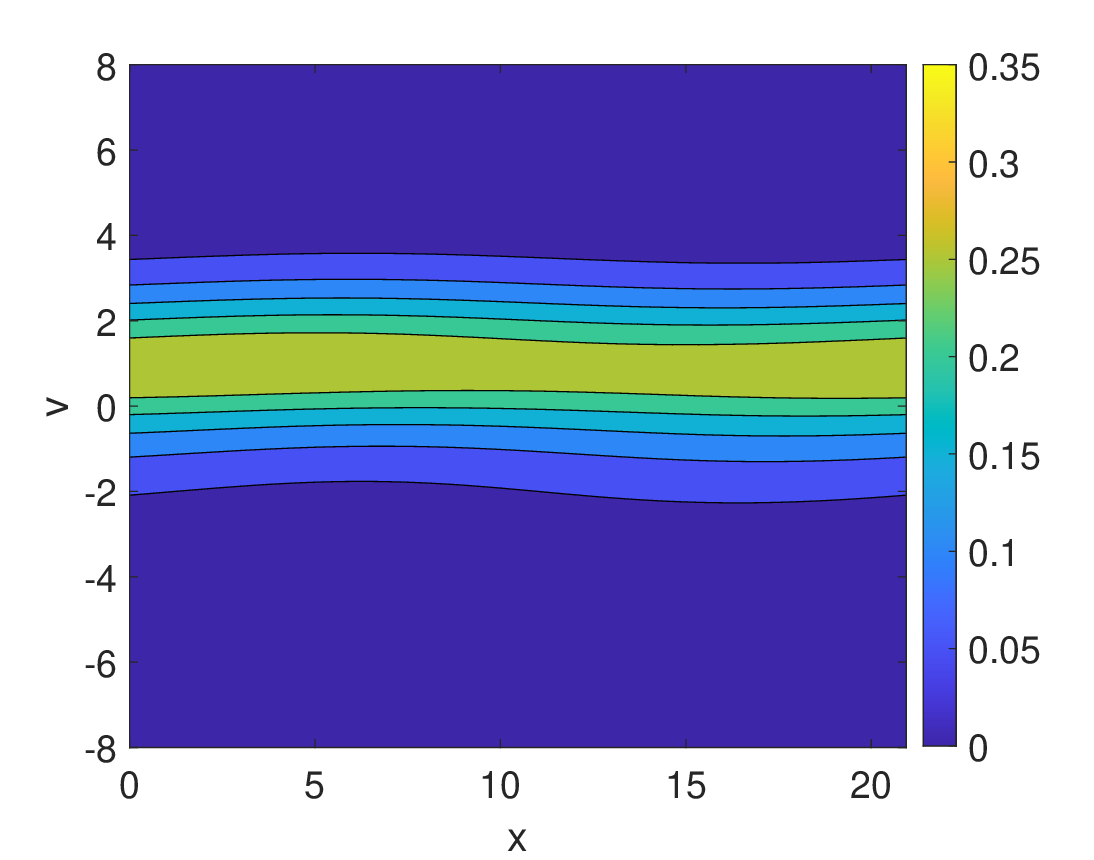}
        \caption{T = 1.}
        \label{fig7:subim4}
    \end{subfigure}
    \caption{Bump-on-tail instability (Example 5.4) with collision frequency $\nu = 1.$ Phase-space contour plots are shown at different times. The strong collisional effects rapidly suppress the instability associated with the high-energy tail and drive the plasma toward equilibrium. As a result, the filamentary structures are eliminated in a short amount of time.
    }
    \label{fig7}
\end{figure}

\begin{figure}
    \centering
    \begin{subfigure}[t]{.45\textwidth}
        \centering
        \includegraphics[scale=.4]{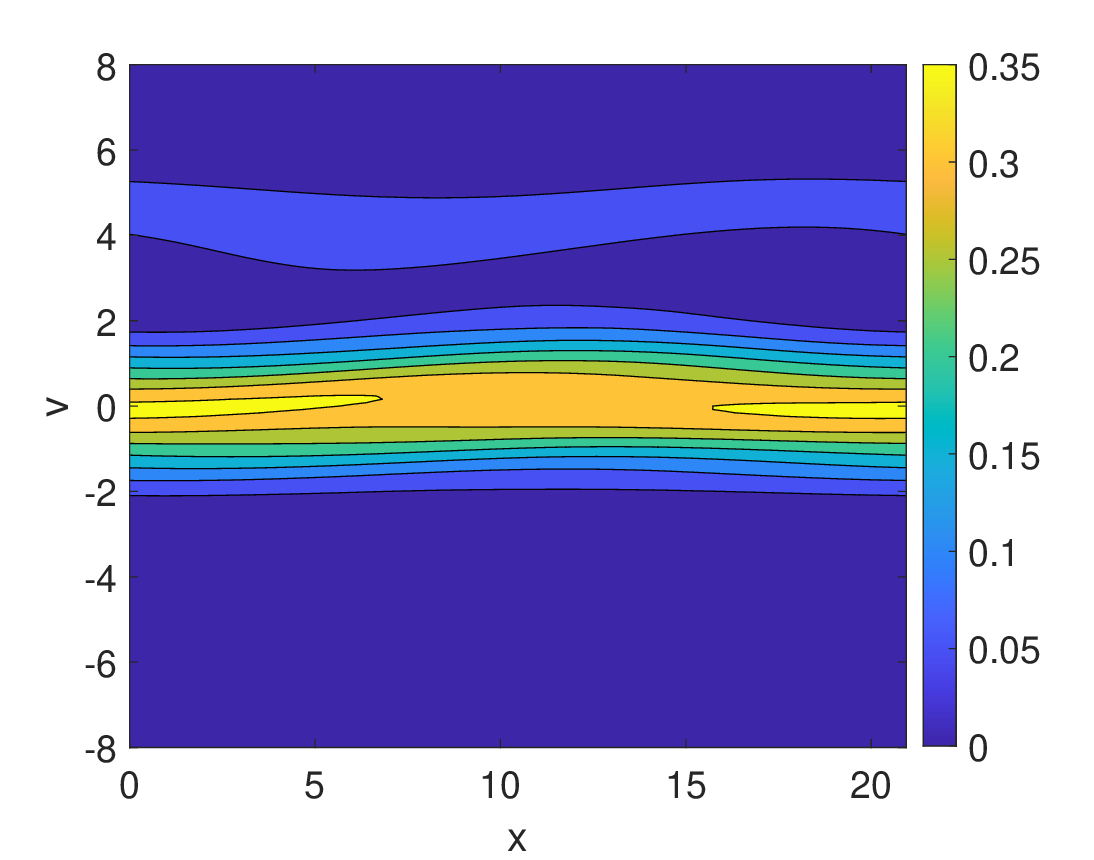}
        \caption{T = 10.}
        \label{fig8:subim1}
    \end{subfigure}
    \begin{subfigure}[t]{.45\textwidth}
    \centering
        \includegraphics[scale=.4]{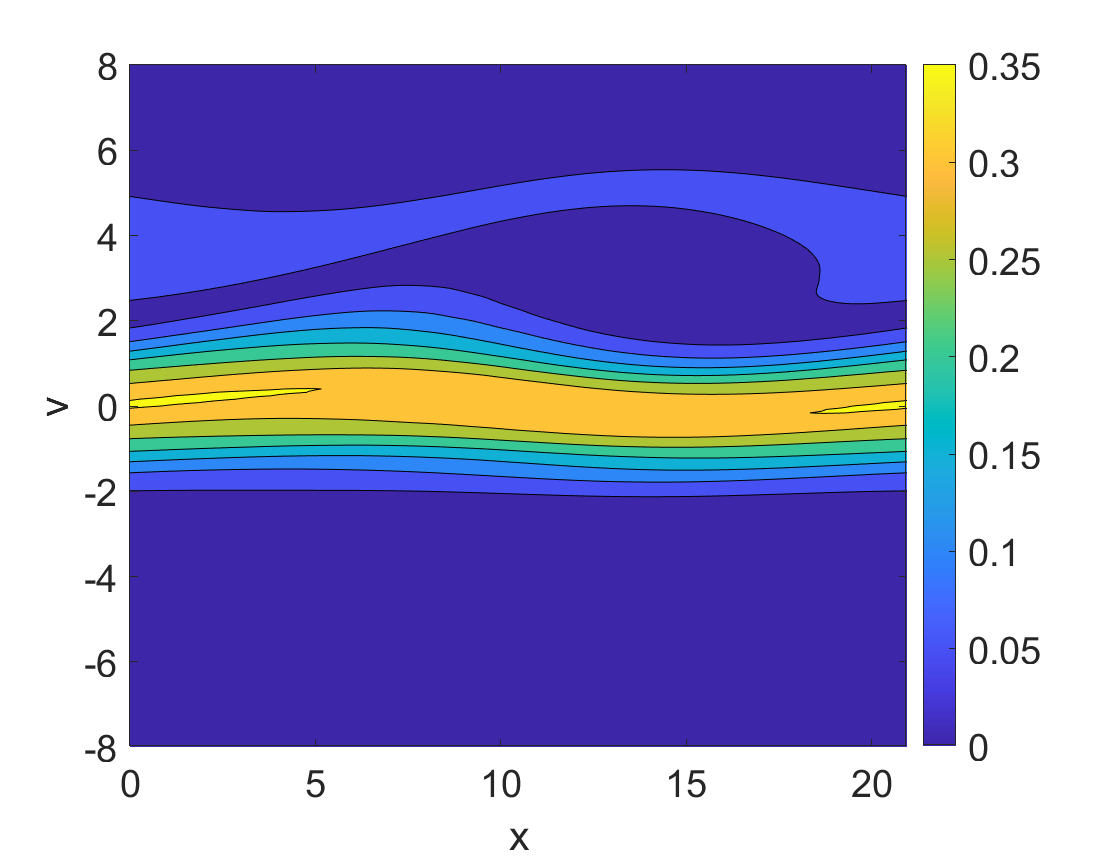}
        \caption{T = 15.}
        \label{fig8:subim2}
    \end{subfigure}
    \hfill
    \begin{subfigure}[t]{.45\textwidth}
    \centering
        \includegraphics[scale=.4]{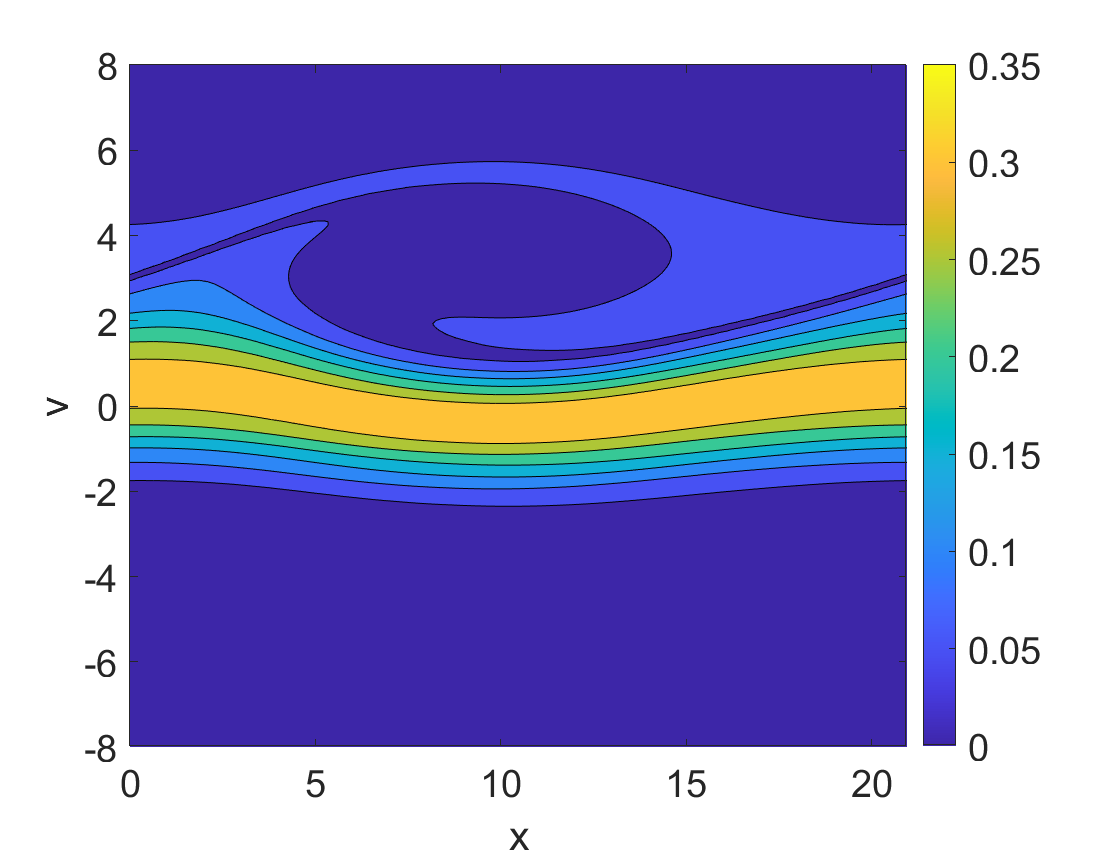}
        \caption{T = 20.}
        \label{fig8:subim3}
    \end{subfigure}
    \begin{subfigure}[t]{.45\textwidth}
    \centering
        \includegraphics[scale=.4]{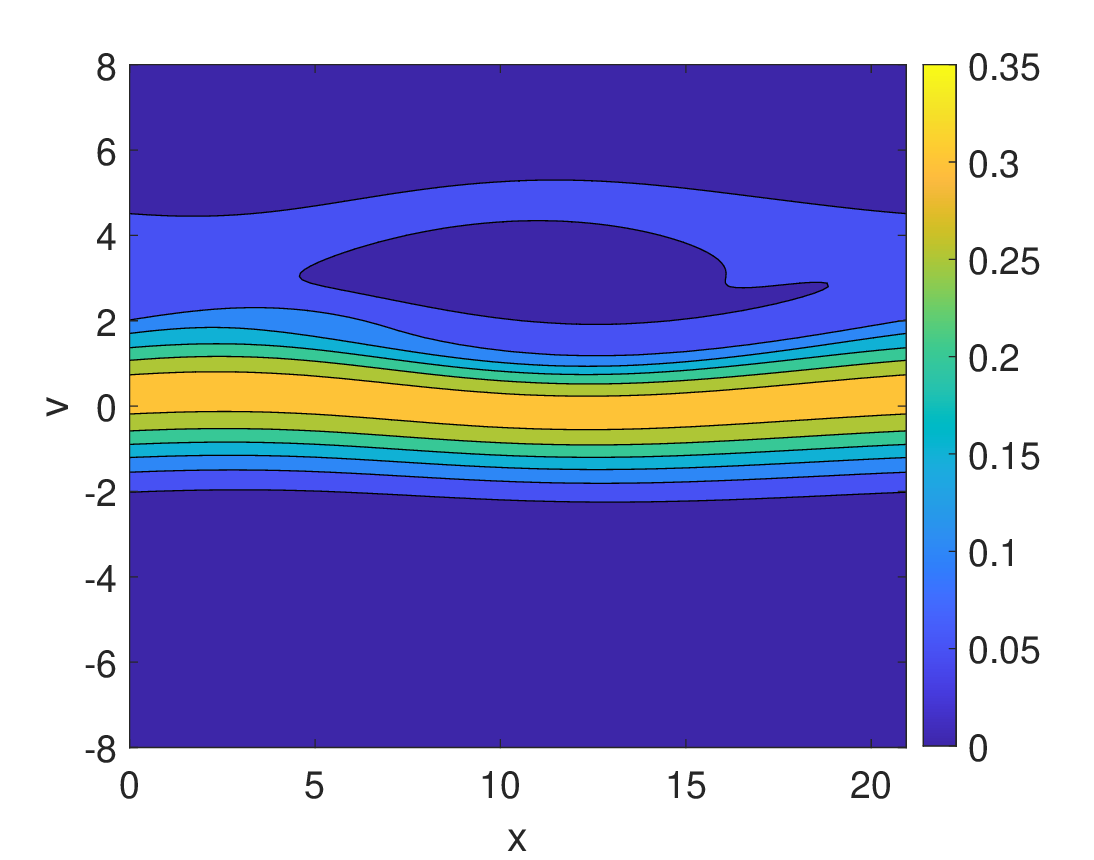}
        \caption{T = 40.}
        \label{fig8:subim4}
    \end{subfigure}
    \hfill
    \begin{subfigure}[t]{.45\textwidth}
    \centering
        \includegraphics[scale=.4]{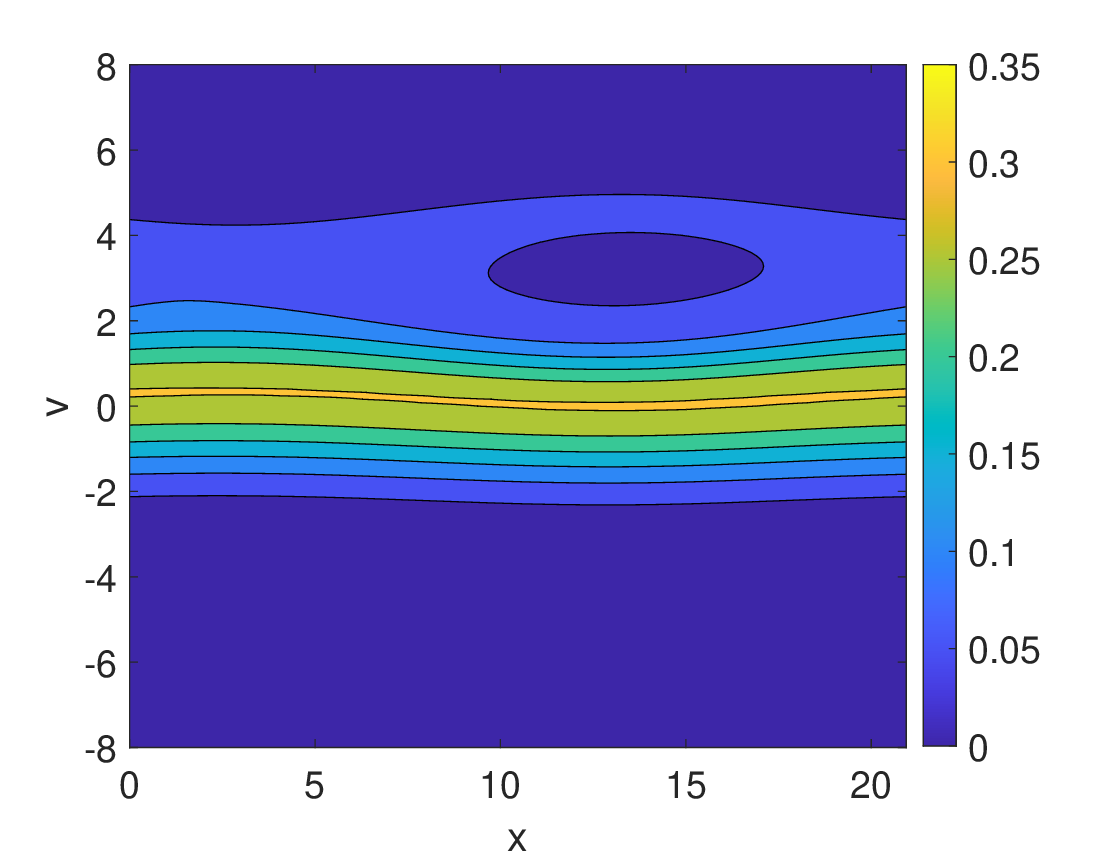}
        \caption{T = 80.}
        \label{fig8:subim5}
    \end{subfigure}
    \begin{subfigure}[t]{.45\textwidth}
    \centering
        \includegraphics[scale=.4]{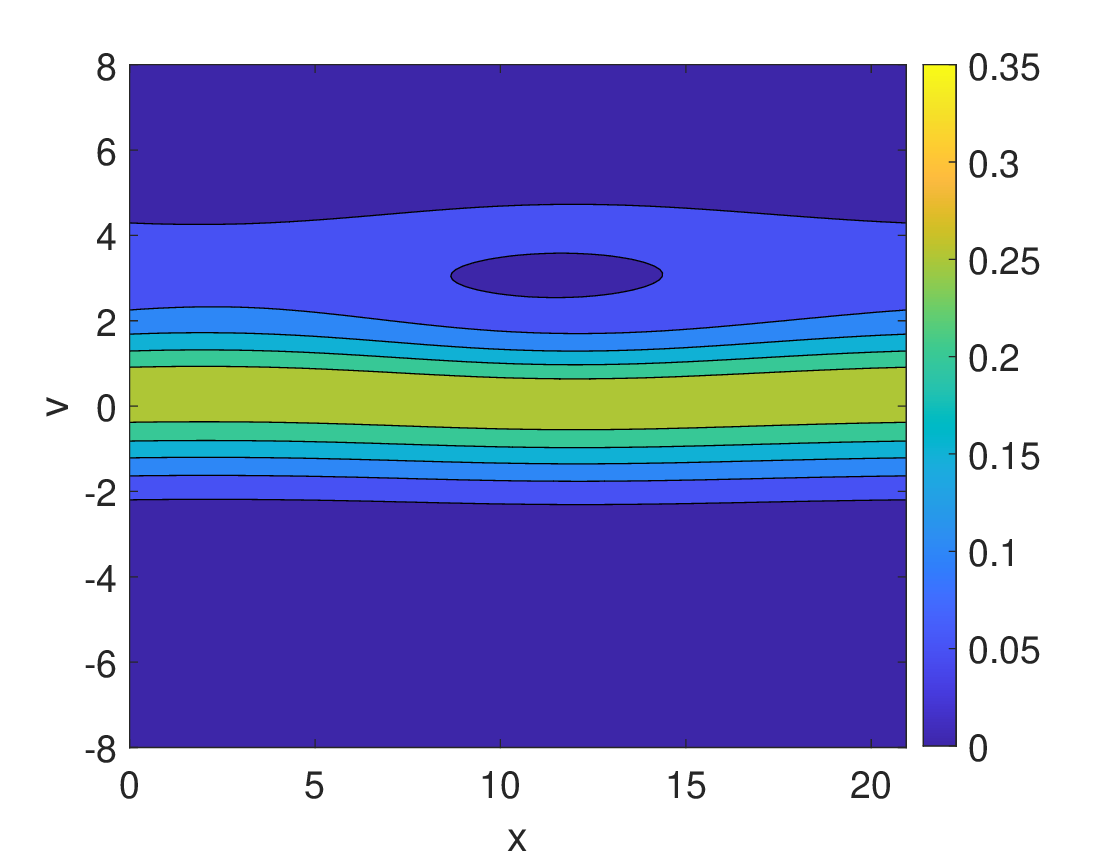}
        \caption{T = 100.}
        \label{fig8:subim6}
    \end{subfigure}
    \caption{Bump-on-tail instability (Example 5.4) with collision frequency of $\nu = .001.$ Phase-space contour plots are shown at different times.  Since the collisional effects are weak, the instability associated with the high-energy tail is only mildly suppressed. At time T = 20 (c), filamentary structures emerge as a results of the instability. As time progresses, collisional effects gradually drive the plasma toward equilibrium, leading to a slow merging of the high-energy tail with the background plasma distribution.}
    \label{fig8}
\end{figure}

\begin{figure}
    \centering
    \begin{subfigure}[t]{.45\textwidth}
        \centering
        \includegraphics[scale=0.4]{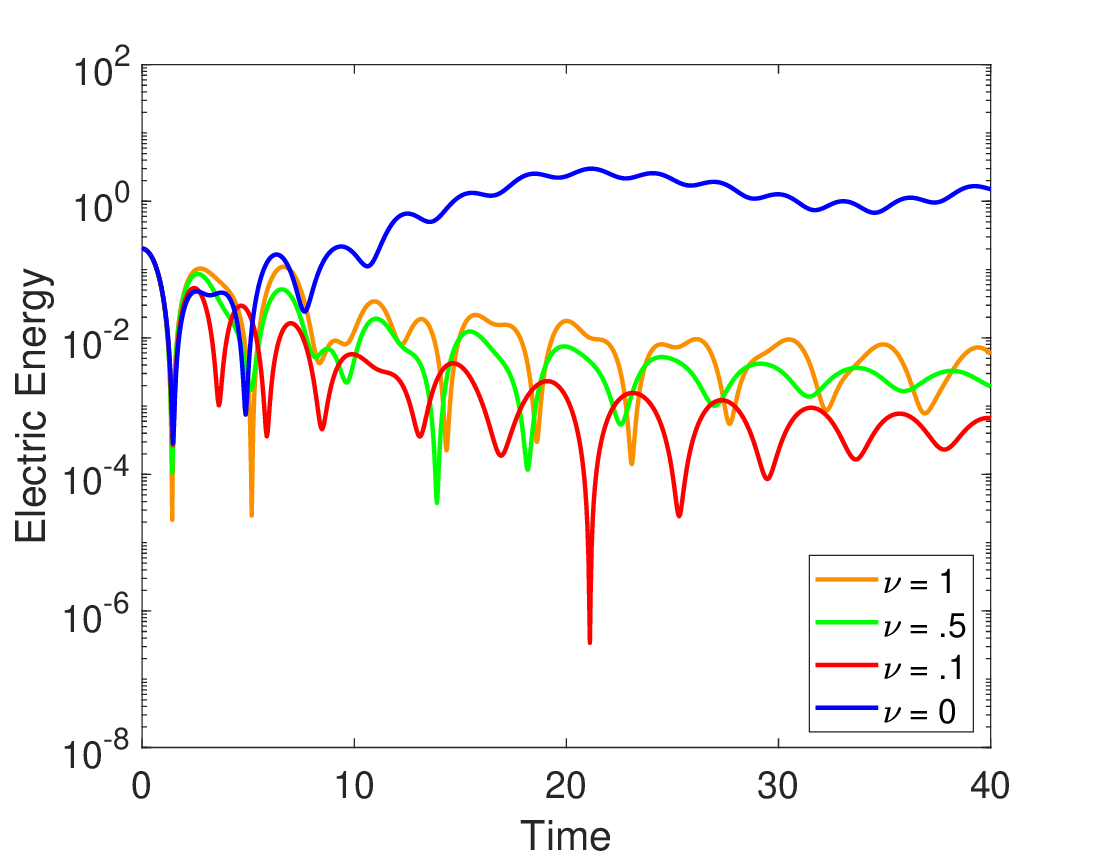}
        \caption{}
        \label{fig9:subim1}
    \end{subfigure}
    \begin{subfigure}[t]{0.45\textwidth}
        \centering
        \includegraphics[scale=0.4]{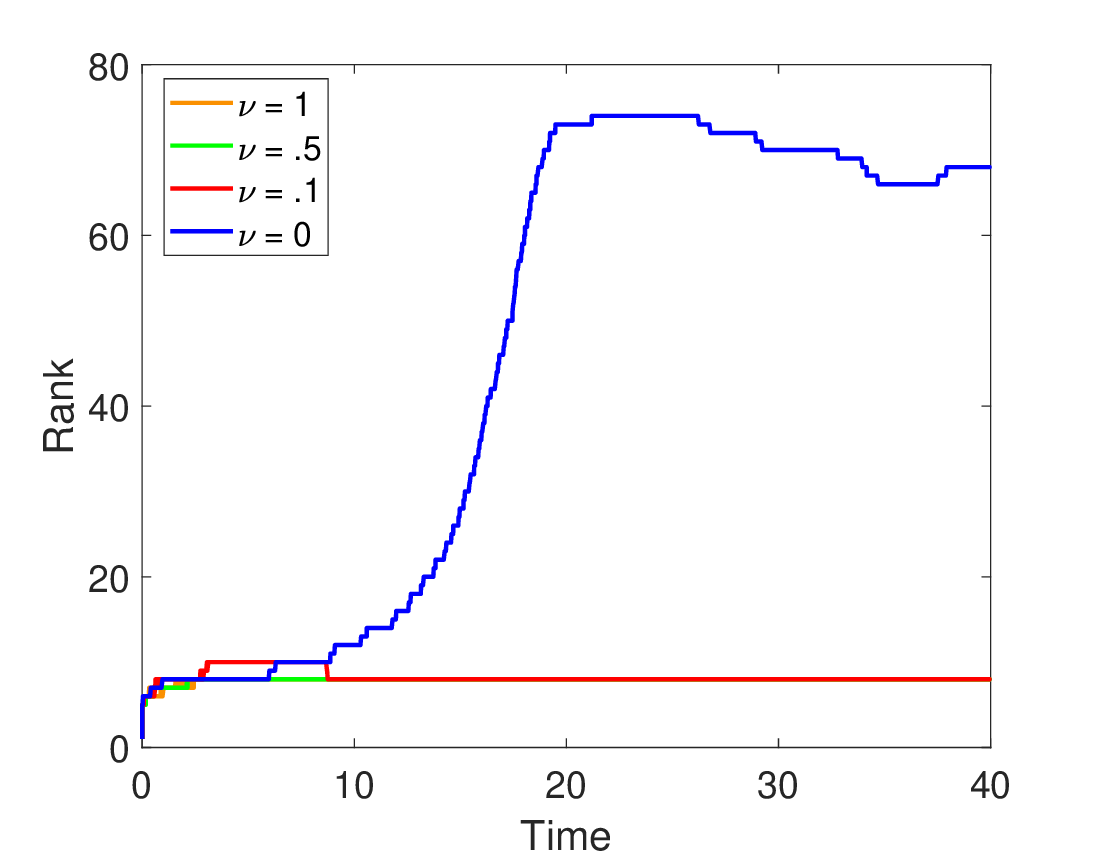}
        \caption{}
        \label{fig9:subim2}
    \end{subfigure}
    \hfill
    \begin{subfigure}[t]{0.45\textwidth}
        \centering
        \includegraphics[scale=0.4]{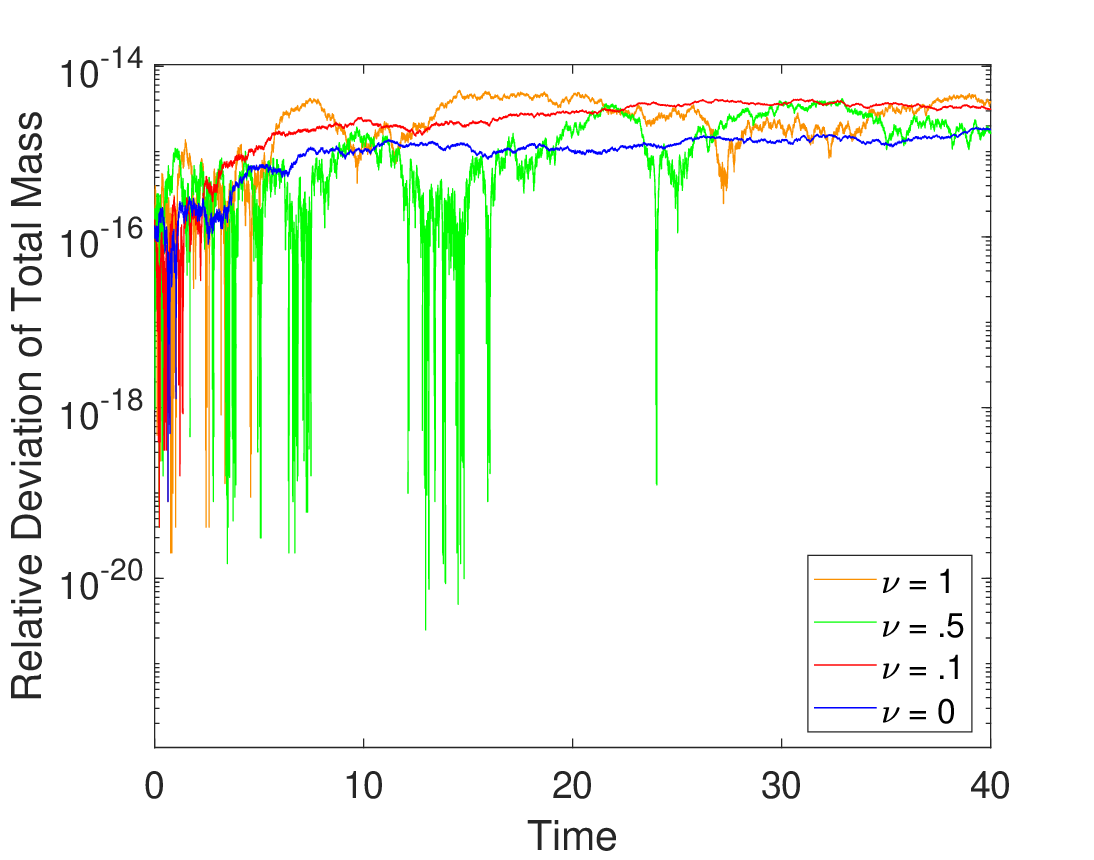}
        \caption{}
        \label{fig9:subim3}
    \end{subfigure}
    \begin{subfigure}[t]{0.45\textwidth}
        \centering
        \includegraphics[scale=0.4]{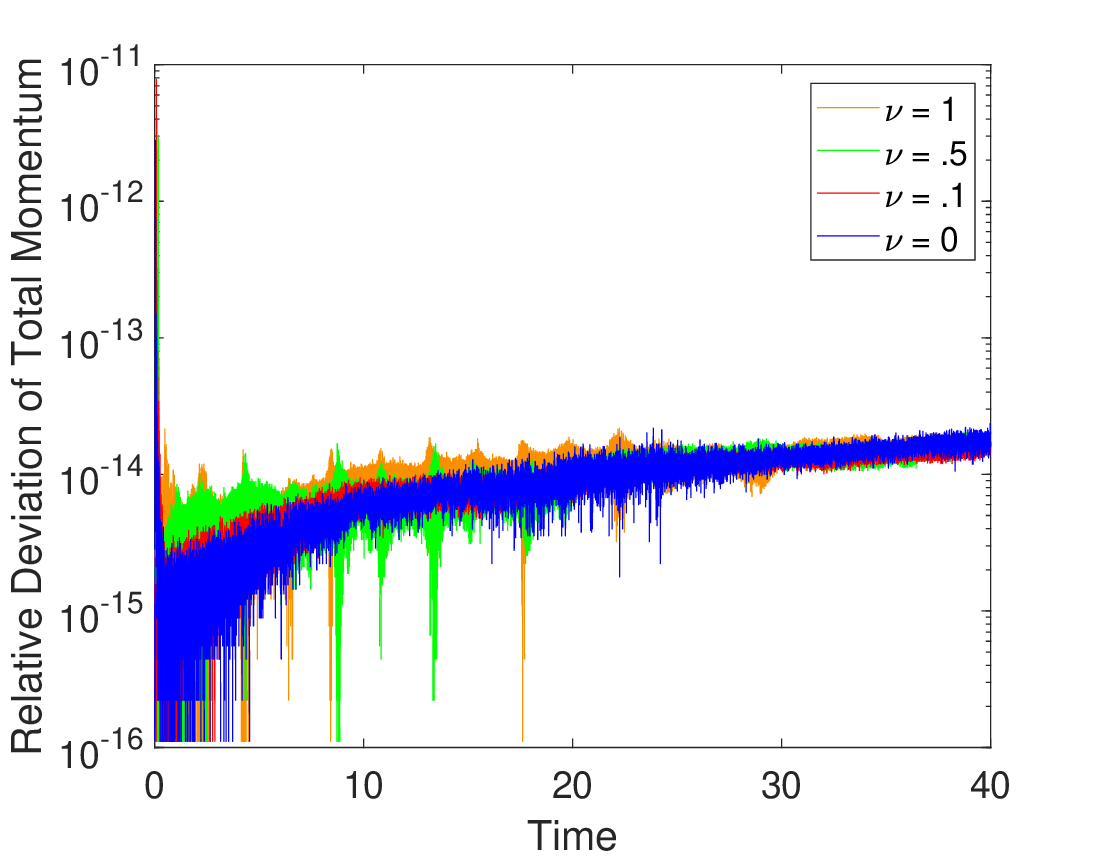}
        \caption{}
        \label{fig9:subim4}
    \end{subfigure}
    \hfill
    \begin{subfigure}[t]{0.45\textwidth}
        \centering
        \includegraphics[scale=0.4]{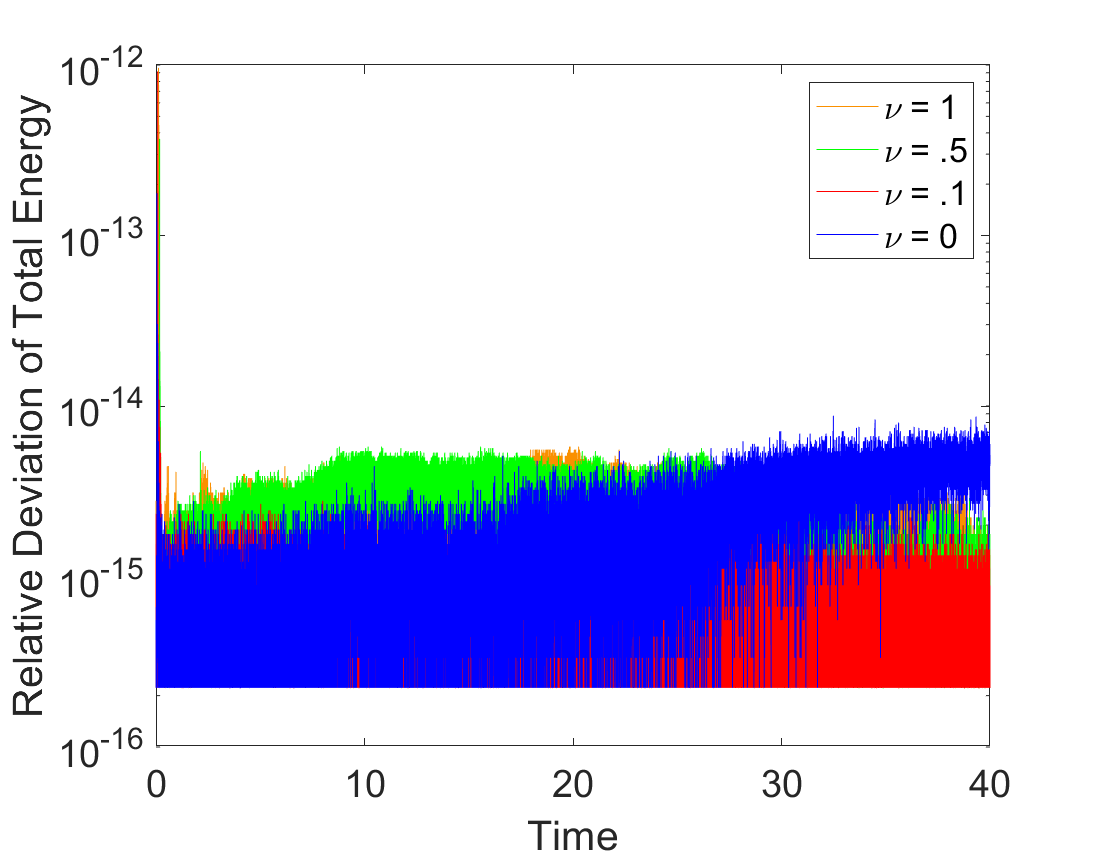}
        \caption{}
        \label{fig9:subim5}
    \end{subfigure}
    \hfill
    \caption{Bump-on-tail instability (Example 5.4). Time histories of the electric energy and numerical rank for $\nu = 0,0.1,0.5,$ and $1$ are shown in (a)-(b). Increasing the collision frequency enhances the damping of the electric energy. Most notably, the included collisional effects greatly suppress the bump-on-tail instability associated with the high-energy tail and limit rank growth. The relative deviations of total mass, momentum, and energy are plotted in (c)-(e), confirming the conservation properties of the method.
    \label{fig9}}
\end{figure}

    Lastly, we simulate the bump-on-tail instability problem.
    \begin{gather*}
        f(x,v,t=0) = (1 + \alpha \text{cos}(\kappa x))\left (  n_p \text{exp}\left (-\frac{v^2}{2}\right ) + n_b \text{exp}\left (- \frac{(v - u)^2}{2v_t} \right )  \right ),
    \end{gather*}
    where $\alpha = 0.04$, $\kappa = 0.3$, $n_p = \frac{9}{10\sqrt{2\pi}},$ and $n_p = \frac{2}{10\sqrt{2\pi}}$, $u = 4.5$, $v_t = 0.5$. The domain is set to be $[0,L_x]\times [-L_v,L_v]$ with $L_x = 2\pi /k$ and $L_v = 8.$ The truncation threshold is set to $\epsilon = 10^{-5}$. Similar to the strong Landau damping test, the filamentation greatly influences the plasma in the collisionless case. Once collisions are added, this influence greatly diminishes as the higher velocity electrons are pushed towards equilibrium. And, similar to the two-stream test, for sufficiently large collision frequency, the instability from the filamentation is entirely eliminated. This is illustrated in Figures \hyperref[fig7]{7} and \hyperref[fig8]{8} for $\nu = 1$ and $\nu = 0.001.$ In Figures \hyperref[fig9:subim1]{9a-9b} we show the time history of the electric energy and the numerical ranks. Lastly, we show the low-rank DG solutions and conservation of mass, momentum, and energy of the system in Figures \hyperref[fig9:subim3]{9c-9e}.

%%%%%%%%%%%%%%%%%%%%%%%%%%%%%%%%%%%%%%%%%
%%% Bump-On-Tail nu = 1.
%%%%%%%%%%%%%%%%%%%%%%%%%%%%%%%%%%%%%%%%%

%%%%%%%%%%%%%%%%%%%%%%%%%%%%%%%%%%%%%%%%%
%%% Bump-On-Tail nu = 0.001.
%%%%%%%%%%%%%%%%%%%%%%%%%%%%%%%%%%%%%%%%%

%%%%%%%%%%%%%%%%%%%%%%%%%%%%%%%%%%%%%%%%%
%%% Bump-On-Tail Stats
%%%%%%%%%%%%%%%%%%%%%%%%%%%%%%%%%%%%%%%%%

\section{Conclusion}
In this paper, we extended the Local Macroscopic Conservative (LoMaC) method to the Vlasov–Poisson system with the Dougherty–Fokker–Planck collision operator. The proposed method exploits the low-rank structures induced by Coulomb collisions to construct a low-rank approximation while locally preserving the mass, momentum, and total energy of the system at the discrete level. The method combines a discontinuous Galerkin discretization with a macroscopic conservative decomposition, allowing the physical accuracy of the solution to be maintained while removing redundant information through low-rank truncation.

Numerical experiments, including weak and strong Landau damping, two-stream instability, and bump-on-tail instability, demonstrate the accuracy and robustness of the proposed approach. In particular, the results show that collisional effects suppress the formation of filamentary structures and drive the plasma toward low-rank equilibria, substantially limiting rank growth and enabling efficient low-rank representations of the solution.

The framework can be naturally extended to higher-dimensional problems using hierarchical Tucker representations. Immediate future work will focus on developing a low-rank IMEX scheme for the VP–DFP system to better handle the stiffness introduced by the collision operator.

\section*{Data Availability} The numerical data and source code used to generate the results presented in this study are available upon reasonable request.

\section*{Competing Interests} The authors declare that they have no competing interests.

\section*{Acknowledgements} The authors thank the Pulse Tensor Network LDRD collaboration (Lawrence Livermore National Laboratory, Sandia National Laboratories, and Los Alamos National Laboratory) for supporting this research. The authors also thank I. Joseph for insightful discussions and helpful comments on the manuscript.

\printbibliography

\end{document}